\documentclass[11pt,a4paper,leqno]{article}
\usepackage{amsmath, amssymb, amsthm, amsfonts}
\usepackage{mathrsfs}

\theoremstyle{plain}
\newtheorem{thm}{Theorem}
\newtheorem{cor}{Corollary}
\newtheorem{lem}{Lemma}

\theoremstyle{remark}
\newtheorem{rem}{Remark}

\renewcommand{\Re}{{\rm Re\,}}
\renewcommand{\Im}{{\rm Im\,}}
\renewcommand{\a}{\alpha}
\renewcommand{\b}{\beta}

\newcommand{\Dsum}{\mathop{{\sum_{l=0}^J \sum_{m=0}^J}}}
\DeclareMathOperator*{\Res}{Res}

\allowdisplaybreaks

\numberwithin{equation}{section}

\begin{document}

\title{On the Tong-type identity and the mean square of the error term for an extended Selberg class }
\author{Xiaodong Cao,  Yoshio Tanigawa and Wenguang Zhai}
\date{}
\maketitle

\footnote[0]{This work is  supported by the National Key Basic Research Program of China(Grant No. 2013CB834201),
the National Natural Science Foundation of China(Grant No. 11171344),  the Natural
Science Foundation of Beijing(Grant No. 1112010) and the Fundamental Research Funds for the
Central Universities in China(2012YS01).}

\begin{abstract}
In 1956, Tong established an asymptotic formula for the mean square of the error term in the summatory function of the
Piltz divisor function $d_3(n).$  The aim of this  paper is to generalize Tong's method to a class of Dirichlet series that satisfy a
functional equation.   As an application, we can establish the asymptotic formulas for the mean square of the error terms
 for a class of functions in the well-known Selberg class. The  Tong-type identity and formula established in this paper can be viewed as
 an analogue of the well-known  Vorono\"{\i}'s formula.

\end{abstract}

\tableofcontents


\section{Introduction and main results}


Let $a(n)(n\geq 1)$ be a sequence of complex numbers, which is an arithmetic function. One of the most basic goals of analytic number theory is to
 establish   the  asymptotic formula  for the summatory function of $a(n),$  as accurate as possible.
Especially it is important to study the properties of the so-called error term  of this asymptotic formula,
such as the upper bound, the moments, the sign changes and
$\Omega$-results, etc.

 One of the important tools in this
area is Vorono\"{\i}'s formula of the error term when the sequence $a(n)(n\geq 1)$ satisfy
good conditions, for example, they are the coefficients of $L$-functions of any degree.
When $a(n) \ (n\geq 1)$ are the  coefficients of an $L$-function of degree two, Vorono\"{\i}'s formula of the corresponding
error term is a very strong tool to study the properties of  the error term. A well-known
 example is the Dirichlet divisor problem. For a survey of Vorono\"{\i}'s formula, see for example, Steuding
 \cite{St2}.

However, when $a(n) \ (n\geq 1)$ are the  coefficients of an $L$-function of degree $\geq 3$, Vorono\"{\i}'s formula is not  strong enough
to get good results for the properties of  the error term. Especially, it fails to give an  asymptotic formula
for the mean square of the error term even when the degree is $3.$

In 1956, Tong \cite{To2} established an asymptotic formula for the mean square of the error term in the summatory function of the
Piltz divisor function $d_3(n).$ This is the first result in this area for the case of degree $3$. Tong's main ingredient  is to replace the
error term by a corresponding integral such that the difference between
 the error term and this integral is very small on average.

The aim of this  paper is to generalize Tong's method to a class of Dirichlet series that satisfy a
functional equation.  Especially as an application, we   establish the asymptotic formulas for the mean square of the error terms
 for a class of functions in the well-known Selberg class.

\subsection{The Selberg class }

To study a summatory function of $a(n)$ and its error term, we consider the Dirichlet series with coefficients $a(n)$ which satisfies
a general functional equation.
There are many results in this direction, see for example, Chandrasekharan and Narasimhan \cite{CN1}, Redmond \cite{Red},
Hafner \cite{Ha1,Ha2}, Ivi\'c \cite{I2}, Meurman \cite{Me}, Lau \cite{La}, Kanemitsu, Sankaranarayanan and Tanigawa \cite{KST},
Friedlander and Iwaniec \cite{FI}. We shall study such a problem for $L$-functions in the so-called Selberg class.

The well-known Selberg class ${\cal S}$ (see for example \cite{Kac, Sel, St}) consists of non-vanishing Dirichlet series
$$
{\cal L}(s):=\sum_{n=1}^\infty\frac{a(n)}{n^s},
$$
which satisfies the following hypotheses:

\noindent
I. {\bf Ramanujan's conjecture}: $a(n)\ll n^\varepsilon$ for any $\varepsilon>0$.  \\[1ex]
II. {\bf Analytic continuation}: There exists a non-negative integer $m_{\cal L}$ such
that $(s-1)^{m_{\cal L}}{\cal L}(s)$ is an entire function of finite order. \\[1ex]
III. {\bf Functional equation}: ${\cal L}(s)$ satisfies a functional equation of type
\begin{equation} \label{functional-equation}
\Lambda_{{\cal L}}(s)=\omega\overline{\Lambda_{{\cal L}}(1-\bar{s})},
\end{equation}
where
\begin{equation} \label{gamma-factors}
\Lambda_{\cal L}(s):={\cal L}(s)Q^s\prod_{j=1}^{L}\Gamma(\a_j s+\b_j),
\end{equation}
and $Q>0, |\omega|=1$ and $\a_j>0, \b_j \in \mathbb{C}$ with $\Re \b_j\ge 0$ for all $1 \leq j \leq L$.
The number $d=2\sum_{j}\a_j$ is called the degree of ${\cal L}(s).$ \\[1ex]
IV. {\bf Euler product}: ${\cal L}(s)$ satisfies
$$
{\cal L}(s)=\prod_p \exp\left(\sum_{n\ge 1}\frac{b(p^n)}{p^{ns}}\right)
$$
with suitable coefficients $b(p^n)$ satisfying $b(p^n)\ll p^{nc}$ for some $c<1/2.$

Many well-known functions are contained in the Selberg class ${\cal S}$. We recall some examples.
The well-known Riemann zeta-function $\zeta(s)$ and Dirichlet $L$-functions are functions in ${\cal S}$ of
degree $1.$ The product of any $\ell$ functions in ${\cal S}$ of degree $1$ is a function in ${\cal S}$ of degree $\ell.$
So $\zeta^2(s)$ is a function in ${\cal S}$ of degree $2$ and $\zeta^3(s)$ is of degree $3.$
Let $f(z)$ be a holomorphic cusp form with respect to $SL_2({\Bbb Z})$. If $f(z)$ is an eigenform of all Hecke operators,
the automorphic $L$-function attached to $f(z)$ is a function in ${\cal S}$ of degree $2$.
The Dedekind zeta-function over an algebraic number field of degree $\kappa\ge 2$ is a function in ${\cal S}$ of degree $\kappa.$

The extended Selberg class ${\cal S}^{\#}$ (see  \cite{Kac, KP1, KP2} for an introduction)
 consists of all Dirichlet series $\sum_{n\ge 1}a(n)n^{-s}$ which satisfy  the conditions  $\mathrm{I}^{*}$, II and III, where
 $\mathrm{I}^{*}$ means

\medskip

 \noindent $\mathrm{I}^{*}$: $\sum_{n\ge 1}a(n)n^{-s}$ is absolutely convergent for $\sigma>1.$

\medskip

Consider the sum
$$
A(y):={\sum_{n\leq y\,}}^{\prime} a(n).
$$
Define $Q(y)$ as the sum of the residues of the function ${\cal L}(s)y^s s^{-1}$ and the error term
$E(y) $  is defined by
$$E(y):=A(y)-Q(y).$$

Estimates for $E(x)$ and $\int_0^T|E(y)|^2dy$ were first studied by Vorono\"i in 1904 and Cram\'er in 1922
for the special case of the Dirichlet  divisor problem respectively, see \cite{V} and \cite{Cr}.  Since then this has been
generalized for larger classes of Dirichlet series.
There are many results in this direction, see for example, \cite{CN4, CN5, Ha1, Ha2, La, Red, Ro}.

Roton \cite{Ro} proved that if $L\in {\cal S}^{\#}$ is a function of degree $d\geq 2.$ Then
\begin{equation*}E(y)\ll y^{(d-1)/(d+1)+\varepsilon},\end{equation*}
which is a  generalization of Landau's classical result \cite{Lan}.
She also proved that if $$\sum_{n\leq y}|a(n)|^2\ll y^{1+\varepsilon},$$
 then
\begin{eqnarray*}
\int_0^T E^2(y)dy\ll \left\{\begin{array}{ll}
T^{2-1/d},&\mbox{if $0<d<3,$}\\
T^{3-4/d+\varepsilon},& \mbox{if $d\ge 3$.}
\end{array}\right.
\end{eqnarray*}


\subsection{Statements of main results}


Suppose $0\le \theta <1$  is a real number. Let ${\cal S}^{\theta}$ denote the set of all Dirichlet series
${\cal L}(s)=\sum_{n\ge 1}a(n)n^{-s}$ which satisfy II, III and
\begin{equation}\label{meanC}
{\mathrm I}^{\prime}:\ \ \ |a(n)|\ll n^{\theta+\varepsilon},\ \ \ \sum_{n\leq y}|a(n)|^2\ll y^{1+\varepsilon}.
\end{equation}
In this paper we assume that $\sum_{j=1}^{L} \b_j$ is real for the sake of simplicity.
Obviously ${\cal S}\subset  {\cal S}^{\theta} \subset {\cal S}^{\#}.$ Let $  {\cal S}_{real}^{\theta} $
denote the set of all functions in ${\cal S}^{\theta}$ with real coefficients.

Now suppose  $a(n)(n\geq 1)$  is a sequence of real numbers such that
its corresponding Dirichlet series ${\cal L}(s)\in {\cal S}_{real}^{\theta} $  and is of degree $d\ge 2.$
Without loss of generality, we suppose that $d$ is an integer.
Let $1/2\le \sigma^{*}<1 $ denote a fixed real number such that the estimate
\begin{equation}\label{meanL}
 \int_0^T|{\cal L}(\sigma^{*}+it)|^2dt\ll T^{1+\varepsilon}
\end{equation}
holds for any $\varepsilon>0.$ We suppose that $\sigma^{*}$ satisfies the condition
\begin{equation}  \label{sigmastarcondition}
\sigma^{*}<(d+1)/2d,
\end{equation}
which plays an important role in our paper.


\begin{thm} Suppose that $d\ge 2$ is a fixed integer, $0\le \theta\le 1/2-1/2d$ is a real number.
Suppose that ${\cal L}(s)\in {\cal S}_{real}^{\theta}$ is a function of degree
 $d\ge 2$ such that   \eqref{meanL} and \eqref{sigmastarcondition} hold. Then we have
\begin{equation}
\int_1^T E^2(y)dy=C_dT^{2-1/d}+O(T^{2-\frac{3-4\sigma^{*}}{2d(1-\sigma^{*})-1}+\varepsilon}),
\end{equation}
where $C_d>0$ is a positive constant.
\end{thm}


\begin{cor} \label{cor-1}
Suppose $0\le \theta\le 1/4$ is a real number and  ${\cal L}(s)\in {\cal S}_{real}^{\theta}$ is of degree $2,$ then we have
\begin{equation} \label{meanE}
\int_1^T E^2(y)dy=C_2T^{3/2}+O_\varepsilon(T^{1+\varepsilon}).
\end{equation}
\end{cor}

\begin{rem}
When ${\cal L}(s)\in {\cal S}_{real}^{0}$ is of degree 2, the result obtained by Vorono\"i's formula
is sometimes stronger. See for example, Meurman \cite {Me} or Lau and Tsang \cite{LT2} in the case of Dirichlet divisor problem.
We note that the error term in the asymptotic formula \eqref{meanE} is best possible when disregarding $\varepsilon.$
\end{rem}


\begin{cor} \label{cor-2}
Suppose $0\le \theta\le 1/3$ is a real number.
Let ${\cal L}(s)\in {\cal S}_{real}^{\theta}$ is a function of degree $3$ such that
${\cal L}(s)={\cal L}_1(s){\cal L}_2(s),$ where ${\cal L}_1(s)\in {\cal S}_{real}^{\theta}$
is a function of degree $1$, and ${\cal L}_2(s)\in {\cal S}_{real}^{\theta}$ is a function of degree $2$.
Then we have
\begin{equation} \label{unconditional}
\int_1^T E^2(y)dy=C_3T^{5/3}+O_\varepsilon(T^{8/5+\varepsilon}).
\end{equation}
Furthermore if we assume that
\begin{equation} \label{assump-d2}
\int_0^T|\mathcal{L}_2(1/2+it)|^6dt \ll T^{2+\varepsilon},
\end{equation}
then we have
\begin{equation}  \label{conditional}
\int_1^T E^2(y)dy=C_3T^{5/3}+O_\varepsilon(T^{14/9+\varepsilon}).
\end{equation}
\end{cor}

\begin{rem}
The assumption \eqref{assump-d2} for degree $2$ $L$-functions  is very natural. The first result in this direction is  the twelfth power moment of the
Riemann zeta-function $\zeta(s)$ (the sixth moment of $\zeta^2(s)$ over the critical line) due to Heath-Brown \cite{He2}.  Meurmann \cite{Me3} proved
the twelfth power moment for any Dirichlet L-functions. From Heath-Brown \cite{He2} and Meurmann \cite{Me3} we see that \eqref{assump-d2} holds for
  Dedekind zeta functions of   quadratic fields. From \cite{J2} we know that \eqref{assump-d2} holds for
 the automorphic $L$-functions attached to holomorphic cusp forms.
\end{rem}

\begin{rem}
Tong \cite{To2} first established Theorem 1 for the Piltz divisor problem of dimension $d\ge 3,$
which becomes a true asymptotic formula when $d=3.$ Fomenko \cite{Fo1} followed Tong's method to study the case of the
Dedekind zeta-function of cubic fields.
\end{rem}


\begin{thm} Suppose that $d\ge 2$ is an integer, $0\le \theta\le 1/2-1/2d $ is a real number.
Suppose that ${\cal L}(s)\in {\cal S}_{real}^{\theta}$ is a function of degree
 $d$ such that  \eqref{meanL} and \eqref{sigmastarcondition} hold. Then the error term $E(t)$ has a distribution
function $f(\alpha)$   in the sense that, for any interval
$I\subset {\Bbb R}$ we have
$$T^{-1} mes \{t\in [1,T]: t^{-(d-1)/2d)}E(t)\in I\}\rightarrow\int_If(\alpha)d\alpha$$
as $T\rightarrow \infty.$ The function $f(\alpha)$ and its
derivatives satisfy
$$\frac{d^k}{d\alpha^k}f(\alpha)\ll_{A,k}(1+|\alpha|)^{-A} $$
for $k=0,1,2,\cdots$ and  $f(\alpha)$ can be extended to an entire
function.
\end{thm}

\begin{cor} 
Suppose   $0\le \theta\le 1/2-1/2d $ is a real number.
Suppose that ${\cal L}(s)\in {\cal S}_{real}^{\theta}$ is a function of degree
$2,$ or is a function of degree $3$ such that    it can be written as a product of
a function of degree $1$ and a function of degree $2$,  then Theorem 2 holds.
\end{cor}


\begin{thm} Suppose that $d\ge 2$ is an integer, $0\le \theta\le 1/2-1/2d $ is a real number.
Suppose that ${\cal L}(s)\in {\cal S}_{real}^{\theta}$ is a function of degree
 $d$ such that   \eqref{meanL} and \eqref{sigmastarcondition} hold. Then for any real number $0\leq u\leq 2,$  the mean
value
$$\lim_{T\rightarrow\infty}T^{-1-\frac{(d-1)u}{2d}}\int_1^T |E(t)|^udt$$ converges to
a finite limit as $T$ tends to infinity.
\end{thm}

\begin{cor} 
Suppose   $0\le \theta\le 1/2-1/2d $ is a real number.
Suppose that ${\cal L}(s)\in {\cal S}_{real}^{\theta}$ is a function of degree
 $2,$ or is a function of degree $3$ such that    it can be written as a product of
 a function of degree $1$ and a function of degree $2$,  Then Theorem 3 holds.
\end{cor}


\subsection{A remark on Tong's method}


The most important tool up to now to study the behaviour of $E(y)$ is Vorono\"i's formula, which is usually of the form
\begin{equation} \label{truncated-V}
E(y)=c_1y^{\frac{d-1}{2d}}\sum_{n\leq N}\frac{a(n)}{n^{\frac{d+1}{2d}}}\cos(c_2(ny)^{1/d}+c_3)+O(y^{\frac{d-1}{d}+\varepsilon}N^{-\frac 1d}),
\end{equation}
where $c_1, c_2, c_3$ are constants and $1\ll N\ll y.$ For a proof of the above general formula, see for example
Friedlander and Iwaniec \cite{FI}.
Roughly speaking, this formula gives a good approximation to the error term, which is usually
a finite exponential sum, plus a ``permissible small" error.

When $d=2,$ the formula \eqref{truncated-V} is strong enough for us to study the properties of $E(x).$ For example,
 its upper bound, power moments, sign changes, etc.
 A good example is the  Dirichlet divisor problem, see for example, the survey paper Tsang \cite{TS}.

 However when $d\geq 3, $ the error term in \eqref{truncated-V} is always $\gg y^{1-2/d+\varepsilon}$, which is much larger than the expected order
 $y^{(d-1)/2d}$.  So \eqref{truncated-V} can't be used to study the asymptotic behaviour of power moments of $E(y),$ even for the mean square.

In \cite{To1, To2}, Tong developed a method to study the mean square of the error term in the summatory function of the
Piltz divisor function $d_{\ell}(n),$ which denotes the number of ways such that $n$ can be written as a product of $\ell$
natural numbers. When $\ell=3,$ Tong's method gives a true asymptotic formula of the mean square of the error term $\Delta_3(y)$.
Tong's main ingredient  is to replace the error term by a corresponding integral such that the difference between
 the error term and this integral is very small on average.

The aim of this paper is to generalize Tong's method to the general case. Actually we shall generalize Tong's method to
 a class of functions much more general than ${\cal S}_{real}^{\theta}$.    As applications, we    establish
 the asymptotic formula of the mean square of $E(y)$ for   functions in ${\cal S}_{real}^{\theta}$.


\subsection{Organization of this paper}


The organization of this paper is as follows. In Section 2 we shall introduce  a class of more general
functions. In Section 3 we give some preliminary lemmas about the so-called generalized Bessel functions.
 In Section 4 we shall establish a Tong-type identity for the
corresponding error term $E(y)$. As  some simple applications of the Tong's identity,
in Section 5 we give a lower bound for integrals involving $E(y)$ and study the small values of
$E(y).$ In Section 6 we shall establish a truncated Tong-type formula for an integral involving   $E(y).$
In Section 7 we shall estimate some exponential integrals, which are important in the proof.
In Section 8 we shall prove Theorem 1. In Section 9 we shall prove Theorem 2 and Theorem 3.
In Section 10, we give some examples.


\section{A class of more general arithmetic functions}

For future applications in mind (e.g.~\cite{CTZ}), we shall derive a Tong-type identity of the error term for a more general
class of functions than that of Selberg class. In this section, following basically to Chandrasekharan and Narasimhan \cite{CN1,CN4}
and Hafner \cite{Ha1}, we introduce such a class of functions.

\smallskip

\noindent {\bf Definition. } Let $\{a(n)\}$ and $\{b(n)\}$ be two sequences of complex numbers, not identically zero.
Let $\{\lambda_n\}$ and $\{\mu_n\}$ be two strictly increasing sequences of positive numbers tending to infinity. Suppose that the series
\begin{align*}
\varphi(s)=\sum_{n=1}^{\infty}a(n)\lambda_n^{-s},\quad
\psi(s)=\sum_{n=1}^{\infty}b(n)\mu_n^{-s}
\end{align*}
converge in some half-plane and have abscissae of absolute convergence $\sigma_a^*$ and $\sigma_b^*$, respectively.
Define two gamma factors
\begin{align}
\Delta_1(s)=\prod_{j=1}^{N}\Gamma(\alpha_j s+\beta_j) \label{gamma-1}
\intertext{and}
\Delta_2(s)=\prod_{h=1}^{N'}\Gamma(\alpha_h' s+\beta_h') \label{gamma-2},
\end{align}
where $\alpha_j$ and $\alpha_h'$ are positive real numbers and $\beta_j$ and $\beta_h'$ are complex numbers.
We assume throughout that
\begin{equation}  \label{def-degree}
\a:=\sum_{j=1}^N \alpha_j = \sum_{h=1}^{N'} \alpha_h'.
\end{equation}

Let $r$ be real. We say that $\varphi$ and $\psi$ satisfy the functional equation
\begin{align}
\Delta_1(s)\varphi(s)=\Delta_2(r-s)\psi(r-s) \label{feq}
\end{align}
if there exists in the $s$-plane a domain $D$ that is the exterior of a compact set $\mathfrak{S} $ (we call it a ``singularity set")
 and on which there exists a holomorphic function $\chi(s)$ ($s=\sigma+it$, $\sigma$ and $t$ real) such that
\begin{align} \texttt{(i)}\quad \lim_{|t|\rightarrow\infty}\chi(\sigma+it)=0\nonumber
\end{align}
uniformly in every interval $-\infty<\sigma_1\le \sigma\le \sigma_2<+\infty$, and
\begin{align}
\texttt{(ii)}\quad \chi(s)=\begin{cases}
\Delta_1(s)\varphi(s) & \mbox{for $\sigma>\sigma_a^*$,}\\
\Delta_2(r-s)\psi(r-s) & \mbox{for $\sigma<r-\sigma_b^*$.}
\end{cases}.\nonumber
\end{align}
Finally suppose that both $\varphi(s)$ and $\psi(s)$ have only a finite number of poles on the complex plane ${\Bbb C}.$

\begin{rem} When $\Delta_1(s)=\Delta_2(s),$ we get the class of functions defined in Chandrasekharan and Narasimhan \cite{CN1, CN4}.
 Clearly any function in the extended Selberg class ${\cal S}^{\#}$ satisfies the functional equation of the form \eqref{feq}.
\end{rem}

\begin{rem} There are many other arithmetic functions in the above class but $\Delta_1(s)\not=\Delta_2(s).$ One example is the well-known
asymmetric many-dimensional divisor problem. For a survey of the asymmetric many-dimensional divisor problem, see Kr\"atzel \cite{Kra} or
Ivi\'c etc. \cite{IKN}.
A forthcoming paper \cite{CTZ} in this direction is in preparation through the approach of this paper.
\end{rem}

From now on we always assume that $a(n), b(n)$ are real. We suppose that there are infinitely many $n$ such
that $|a(n)|\gg 1$ and the same holds for $b(n).$

For $y>0$ and a real number $\varrho $, we define the summatory  function $A_\varrho(y)$ of the arithmetical function $a(n)$ by
\begin{align*} 
A_\varrho(y)=\frac {1}{\Gamma(\varrho+1)}{\sum_{\lambda_n\le y}}^{\prime }a(n)(y-\lambda_n)^\varrho,
\end{align*}
where the symbol $\prime$ indicates that the last term has to be halved if $\varrho=0$ and $y=\lambda_n$.
When $\varrho$ is negative, $A_\varrho(x)$ is defined only for those positive $y$ not equal to any $\lambda_n$.

Let $s_0=\sup\{|s|:s\in \mathfrak{S}  \},$ where $\mathfrak{S}  $ is the "singularity set" in Definition,
and $t_0=\max\{|\frac{\beta_j}{\alpha_j}|, |\frac{\beta_h'}{\alpha_h'}| : j=1,2,\ldots,N, h=1,2,\ldots, N' \}$.
Choose two constants $c>\max\{\sigma_a^*,\sigma_b^*,s_0,t_0\}$ and $R>\max\{s_0,t_0\}$.
Choose the third constant $b\not\in \mathbb{Q}$ such that $r-b$ is not an integer and $b>\max\{c,r\}$.
We choose $a \leq \min\{\sigma_b^*-\frac {2}{\alpha},\frac{r}{2}-\frac{1}{2\alpha}\}$. In fact $a$ will
be chosen to be small so that the relevant integral converges absolutely.

Let $\mathcal{C}$ be the rectangle with vertices $c\pm iR$ and $ r-b\pm iR$,
taken in the counter-clockwise direction. For convenience, we shall use $\mathcal{C}_{u,v}$ to denote
the oriented polygonal path with vertices $u-i\infty, u-iR$, $v-iR, v+iR, u+iR$, and $ u+i\infty$ in this order.

We define the residual function or the main term
\begin{align} \label{eq2-7}
Q_\varrho (y)=\frac {1}{2\pi i}\int_{\mathcal{C}}
\frac{\Gamma (s)\varphi(s)y^{\varrho+s}}{\Gamma(s+\varrho+1)}ds.
\end{align}

From our choices of $b,c$ and $R$, the path $\mathcal{C}$ encircles all of $S$. By the residue theorem we have
\begin{equation} \label{eq2-8}
Q_\varrho(y)=\sum_{s_j}y^{\varrho+s_j}P_{s_j}(\log y),
\end{equation}
where $s_j$ runs over all the poles of $\Gamma(s)\varphi(s)y^{\varrho+s}/\Gamma(s+\varrho+1)$ inside $\mathcal{C}$,
$P_{s_j}(t)$ is a polynomial of $t$ such that its degree is  the order of the pole $s_j$ minus $1$.

We define the error term in the asymptotic formula for $A_\varrho(y)$ as
\begin{equation*}
E_\varrho(y)=A_\varrho(y)-Q_\varrho(y).
\end{equation*}
As for $\psi(s)$ we define similarly that
\begin{align*}
A_\varrho^{*}(y)&=\frac {1}{\Gamma(\varrho+1)}{\sum_{\mu_n\le y}}^{\prime }b(n)(y-\mu_n)^\varrho,\\
Q_\varrho^{*} (y)&=\frac {1}{2\pi i}\int_{\mathcal{C}} \frac{\Gamma (s)\psi(s)y^{\varrho+s}}{\Gamma(s+\varrho+1)}ds,\\[1ex]
E_\varrho^{*}(y)&=A_\varrho^{*}(y)- Q_\varrho^{*} (y).
\end{align*}
Similarly to \eqref{eq2-8} we have
 \begin{equation*}
Q^{*}_\varrho(y)=\sum_{s_j}y^{\varrho+s_j}P_{s_j}^{*}(\log y),
\end{equation*}
where $s_j$ runs over all the poles of $\Gamma(s)\psi(s)y^{\varrho+s}/\Gamma(s+\varrho+1)$ inside $\mathcal{C}$,
$P_{s_j}^{*}(t)$ is a polynomial such that its degree is  the order of the pole $s_j $ minus $1$.

When $\varrho=0,$ for simplicity, we also use the following notation
\begin{align*}
A(y)&=A_0(y), & Q(y)&=Q_0(y), & E(y)&=E_0(y),\\
A^{*}(y)&=A_0^{*}(y),& Q^{*}(y)&=Q_0^{*}(y), & E^{*}(y)&=E_0^{*}(y).
 \end{align*}

We suppose that $\sigma^*<\sigma_b^{*}$ is a real number such that the estimate
\begin{align}\label{psi-condition}
\int_{-T}^{T}|\psi(\sigma^*+it)|^2 dt\ll T^{1+\varepsilon}
\end{align}
holds. Since there are infinitely many $n$ such that $b(n)\gg 1,$ we have $\sigma^{*}\geq 0.$
We also suppose that $\psi(s)$ doesn't have any poles in the region $\sigma\leq \sigma^{*}.$


\section{Some preliminary lemmas}


We shall prove a fundamental lemma of the asymptotic expansion of the integral which we need later.
First we recall the well-known Stirling's formula for gamma function.

\begin{lem}[Stirling's formula]
Let $c$ be a constant. Then there exist some constants $\gamma_j$ and $\gamma_j'$ depending on $c$ such that for any positive integer $m$,
\begin{align}
\log \Gamma(s+c) = \left(s+c-\frac12\right) \log s -s+\frac12 \log 2\pi +\sum_{j=1}^m \frac{\gamma_j}{s^j}
+O\left(\frac{1}{|s|^{m+1}}\right) \label{log-Stirling}
\end{align}
and
\begin{align}
\Gamma(s+c) = \sqrt{2\pi}e^{\left(s+c-\frac12\right)\log s-s}
\left(1+\sum_{j=1}^m \frac{\gamma_j'}{s^j}+O\left(\frac{1}{|s|^{m+1}}\right)\right)  \label{Stirling}
\end{align}
uniformly for $|\arg s|<\pi-\delta$ for fixed $\delta>0$, as $|s| \to \infty$.
\end{lem}
The constants $\gamma_j$ in \eqref{log-Stirling} are given explicitly by
$$
\gamma_j=\frac{(-1)^{j-1}B_{j+1}(c)}{j(j+1)},
$$
where $B_j(x)$ is the Bernoulli polynomials of degree $j$.
See for example,  Wang and Guo \cite[p.~123 (5)]{WG}  and see also A. Erd\'elyi
et al. \cite{E}.

\medskip

For two gamma factors $\Delta_1(s)$ and $\Delta_2(s)$ defined by \eqref{gamma-1} and \eqref{gamma-2}, respectively, we define
\begin{align}
\mu&=\sum_{j=1}^N\left(\b_j-\frac12\right)+\frac12, & \mu'&=\sum_{h=1}^{N'}\left(\b_h'-\frac12\right)+\frac12,  \label{def-mu} \\
\nu&=\sum_{j=1}^N\left(\b_j-\frac12\right)\log \a_j, & \nu'&=\sum_{h=1}^{N'}\left(\b_h'-\frac12\right)\log \a_h', \label{def-nu} \\
\tau&=\sum_{j=1}^N\a_j \log \a_j, & \tau'&=\sum_{h=1}^{N'}\a_h' \log \a_h'. \label{def-lambda}
\end{align}
Furthermore we define
\begin{equation} \label{def-thetarho}
\theta_\varrho=\frac{r}{2}-\frac{1}{4\a}+\varrho\left(1-\frac{1}{2\a}\right)+\frac{\mu'-\mu}{2\a}
\end{equation}
and
\begin{equation} 
h=2\a e^{-\frac{\tau+\tau'}{2\a}}. \nonumber
\end{equation}

Let $\mathcal{C}_{a,b}$ be a curve defined in the previous section such that all poles of $\Delta_2(s)$ lie in the left hand side of
$\mathcal{C}_{a,b}$. If $\b_j$ and $\b_h'$ are real, this condition means that $b$ is greater than
the maximal pole of $\Delta_2(s)$.

\begin{lem}
We assume that
\begin{equation} \label{gamma-katei}
\overline{\Delta_1(s)}=\Delta_1(\bar{s}) \quad \mbox{and} \quad \overline{\Delta_2(s)}=\Delta_2(\bar{s}).
\end{equation}
Let $\omega$ be a real number and $M$ an integer. Let $a$ be a real number such that
\begin{equation} \label{a-jyouken}
a<\frac{r}{2}-\frac{\mu'-\mu+\omega+M+1}{2\a}.
\end{equation}
Let $\mathcal{D}$ be a domain such that
$$
\mathcal{D}=\mathbb{C} \backslash \{s \in \mathbb{C}\,|\, \Re s < b', |\Im s|<R'\}
$$
with some constant $b'$ and $R'$. We suppose that $\mathcal{D}\supset \mathcal{C}_{a,b}$, and furthermore if $\omega \neq 0$
we assume that $b'>0$.
Suppose that a function $F(s)$ is regular in the domain $\mathcal{D}$ and
has an asymptotic expansion
\begin{equation}  \label{def_function_F}
F(s) \sim s^M \sum_{n=0}^{\infty}\frac{d_n}{s^n}  \ \ \ (d_0 \neq 0)
\end{equation}
with real coefficients $d_n$ as $|s| \to \infty$. Let $\mathcal{I}(x)$ be the function defined by
\begin{align}  \label{main-integral}
\mathcal{I}(x):=\int_{\mathcal{C}_{a,b}}\frac{\Delta_2(s)}{\Delta_1(r-s)}s^{\omega}F(s)x^{-s}ds.
\end{align}
Then for any positive integer $m$ we have an asymptotic expansion
\begin{align*}
\mathcal{I}(x)&=2i C \sum_{k=0}^{m} c_k x^{\frac{\tilde{M}-k+\frac12}{2\a}}\cos\left(hx^{\frac{1}{2\a}}-a_k\pi\right)
+O\left(x^{\frac{\tilde{M}-m-\frac12}{2\a}}\right),
\end{align*}
where
\begin{align}
\tilde{M}&=\mu'-\mu-\a r+\omega+M, \label{Mtilde} \\
a_k&=\frac{1}{2}\left(M+\omega-k+\a r+\mu+\mu'-\frac12\right), \nonumber \\
C  &=(2\pi)^{\frac{N'-N}{2}}e^{\nu'-\nu-\tau r-\frac{(\tau+\tau')(\tilde{M}+\frac12)}{2\a}} \nonumber
\end{align}
and $c_k$ are real constants which do not depend on $x$.
In particular
$$
c_0=\frac{\sqrt{\pi}}{\sqrt{\a}}d_{0}
$$
and
$$
c_1=\frac{\sqrt{\pi}}{\sqrt{\a}}\left\{d_{0}\left(-\frac{1}{2\a}\left(\frac{\tilde{M}^2}{2}-\frac{1}{24}\right)
+\frac{1}{2}\left(\sum_{h=1}^{N'}\frac{B_2(\b_h')}{\a_h^{\prime}}+\sum_{j=1}^N\frac{B_2(\alpha_j r+\beta_j)}{\a_j}\right)\right)+d_{1}\right\},
$$
where $B_2(x)$ is the Bernoulli polynomial of degree 2.
\end{lem}

\begin{rem}
This lemma is a generalization of Theorem 3 in Tong \cite{To1}. But Tong omitted the  proof of his Theorem 3.
\end{rem}

\proof  We note first that the integral \eqref{main-integral} is  convergent absolutely under the condition \eqref{a-jyouken}.

We derive an asymptotic expansion of the integrand of \eqref{main-integral}. Here for simplicity we use the symbol $\sim$ to denote that
the right hand side of $\sim$ is an asymptotic expansion of a function in the left hand side.

Consider $\Delta_2(s)$ first.
By \eqref{Stirling}, there exist constants $b_h'$ such that
\begin{align*}
\Delta_2(s) = \prod_{h=1}^{N'}\Gamma(\a_h's+\b_h') \sim (2\pi)^{N'/2} \exp(g_2(s))\left(1+\frac{b_1'}{s}+\frac{b_2'}{s^2}+\cdots\right), \label{delta_2-1}
\end{align*}
where $g_2(s)$ is given by
\begin{align}
g_2(s)&=\sum_{h=1}^{N'}\left\{\left(\a_h's+\b_h'-\frac12\right)(\log s+\log \a_h')-\a_h's\right\} \notag \\
&=\left(\a s+\mu'-\frac12\right)\log s+(\tau'-\a)s+\nu'. \nonumber 
\end{align}

Similarly we have, with some constants $b_j''$,
\begin{align*}
\Delta_1(r-s)&=\prod_j\Gamma(-\a_js+\a_jr+\b_j) \\
&\sim (2\pi)^{N/2}\exp(g_1(s))\left(1+\frac{b_1''}{s}+\frac{b_2''}{s^2}+\cdots\right),
\end{align*}
where $g_1(s)$ is given by
\begin{align}
g_1(s)&=\sum_{j=1}^{N}\left\{\left(-\a_js+\a_jr+\b_j-\frac12\right)(\log(-s)+\log \a_j)+\alpha_j s \right\}  \notag \\
& =\left(-\a s+ \a r+\mu-\frac12\right)\log(-s)-(\tau-\a)s+\tau r+\nu. \nonumber
\end{align}
In particular
$$
b_1'=\frac12 \sum_{h=1}^{N'}\frac{B_2(\b_h')}{\a_h^{\prime}} \qquad \mbox{and} \qquad b_1''=-\frac12\sum_{j=1}^{N}\frac{B_2(\a_jr+\b_j)}{\a_j}.
$$
By using
$$
\log(-s)=\log s - \pi i \, {\rm sgn}(t)
$$
for  non-zero and non-negative $s$, where $t=\Im s$, we can see easily that
\begin{align}
g(s):&=g_2(s)-g_1(s) \notag \\
&=(2\a s+\mu'-\mu- \a r)\log s -\pi i \a s \, {\rm sgn}(t) \notag \\
& \quad +\pi i \left(\a r+\mu-\frac12\right){\rm sgn}(t) +(\tau+\tau'-2\a)s+\nu'-\nu-\tau r. \notag
\end{align}

Combing these formulas we find that with some constants $p_n$,
\begin{align*}
&\frac{\Delta_2(s)}{\Delta_1(r-s)}s^{\omega}F(s)x^{-s} \\
&\sim (2\pi)^{(N'-N)/2}\exp(g(s))e^{-s\log x}\left(1+\frac{b_1'}{s}+\cdots\right)\left(1+\frac{b_1''}{s}+\cdots \right)^{-1}\\
& \qquad \times  s^{\omega+M} \left(d_{0}+\frac{d_{1}}{s}+\cdots\right)\\
&=(2\pi)^{(N'-N)/2}\exp(g(s))e^{-s\log x}s^{\omega+M}\sum_{n=0}^{\infty}\frac{p_n}{s^{n}}.
\end{align*}


By \eqref{gamma-katei} and the assumption that $d_n$ are real, we have
$$
\overline{\int_{\mathcal{C}_{a,b}\cap \{s | t >0\}} \frac{\Delta_2(s)}{\Delta_1(r-s)}s^{\omega}F(s)x^{-s}ds}
=- \int_{\mathcal{C}_{a,b}\cap \{s | t <0\}} \frac{\Delta_2(s)}{\Delta_1(r-s)}s^{\omega}F(s)x^{-s}ds,
$$
hence in order to evaluate $\mathcal{I}(x)$ it is enough to consider the case $t>0$ and take the imaginary part.
So we suppose that $t>0$. Let
\begin{equation*} 
f(s)=2\a s\log s+(\tau+\tau'-2\a)s-\pi i \a s -s\log x.
\end{equation*}
Then
\begin{align}
\frac{\Delta_2(s)}{\Delta_1(r-s)}s^{\omega}F(s)x^{-s} \sim & C_0 e^{\pi i(\a r+\mu-1/2)} s^{\mu'-\mu-\a r+\omega+M}
\exp(f(s)) \sum_{n=0}^{\infty}\frac{p_n}{s^n}  \label{integrad}
\end{align}
where we put
$$
C_0=(2\pi)^{(N'-N)/2}e^{\nu'-\nu-\tau r}.
$$

Now let
$$
y:=xe^{-(\tau+\tau')}.
$$
We define the new parameters $w$, $\xi$ and $\eta$ by
\begin{align*}
&s=y^{\frac{1}{2\a}}w \\
&w=i(1+\xi), \ \ \mbox{$|\xi|$ small} \\
&\eta^2=-2\a i\left((1+\xi)\log(1+ \xi)-\xi\right).
\end{align*}
The branch of $\eta$ will be taken as
$$
\xi=\xi(\eta)=\frac{i^{1/2}\eta}{\sqrt{\a}}+\frac{1}{6}\left(\frac{i^{1/2}\eta}{\sqrt{\a}}\right)^2
-\frac{1}{72}\left(\frac{i^{1/2}\eta}{\sqrt{\a}}\right)^3 +\cdots.
$$

Let $\delta$ be a small constant. 
We put
$$
w_1=u_1+iv_1=i(1+\xi(-\delta)) \ \ w_2=u_2+iv_2=i(1+\xi(\delta)).
$$
Define the paths of integration by
\begin{align*}
\mathcal{L}_0&=\left\{s \ | \ s=y^{\frac{1}{2\a}}i(1+\xi(\eta)),  -\delta \leq \eta \leq \delta \right\} \\
\mathcal{L}_1&=\left\{s \ | \ s=y^{\frac{1}{2\a}}(u_1+iv), 0 \leq v \leq v_1 \right\} \\
\mathcal{L}_2&=\left\{s \ | \ s=y^{\frac{1}{2\a}}(u_2+iv), v \geq v_2 \right\}.
\end{align*}
By Cauchy's theorem and the remark above,
$$
\mathcal{I}(x)=2i \sum_{j=0}^2 \Im \int_{\mathcal{L}_j} \frac{\Delta_2(s)}{\Delta_1(r-s)}s^{\omega}F(s)x^{-s}ds.
$$

First we consider the integral
$$
J_0=\int_{\mathcal{L}_0} \frac{\Delta_2(s)}{\Delta_1(r-s)}s^{\omega}F(s)x^{-s}ds.
$$
From the above choice of parameters, we have
\begin{align*}
f(s)&=2\a s \log s -2\a s-\pi i \a s-s\log y \\
&=2\a y^{\frac{1}{2\a}}i(1+\xi)\Bigl(\log(1+\xi)-1\Bigr) \\
&=-2\a y^{\frac{1}{2\a}}i+2\a y^{\frac{1}{2\a}}i\Bigl((1+\xi)\log(1+\xi)-\xi\Bigr).
\end{align*}
Therefore we have
\begin{align}
\frac{\Delta_2(s)}{\Delta_1(r-s)}s^{\omega}F(s)x^{-s}
=C_0 e^{-2\a iy^{\frac{1}{2\a}}} e^{\pi i (\a r+\mu-1/2)} \exp(-y^{\frac{1}{2\a}}\eta^2)
 \sum_{n=0}^{\infty}p_n s^{\tilde{M}-n}, \notag 
\end{align}
where $\tilde{M}$ is defined by \eqref{Mtilde}. In terms of the parameter $\eta$, we have
\begin{align*}
&\sum_{n=0}^{\infty} p_n s^{\tilde{M}-n}=\sum_{n=0}^{\infty}p_n\left(iy^{\frac{1}{2\a}}(1+\xi)\right)^{\tilde{M}-n} \\
&=\sum_{n=0}^{\infty}p_n(iy^{\frac{1}{2\a}})^{\tilde{M}-n}
  \left(1+(\tilde{M}-n)\xi+\frac{(\tilde{M}-n)(\tilde{M}-n-1)}{2}\xi^2 +\cdots\right) \\
&=\sum_{n=0}^{\infty}p_n(iy^{\frac{1}{2\a}})^{\tilde{M}-n}
  \left(1+(\tilde{M}-n)\frac{i^{\frac12}\eta}{\sqrt{\a}}+\frac{(\tilde{M}-n)(\tilde{M}-n-\frac23)}{2}
\left(\frac{i^{\frac12}\eta}{\sqrt{\a}}\right)^2 +\cdots\right)
\end{align*}
and
$$
ds=iy^{\frac{1}{2\a}}d\xi=iy^{\frac{1}{2\a}}\frac{i^{\frac12}}{\sqrt{\a}}\left(1+\frac{i^{\frac12}\eta}{3\sqrt{\a}}-\frac{i\eta^2}{24\a}+\cdots\right)d\eta,
$$
hence
\begin{align}
\left(\sum_{n=0}^{\infty}p_n s^{\tilde{M}-n}\right)ds
=\sum_{n=0}^{\infty}p_n(iy^{\frac{1}{2\a}})^{\tilde{M}-n+1} i^{\frac12}\sum_{j=0}^{\infty}q_{n,j}(i^{\frac12}\eta)^jd\eta \notag 
\end{align}
with some real constants $q_{n,j}$.

Now take the integral of $\eta$ over $[-\delta, \delta]$. The terms of odd powers of $\eta$ vanish.
By using the well known formula
$$
\int_{-\delta}^{\delta}e^{-a\eta^2}\eta^{2h}d\eta=a^{-h-1/2}\Gamma(h+1/2)+O_{\delta, h}\left(\frac{e^{-a\delta^2}}{a}\right) \ \ a>0,
$$
we find that
\begin{align*}
&\int_{\mathcal{L}_0}\left(\exp\left(-y^{\frac{1}{2\a}}\eta^2\right)\sum_{n=0}^{\infty}p_n s^{\tilde{M}-n}\right)ds \\
&\sim \sum_{n=0}^{\infty}\sum_{h=0}^{\infty}p_n q_{n,2h}\left(iy^{\frac{1}{2\a}}\right)^{\tilde{M}-n+1}i^{\frac12+h}
\int_{-\delta}^{\delta}e^{-y^{\frac{1}{2\a}}\eta^2}\eta^{2h}d\eta\\
&=\sum_{n=0}^{\infty}\sum_{h=0}^{\infty}p_n q_{n,2h}\left(iy^{\frac{1}{2\a}}\right)^{\tilde{M}-n+1}i^{\frac12+h}
\left(y^{-\frac{h+\frac12}{2\a}}\Gamma(h+1/2)+O\left(\frac{e^{-y^{\frac{1}{2\a}}\delta^2}}{y^{\frac{1}{2\a}}}\right) \right)\\
&=\sum_{n=0}^{\infty}\sum_{h=0}^{\infty}p_n q_{n,2h} \Gamma(h+1/2) i^{\tilde{M}-n+\frac32+h}
 y^{\frac{\tilde{M}-n-h+\frac12}{2\a}}+O\left(y^{\frac{\tilde{M}}{2\a}}e^{-y^{\frac{1}{2\a}}\delta^2}\right).
\end{align*}
Hence for any $m>0$, we get
\begin{equation} 
J_0=C_0e^{-2\a iy^{\frac{1}{2\a}}}e^{\pi i (\a r+\mu-\frac12)}\sum_{k=0}^m c_k' i^{\tilde{M}-k+\frac32} y^{\frac{\tilde{M}-k+\frac12}{2\a}}
+O\left(y^{\frac{\tilde{M}-m-1/2}{2\a}}\right), \notag
\end{equation}
where we put
$$
c_k'=\sum_{n+h=k}p_n q_{n,2h}\Gamma(h+1/2)(-1)^h.
$$

Next we consider the integral on $\mathcal{L}_2$. Let
$$
J_2=\int_{\mathcal{L}_2}\frac{\Delta_2(s)}{\Delta_1(r-s)}s^{\omega}F(s)x^{-s}ds.
$$
From the definition of $\mathcal{L}_2$, we have
\begin{align*}
|J_2| & \leq \int_{v_2}^{\infty} \left|\frac{\Delta_2(s)}{\Delta_1(r-s)}s^{\omega}F(s)x^{-s}\right|y^{\frac{1}{2\a}}dv
\end{align*}
with $s=y^{\frac{1}{2\a}}w$, $(w=u_2+iv)$. By using \eqref{integrad}, the above integrand is bounded as
$$
\ll x^{-u_2}y^{\frac{\tilde{M}+1}{2\a}}|w|^{\tilde{M}}|\exp(f(s))|.
$$
Since
\begin{align}
f(y^{\frac{1}{2\a}}w)&=2\a y^{\frac{1}{2\a}}w\left(\frac{1}{2\a}\log y+\log w\right)-2\a y^{\frac{1}{2\a}}w  
-\pi i \a y^{\frac{1}{2\a}}w-y^{\frac{1}{2\a}}w\log y \notag \\
&= 2\a y^{\frac{1}{2\a}}\left(w\log w - w -\frac{\pi i }{2}w\right), \notag 
\end{align}
we find that 
\begin{align*}
|J_2| &\ll x^{-u_2}y^{\frac{\tilde{M}+1}{2\a}} \max_{v_2 \leq v < \infty} G(v) \int_{v_2}^{\infty}|w|^{-2}dv,
\end{align*}
where we put
\begin{align*}
G(v)=v^{\tilde{M}+2}\exp\left(2\a y^{\frac{1}{2\a}}\Re\left(w\log w - w -\frac{\pi i }{2}w\right)\right).
\end{align*}
We shall see that $G(v)$ is a decreasing function. Differentiate $G(v)$ with respect to $v$, then we have
\begin{align*}
\frac{d}{dv}G(v)=&v^{\tilde{M}+1}\exp\left(2\a y^{\frac{1}{2\a}}\Re\left(w \log w -w -\frac{\pi i }{2}w\right)\right) \\
                 & \times \left(\tilde{M}+2-2\a y^{\frac{1}{2\a}}v\left(\arg w -\frac{\pi}{2}\right)\right).
\end{align*}
Note that
$$
\lim_{v \to \infty} v\left(\arg w -\frac{\pi}{2}\right)=-u_2>0.
$$
Therefore for sufficiently large $x$ (hence for sufficiently large $y$), $G(v)$ is a decreasing function.
Thus $G(v)$ attains its maximal value at $v=v_2$, namely, we have
\begin{align*}
|J_2| \ll & x^{-u_2}y^{\frac{\tilde{M}+1}{2\a}}\exp\left(2\a y^{\frac{1}{2\a}}
\Re\left(w_2\left(\log w_2 - 1-\frac{\pi i }{2}\right)\right)\right).
\end{align*}
Let
$$
\delta_2=\arg w_2 >\frac{\pi}{2}.
$$
Then
$$
\Re\left(w_2\left(\log w_2 - 1-\frac{\pi i }{2}\right)\right)=u_2(\log|w_2|-1)-v_2\left(\delta_2-\frac{\pi}{2}\right),
$$
which is a constant smaller than 0 (since $\delta$ is also small). Consequently $J_2$ is exponentially
decayed as $x \to \infty$.

Similarly it is shown that $J_1$ is also exponentially decayed as $x \to \infty$. 

Collecting these formulas, we get
\begin{align}
J:&=\int_{\mathcal{L}_0+\mathcal{L}_1+\mathcal{L}_2}\frac{\Delta_2(s)}{\Delta_1(r-s)}s^{\omega}F(s)x^{-s}ds \notag \\
&=C_0e^{-2\a iy^{\frac{1}{2\a}}}e^{\pi i (\a r+\mu-\frac12)}\sum_{k=0}^m c_k' i^{\tilde{M}-k+\frac32}
 y^{\frac{\tilde{M}-k+\frac12}{2\a}}+ O\left(y^{\frac{\tilde{M}-m-1/2}{2\a}}\right) \notag
\end{align}
for any $m$. Noting that $\mathcal{I}(x)=J-\bar{J}$ and $y=xe^{-(\tau+\tau')}$, we conclude that
\begin{align}
&\mathcal{I}(x)
=2i C_0 \sum_{k=0}^{m} c_k' y^{\frac{\tilde{M}-k+\frac12}{2\a}}
 \sin\left(-2\a y^{\frac{1}{2\a}}+ \left(\tilde{M}-k+2\a r+2\mu+\frac12\right)\frac{\pi}{2}\right) \notag \\
& \quad +O\left(y^{\frac{\tilde{M}-m-\frac12}{2\a}}\right) \notag \\
&=2i C_0 e^{-\frac{(\tau+\tau')(\tilde{M}+\frac12)}{2\a}}\sum_{k=0}^{m} c_k x^{\frac{\tilde{M}-k+\frac12}{2\a}}
 \cos\left(hx^{\frac{1}{2\a}}-\left(M+\omega-k+\a r+\mu'+\mu-\frac12\right)\frac{\pi}{2}\right) \notag \\
& \quad +O\left(x^{\frac{\tilde{M}-m-\frac12}{2\a}}\right), \notag
\end{align}
where we write
$$
c_k=c_k'e^{\frac{(\tau+\tau')}{2\a}k}.
$$
This completes the proof of Lemma 2.
\qed

\medskip

\begin{lem}
We assume \eqref{gamma-katei} in Lemma 2. Let $F(s)$ be a function which satisfies the same assumptions in Lemma 2.
Define the function $\mathcal{J}(x)$ by
\begin{equation*}
\mathcal{J}(x)=\frac{1}{2\pi i}\int_{\mathcal{C}_{a,b}}\frac{\Gamma(r-s)\Delta_2(s)}{\Gamma(r+\rho-s)\Delta_1(r-s)}F(s)x^{-s}ds,
\end{equation*}
with real $\varrho$. Then there exist constants $a_l$ and $c_l$ such that for any positive integer $m$, we have
\begin{align}
\mathcal{J}(x)=&\sum_{l=0}^m a_l  x^{\frac{\tilde{M}-l+1/2}{2\a}}\cos(hx^{\frac{1}{2\a}}+c_l \pi)+O(x^{\frac{\tilde{M}-m-1/2}{2\a}})
+O(x^{-b}), \label{J-asym}
\end{align}
where
$$
\tilde{M}=\mu'-\mu-\a r-\varrho+M.
$$
If we take $b$ large, we have
\begin{equation} \label{J-b}
\mathcal{J}(x) \ll x^{\frac{\mu'-\mu-\a r-\varrho+M+1/2}{2\a}}.
\end{equation}

\end{lem}

\proof  The right hand side except the last error term in \eqref{J-asym} are obtained by putting $\omega=-\rho$ in Lemma 2.
When we deform the path of integration $\mathcal{C}_{a,b}$ to the path  $\mathcal{L}_0 \cup
\mathcal{L}_1\cup \mathcal{L}_2$ and its complex conjugate, it may pass across the poles of $\Gamma(r-s)$ greater than $b$,
but they are finite depending on $m$. Hence the contribution from these poles is $O(x^{-b})$.
The assertion \eqref{J-b} is obtained by taking $m=0$ in \eqref{J-asym}.
\qed

\medskip

We note that if $F(s)=P(s)/Q(s)$ is a rational function of $s$, we can take $M$ as $M=\deg P -\deg Q$.

\medskip


\begin{lem} Let $x>0$ and $\varrho>\min\{2\alpha a-r\a +\mu'-\mu, -1\}$. Define the function $f_{\varrho}(x)$ by
\begin{equation*} 
f_\varrho(x)=\frac {1}{2\pi i}\int_{\mathcal{C}_{a,b}}G_\varrho(s)x^{r+\varrho-s}ds
\end{equation*}
where
\begin{align*}
G_\varrho(s)=\frac {\Gamma(r-s)\Delta_2(s)}{\Gamma(r+\varrho+1-s)\Delta_1(r-s)}.
\end{align*}
Then
\begin{align*}
\frac {d}{dx}f_{\varrho+1}(x)=f_{\varrho}(x)
\end{align*}
and for any non-negative integer $m$, we have an asymptotic expansion
\begin{align}
f_\varrho(x)=&\sum_{l=0}^{m}\kappa_l x^{\theta_\varrho-\frac{l}{2\alpha}}\cos \left( hx^{\frac{1}{2\alpha}}+c_l\pi\right)
+O\left( x^{\theta_\varrho-\frac{m+1}{2\alpha}}+x^{r+\varrho-b}\right) ,  \label{f-asymptotic}
\end{align}
where
\begin{align*}
c_l=-\frac 12\left(\mu+\mu^{\prime}+r\alpha+ \varrho+\frac 12\right)+\frac l2
\end{align*}
and $\kappa_l$ are constants. In particular,
$$
\kappa_0=(2\pi)^{\frac{N'-N}{2}}\left(\frac{2\alpha}{h}\right)^\varrho \sqrt{\frac{2}{h\pi}}\, e^{\nu'-\nu-\frac{(\tau'+\tau)(\mu'-\mu)}{2\a}
+\frac{r}{2}(\tau'-\tau)},
$$
and
$$
\kappa_1=\kappa_0 \frac{2\a}{h} \left(-\frac{1}{2\a}\left(\frac{(\mu'-\mu-\a r-\varrho-1)^2}{2}-\frac{1}{24}\right)+B\right),
$$
where we put
$$
B=\frac12\left(\sum_{h=1}^{N'}\frac{B_2(\b_h')}{\a_h'}+\sum_{j=1}^{N}\frac{B_2(\a_j r+\b_j)}{\\a_j}-B_2(r)+B_2(r+\varrho+1)\right).
$$
\end{lem}

\proof
\eqref{f-asymptotic} is obtained from Lemma 3. The other assertions are calculated explicitly by using Lemma 2.   \qed

%


\section{The Tong-type identity for the error term}


 In this section we shall prove the Tong-type identity for the integral (or multiple integral) of the error term $E_\varrho(y)$,
which is a generalization of Theorem 1 in Tong \cite{To1} for the error term in the Piltz divisor problem.
The main result of this section is

\begin{thm} Let $x\ge 1$, $L \ge 0$. Suppose that
\begin{equation*}
\varrho+k> 2\alpha\sigma_b^*-\alpha r-\frac12+\mu'-\mu.  \label{rho-k-initial-condition}
\end{equation*}
Then for $x+ky>0$ we have
\begin{align}
&\int_0^y\cdots\int_0^y E_\varrho(x+y_1+\cdots+y_k)dy_1\cdots dy_k \nonumber\\
&=\sum_{\mu_n\le L}\frac{b(n)}{\mu_n^{r+\varrho}}\int_0^y\cdots\int_0^y
f_\varrho\left((x+y_1+\cdots+y_k)\mu_n\right)dy_1\cdots dy_k\nonumber\\
& \quad +\sum_{j=0}^{k}(-1)^{k-j}\binom kj \sum_{\mu_n>L}\frac{b(n)}{\mu_n^{r+\varrho+k}}f_{\varrho+k}\left((x+jy)\mu_n\right).\nonumber
\end{align}
\end{thm}

\medskip

\begin{rem} Cao \cite{Cao} and Fomenko \cite{Fo3} used Tong's idea in \cite{To1} to study the asymmetric many-dimensional
problem and mean value theorems for automorphic $L$-functions, respectively.

\end{rem}

\proof

We assume first that $\varrho \geq 0$ and $\varrho+k>\a(2c-r)+\mu'-\mu$.
Then by  Perron's formula (see, Ivi\'c \cite[(A.8) and (A.4)]{I1}) we have
\begin{align*}
A_\varrho(x)=\frac {1}{2\pi i}\int_{c-i\infty}^{c+i\infty}
\frac {\Gamma(s)\varphi(s)x^{s+\varrho}}{\Gamma(s+\varrho+1)}ds,
\end{align*}
where $c$ is the constant chosen in Section 2. By the residue theorem  and \eqref{eq2-7} we have
\begin{align*}
&E_\varrho(x)=\frac {1}{2\pi i}\int_{\mathcal{C}_{c,r-b}}
\frac {\Gamma(s)\varphi(s)x^{s+\varrho}}{\Gamma(s+\varrho+1)}ds,\\
\intertext{and}
&\int_{\lambda_1}^x E_\varrho(u) du=\frac {1}{2\pi i}\int_{\mathcal{C}_{c,r-b}}
\frac {\Gamma(s)\varphi(s)x^{s+\varrho+1}}{\Gamma(s+\varrho+2)}ds+c_1
\end{align*}
with some constant $c_1$. By induction we deduce
\begin{align*}
\int_0^y E_\varrho(x+y_1) dy_1=&\int_{\lambda_1}^{x+y} E_\varrho(u)du-\int_{\lambda_1}^{x} E_\varrho(u)du\\
=&\frac {1}{2\pi i}\int_{\mathcal{C}_{c,r-b}}
\frac {\Gamma(s)\varphi(s)\left((x+y)^{s+\varrho+1}-x^{s+\varrho+1}\right)}
{\Gamma(s+\varrho+2)}ds,
\end{align*}
and
\begin{align}  \label{averageE}
&\int_0^y\cdots\int_0^y E_\varrho(x+y_1+\cdots+y_k) dy_1\cdots dy_k  \\
&=\frac {1}{2\pi i}\int_0^y\cdots\int_0^ydy_2\cdots dy_k\int_{\mathcal{C}_{c,r-b}}
\frac {\Gamma(s)\varphi(s)\left((x+y+\sum\limits_{j=2}^ky_j)^{s+\varrho+1}-
(x+\sum\limits_{j=2}^ky_j)^{s+\varrho+1}\right)}
{\Gamma(s+\varrho+2)}ds \nonumber\\
&=\frac {1}{2\pi i}\int_{\mathcal{C}_{c,r-b}}
\frac {\Gamma(s)\varphi(s)}{\Gamma(s+\varrho+2)}
\int_0^y\cdots\int_0^y
\left(\Bigl(x+y+\sum_{j=2}^ky_j\Bigr)^{s+\varrho+1}-\Bigl(x+\sum_{j=2}^ky_j\Bigr)^{s+\varrho+1}\right)
dy_2\cdots dy_k \nonumber\\
&=\sum_{j=0}^{k}(-1)^{k-j}\binom kj\frac {1}{2\pi i}\int_{\mathcal{C}_{c,r-b}}
\frac {\Gamma(s)\varphi(s)(x+jy)^{s+\varrho+k}}
{\Gamma(s+\varrho+k+1)}ds. \notag
\end{align}
Furthermore, under the condition $\varrho+k>\a(2c-r)+\mu'-\mu$. we can change the path of integration in \eqref{averageE} to $\mathcal{C}_{r-c,r-b}$,
and by using the functional equation along with a change of variable from $s$ to $r-s$, we have
\begin{align*}
\frac {1}{2\pi i}\int_{\mathcal{C}_{c,r-b}}
\frac {\Gamma(s)\varphi(s)u^{s+\varrho+k}}
{\Gamma(s+\varrho+k+1)}ds=\frac {1}{2\pi i}\int_{\mathcal{C}_{c,b}}
G_{\varrho+k}(s)\psi(s)u^{r+\varrho+k-s}ds \quad (u>0).
\end{align*}
Since $b>c>\sigma_b^*$, $\psi(s)$ is expressed as a Dirichlet series. Exchanging the order of summation and integration, we get
\begin{align}
\sum_{n=1}^{\infty}\frac {b(n)}{\mu_n^{r+\varrho+k}}\frac {1}{2\pi i}\int_{\mathcal{C}_{c,b}}
G_{\varrho+k}(s)(u\mu_n)^{r+\varrho+k-s}ds=\sum_{n=1}^{\infty}\frac {b(n)}{\mu_n^{r+\varrho+k}}f_{\varrho+k}(u\mu_n) \label{series}
\end{align}
for $a<c$. Since $f_{\varrho+k}(u\mu_n) \ll (u\mu_n)^{\theta_{\varrho+k}}$, the series \eqref{series} converges absolutely if
\begin{equation}
\varrho+k> 2\a \sigma_b^{\ast}-\a r-\frac12+\mu'-\mu. \label{conv-condition}
\end{equation}
Thus under this assumption, the exchange of the order of summation and integration is justified and
the equality
\begin{align} \label{range-ext}
&\int_0^y\cdots\int_0^y E_\varrho \left(x+\sum_{j=1}^ky_j\right) dy_1\cdots dy_k  \\
& \qquad =\sum_{j=0}^{k}(-1)^{k-j}\binom kj \sum_{n=1}^{\infty}\frac{b(n)}{\mu_n^{r+\varrho+k}}f_{\varrho+k}\left((x+jy)\mu_n\right) \notag
\end{align}
is valid under \eqref{conv-condition}.

Finally we note that if $\min\{u,u+y\}>0$, then for any non-negative integer $j$ we have
\begin{align} \label{indexshift-0}
&\int_0^yf_{\varrho+j}\left((u+y_1)\mu_n\right)dy_1=
\int_0^ydy_1\frac {1}{2\pi i}\int_{\mathcal{C}_{a,b}}
G_{\varrho+j}(s)\left((u+y_1)\mu_n\right)^{r+\varrho+j-s}ds \\
&=\frac {1}{2\pi i \mu_n}\int_{\mathcal{C}_{a,b}} \frac{\Gamma(r-s)\Delta_2(s)}{\Gamma(r+\varrho+j+2-s)\Delta_1(r-s)} \notag \\
& \hspace{2cm} \times \left(\left((u+y_1)\mu_n\right)^{r+\varrho+j+1-s}-(u\mu_n)^{r+\varrho+j+1-s}\right)ds\nonumber\\
&=\frac {1}{\mu_n}\left(f_{\varrho+j+1}(u+y)-f_{\varrho+j+1}(u)\right).\nonumber
\end{align}
By induction we get
\begin{align}
\int_0^y\cdots\int_0^y f_{\varrho}\left(\Bigl(x+\sum_{j=1}^ky_j\Bigr)\mu_n\right) dy_1\cdots dy_k 
=\frac {1}{\mu_n^k}\sum_{j=0}^{k}(-1)^{k-j}\binom kj f_{\varrho+k}\left((x+jy)\mu_n\right).\label{indexshift}
\end{align}
Theorem 4 follows from \eqref{range-ext} and \eqref{indexshift} at once.
\qed



\section{Large and   small values   of $E_{\varrho}(y)$ in   short intervals}


In this section we shall study the large and small values of the error term $E_\varrho(y)$ for $y$
in some "short" interval. We always suppose that  $\alpha>\frac 12$.

We first state some conditions.
\begin{enumerate}
\item[(A1)] $b(n)\ll _\varepsilon {\mu_n}^{\ell_1+\varepsilon}$ for some $\ell_1\ge 0$ and $b(n)\gg 1$ for infinitely many $n$.

\item[(A2)] $n^{\hbar}\ll \mu_n\ll n^{\hbar} $ for some $0<\hbar $, and $\mu_m^c-\mu_n^c\gg m^{\hbar c}-n^{\hbar c}$
for $0<c\le\frac {1}{\hbar}$ and $m>n$.

\item[(A3)] $a(n)\ll _\varepsilon {\lambda_n}^{\ell+\varepsilon}$ for some $\ell\ge 0$ and $a(n)\gg 1$ for infinitely many $n$.

\item[(A4)] $a(n)\ge 0$ for all  $n\in \mathbb{N}$ and $a(n)\gg 1$ for infinitely many $n$.
\end{enumerate}


\begin{thm} Let $\lambda \ge 1$ be a fixed real number. Assume that  the conditions (A1) and (A2) hold.
Then there exists a constant $B>0$ such that   for $x\ge x_0 $ and $ Bx^{1-\frac {1}{2\alpha}}\le U\le x$ we have
\begin{align}
\int_x^{x+U}|E_\varrho(y)|^{\lambda} dy\gg x^{\lambda\theta_\varrho}U,
\quad \int_1^{x}|E_\varrho(y)|^{\lambda} dy\gg x^{1+\lambda\theta_\varrho},  \label{Eq5-1}
\end{align}
where $\theta_\varrho $ is defined by (3.6).
In particular
\begin{align}
\int_x^{x+U}|E_\varrho(y)| dy\gg x^{\theta_\varrho}U,\quad \int_1^{x}|E_\varrho(y)| dy\gg x^{1+\theta_\varrho}. \label{Eq5-2}
\end{align}
\end{thm}


\begin{thm} Under the conditions (A1) and (A2), there exist two positive constants $B$ and $B'$ such that for $x\ge x_0$,
we have
$$
\max\limits_{x\le y\le x+Bx^{1-\frac {1}{2\alpha}}}\pm E_{\varrho}(y)>B' x^{\theta_\varrho}.
$$
Furthermore we have:
\begin{enumerate}
\item[{\rm (i)}]  If the conditions (A1), (A2) and (A3)  hold, then there exists a point $x^*\in  [x,x+Bx^{1-\frac {1}{2\alpha}}]$
with $|E(x^*)|\ll x^{\ell+\varepsilon}$ and $E(y)$ is continuous at $x^*$.
\item[{\rm (ii)}] If the conditions (A1), (A2) and (A4) hold,  then for any $t$ with $|t|<B' x^{\theta_0}$ there exists
at least a point $x^*\in  [x,x+Bx^{1-\frac {1}{2\alpha}}]$ such that  $E(x^*)=t$ and $E(y)$ is continuous at $x^*$.
\end{enumerate}
\end{thm}

\begin{rem}
Tong \cite{To1} first proved Theorem 6 for the error term in the summatory function of the the general Piltz divisor problem $d_{\ell}(n)$
for any $\ell\geq 2$. A. Ivi\'c \cite{I2}  studied the general case for $\Delta_1(s)=\Delta_2(s).$
Our result here provides a new proof of Ivi\'c's result.
\end{rem}

\begin{rem}
When $\varrho>0$,  $E_{\varrho}(y)$ is a continuous function of $y$. Then for any $t$ with $|t|<B' x^{\theta_\varrho}$
there exists at least a point $x^*\in  [x,x+Bx^{1-\frac {1}{2\alpha}}]$ such that  $E(x^*)=t$ and $E(y)$ is continuous at $x^*$.
\end{rem}

\begin{rem}
For $E_\varrho^{*}(y),$ we have similar results. Suppose that $\{\lambda_n\}(n\ge 1)$ satisfy

\noindent (A2*): $n^{\hbar^{*}}\ll \lambda_n\ll n^{\hbar^{*}} $ for some $0<\hbar^{*} $,
and $\mu_m^{c^{*}}-\mu_n^{c^{*}}\gg m^{\hbar^{*} c^{*}}-n^{\hbar^{*} c^{*}}$ for $0<c^{*}\le\frac {1}{\hbar^{*}}$ and $m>n$.

So we have the following results:
\begin{enumerate}
\item[(1)] If (A3) and (A2*) hold, then Theorem 5 holds for $E_\varrho^{*}(y)$ with $\theta_\varrho$ replaced by
$$\theta_\varrho '=: \frac{r}{2}-\frac{1}{4\a}+\varrho\left(1-\frac{1}{2\a}\right)
+\frac{\mu-\mu'}{2\a}.$$

\item[(2)] If (A3) and (A2*) hold,   there exists two positive constants $B$ and $B'$ such that for $x\ge x_0$,
we have $$\max\limits_{x\le y\le x+Bx^{1-\frac {1}{2\alpha}}}\pm E_{\varrho}^{*}(y)>B' x^{\theta_\varrho ' }.$$

\item[(3)]  If the conditions (A3), (A2*) and (A1)  hold, then there exists a point
 $x^*\in  [x,x+Bx^{1-\frac {1}{2\alpha}}]$ with $|E^{*}(x^*)|\ll x^{\ell_1+\varepsilon}$ and $E^{*}(y)$ is continuous at point $x^*$.

\item[(4)] If the conditions (A3), (A2*) hold,   $b(n)\ge 0(n\ge 1) $  and $b(n)\gg 1$ for infinitely many $n,$
 then for any $t$ with $|t|<B' x^{\theta_0 '}$
there exists at least a point $x^*\in  [x,x+Bx^{1-\frac {1}{2\alpha}}]$ such that  $E^{*}(x^*)=t$ and $E^{*}(y)$ is continuous at point $x^*$.
\end{enumerate}

\end{rem}

\medskip

{\it Proof of Theorem 5.}  We let $k$ be a fixed large integer such that $v=r+\varrho+k-\theta_{\varrho+k}>\sigma_b^*+1/2$.
Therefore the series
 \begin{equation*}
 g(t)=\sum_{n=1}^{\infty}\frac {b(n)}{\mu_n^v}\cos\left(
 h(t\mu_n)^{\frac {1}{2\alpha}}+c_0\pi  \right)
 \end{equation*}
is absolutely convergent. By Lemma 2 of Ivi\'c \cite{I2}, there exist two constants $B'>0, D'>0$ such that for $x\ge x_0$ every
interval $[x,x+D'x^{1-\frac {1}{2\alpha}}]$ contains two points $x_1, x_2$ for which
 \begin{equation*}
g(x_1)>B',\quad g(x_2)<-B'.
 \end{equation*}

 Take $U'=c' x^{1-\frac {1}{2\alpha}}$, where $c'$ be a constant to be specified below, and $x^*=x_j$ for $j=1$ or 2.
Applying Theorem 4 with $L=0$, Lemma 4 and \eqref{indexshift-0} we have
\begin{align}
&\int_{0}^{U'}dy \int_0^y\cdots\int_0^y E_{\varrho}(x^*+y_1+\cdots+y_k)dy_1\cdots dy_k\nonumber\\
&=\int_{0}^{U'}\left(
\sum_{j=0}^{k}(-1)^{k-j}\binom kj
\sum_{n=1}^{\infty}\frac{b(n)}{\mu_n^{r+\varrho+k}}
f_{\varrho+k}\left((x^*+jy)\mu_n\right)
\right)dy\nonumber\\
&=(-1)^k \kappa_0(\varrho+k) U' g(x^*)(x^*)^{\theta_{\varrho+k}}+ O\left(U'x^{\theta_{\varrho+k}-\frac {1}{2\alpha}}\right)\nonumber\\
&\quad +
\sum_{j=1}^{k}\frac{(-1)^{k-j}}{j}\binom kj
\sum_{n=1}^{\infty}\frac{b(n)}{\mu_n^{r+\varrho+k+1}}
\left(f_{k+\varrho+1}\left((x^*+jU')\mu_n\right)
-f_{k+\varrho+1}(x^*\mu_n)
\right).\nonumber
\end{align}
Using Lemma 4 again we obtain
\begin{align}\label{Eq5-5}
&\left|\int_{0}^{U'}dy \int_0^y\cdots\int_0^y E_{\varrho}(x^*+y_1+\cdots+y_k)dy_1\cdots dy_k\right|  \\
&\ge (x^*)^{\theta_{\varrho+k+1}}\left(c' \kappa_0(\varrho+k)B'+O\left(\frac {U'}{x}\right) \right. \nonumber \\
& \qquad \left. -\kappa_0(\varrho+k+1)2^{k+1} \psi_1(r+\varrho+k+1-\theta_{\varrho+k+1})
                 \left(1+\frac {kU'}{x}\right)^{\theta_{\varrho+k+1}} \right), \nonumber
\end{align}
where $\psi_1(s)=\sum_{n=1}^{\infty}\frac {|b(n)|}{\mu_n^s}$.
Now we can choose a constant $c'$ such that
\begin{equation}
c'>\frac {\kappa_0(\varrho+k+1)2^{k+2} \psi_1(r+\varrho+k+1-\theta_{\varrho+k+1})}{\kappa_0(\varrho+k)B'}. \label{Eq5-6}
\end{equation}
Combining  \eqref{Eq5-5} and \eqref{Eq5-6} we have for some constant $B_0>0$ as $x>x_0$
\begin{align}\label{eq5-7}
B_0x^{\theta_{\varrho+k+1}}&\le \left|\int_{0}^{U'}dy \int_0^y\cdots\int_0^y E_{\varrho}(x^*+y_1+\cdots+y_k)dy_1\cdots dy_k\right| \\
&\le \int_{0}^{U'}dy \int_0^y\cdots\int_0^y |E_{\varrho}(x^*+y_1+\cdots+y_k)|dy_1\cdots dy_k\nonumber\\
&\le \int_{0}^{U'}dy\int_{0}^{U}|E_{\varrho}(x+y_1)|dy_1
\int_0^y\cdots\int_0^ydy_2\cdots dy_k\nonumber\\
&\le x^{k(1-\frac {1}{2\alpha})}\int_{x}^{x+U}|E_{\varrho}(y)|dy,\nonumber
\end{align}
here $U=( k+2)U'$.
This proves the first assertion in \eqref{Eq5-2}.

Next, we have
\begin{align*}
\int_{1}^{x}|E_{\varrho}(y)|dy\ge \sum_{j=0}^{\frac {1}{2B}x^{\frac {1}{2\alpha}}}\int_{\frac x2 +jBx^{1-\frac {1}{2\alpha}}}^{\frac x2 +(j+1)Bx^{1-\frac {1}{2\alpha}}}|E_{\varrho}(y)|dy\gg x^{\frac {1}{2\alpha}}
\left(x^{\theta_{\varrho}}x^{1-\frac {1}{2\alpha}}\right)=x^{1+\theta_{\varrho}},
\end{align*}
and this completes the proof of the second assertion in \eqref{Eq5-2}.

By H\"{o}lder's inequality and \eqref{Eq5-2} we immediately get \eqref{Eq5-1}.

\bigskip

{\it Proof of Theorem 6.}
Similarly to \eqref{eq5-7} we see that for $j=1,2$
\begin{align}
x^{\theta_{\varrho+k+1}}&\ll \max\pm\int_{0}^{U'}dy \int_0^y\cdots\int_0^y E_{\varrho}(x_j+y_1+\cdots+y_k)dy_1\cdots dy_k\nonumber\\
&\ll \max_{x\le y\le x+U}\pm E_{\varrho}(y)\times \int_{0}^{U'}dy
\int_0^y\cdots\int_0^ydy_2\cdots dy_k.\nonumber
\end{align}
Hence for $x>x_0$
\begin{equation}
\max_{x\le y\le x+U}\pm E_{\varrho}(y)\gg x^{\theta_{\varrho}}.  \label{eq5-8}
\end{equation}

When $\varrho=0$, for every interval $[x,x+2U]$, it follows from \eqref{eq5-8} that there exists
two points $y',y''\in [x,x+2U]$ such that
\begin{align*}
E(y')>0, E(y'')<0, \quad y'<y''.
\end{align*}

If the condition (A4) holds,  then $Q(y),$ which is  the main term  of the summatory function $a(n),$    is of the form (2.6) with $\varrho=0$,
   always positive and increasing. And
 the function $E(y)$ is a continuous and strictly decreasing in each interval $(\lambda_n,\lambda_{n+1})$.
Moreover, $E(y)$ has a jump up to $a(n)/2$ at the point $\lambda_n \ \ (n=1,2,\cdots)$, i.e.
$E(\lambda_n)=E(\lambda_{n}-0)+a(n)/2$. Hence there exists at least a zero
point $x^*$ of $E(x)$ in the interval $(y',y'')$ with $E(x)$ is continuous at point $x^*$.

Otherwise, we may suppose that for any $x\in (y',y'')$ with $E(x)\ne 0$. Without loss of generality,
we assume that for all $\lambda_n\in (y',y'') $ with $E(\lambda_{n}+0)E(\lambda_{n+1}-0)>0$.
Hence we can find at least a $\lambda_n\in (y',y'')$ such that $E(\lambda_{n}-0)\times E(\lambda_{n}+0)<0$,
in this time we can choose $x^*=\frac {\lambda_{n}+\lambda_{n+1}}{2}$.
This completes the proof of Theorem 6.


\section{The truncated Tong-type formula for an integral involving the error term}


Suppose $x>10$ is a large parameter. For any integrable function $g(y)$, we define
\begin{gather*}
 \tilde{y}:=y+\frac 1x(y_1+\cdots+y_k),\\
\int_{\mathbf{E}_k}g(\tilde{y})d\mathrm{Y}_k:=\int_0^1\cdots \int_0^1 g(\tilde{y})dy_1\cdots dy_k.
\end{gather*}
The integral $\int_{\mathbf{E}_k}E(\tilde{y})d\mathrm{Y}_k$ is very important in Tong's theory.
The aim of this section is to give its truncated expression
under some suitable  conditions.

Recall the definition of $A^{*}(y).$  Let $\hat{Q}(y)$ denote the sum of residues over all poles of the function
$\psi(s)y^ss^{-1}$ except the pole $s=0.$ And define the error term $\hat{E}(y)$ by
$$\hat{E}(y):=A^{*}(y)-\hat{Q}(y).$$
Since $\psi(s)$ doesn't have poles when $\sigma\leq \sigma^{*}$ by assumption, we have immediately that
\begin{equation*}
Q^{*}(y)=\hat{Q}(y)+\psi(0),\ \ E^{*}(y)=\hat{E}(y)-\psi(0).
\end{equation*}

Let $y\ge 1$, $\frac 12\le M\le N<\infty$, $\lambda$ be a real number, and define
\begin{align}\label{I-definition}
I(\lambda,M,N,y)=2\pi i\int_M^Nu^\lambda \hat{E}(u)\exp\left(-ih(uy)^{\frac {1}{2\alpha}}\right)du.
\end{align}

\medskip

We assume in this section that
\begin{enumerate}
\item[(B1)] There exists a constant $1-\sigma_b^{*}\leq \omega_0 < 1$ such that $|b(n)| \ll \mu_n^{\sigma_b^{\ast}-1+\omega_0}$.

\item[(B2)] $\sum_{\mu_n \leq y}|b(n)| \ll y^{\sigma_b^{\ast}} \log^A y, \ y\ge 2.$

\item[(B3)] There exists a constant $0 \leq \omega_1 \leq 1$  such that
$$ \int_1^T |\hat{E}(u)|du \ll T^{\sigma_b^{\ast}+\omega_1}. $$
\end{enumerate}

Throughout this section we let  $k$ to be a fixed  positive integer such that
\begin{equation} \label{k-condition}
k>\max\left\{2\alpha \sigma_b^*-\alpha r-\frac 12 -\varrho+\mu'-\mu,\alpha r \right\}.
\end{equation}
Furthermore we assume that
\begin{align}
\sigma_b^{\ast}+\omega_1-r-\frac{1}{4\a}-\frac{\varrho}{2\a}+\frac{\mu'-\mu}{2\a}-1<0. \label{xi}
\end{align}
We note that \eqref{k-condition} and \eqref{xi} implies
\begin{equation}
\sigma_b^{\ast}+\omega_1-\frac{r}{2}-\frac{1}{4\a}-1-\frac{\varrho}{2\a}+\frac{\mu'-\mu}{2\a}-\frac{k}{\a}<0. \label{xi-1}
\end{equation}


Let $1\le x\le y\le (1+\delta)x$,  $N=[x^{4\alpha-1-\varepsilon}]$ and $J=\left[(4\alpha^2 r+4\alpha)\varepsilon^{-1}\right]$,
where $\delta>0$ is a fixed small positive constant. Without loss of generality, we suppose $N\not= \mu_n$ for any $n.$
Furthermore let
\begin{align*}
\lambda_{0}:=\theta_\varrho+\frac{1}{2\a}-r-\varrho-1=-\frac{r}{2}+\frac{1}{4\a}-\frac{\varrho}{2\a}+\frac{\mu'-\mu}{2\a}-1.
\end{align*}
From Theorem 6 we can find a real number $M$ in every subinterval $[t, t+Bt^{1-\frac{1}{2a}}]$ in $[1,\sqrt N]$
such that $M \neq \mu_n$ for $n=1,2,\ldots$ and $\hat{E}(M)= 0$ if $b(n) \geq 0$ for all $n$ or $\hat{E}(M) \ll M^{\sigma_b^{*}-1+\omega_0}$
otherwise.

The truncated Tong-type formula can be stated as


\begin{thm} Assume the conditions  (A2), (B1), (B2) and (B3) hold.
Let $M \in [1,\sqrt{N}]$ be chosen as above. Then we have 
\begin{align}\label{E-I-expression}
\int_{\mathbf{E}_k}E_{\varrho}(\tilde{y})d\mathrm{Y}_k=
\sum_{j=1}^{6}R_j(y;M),
\end{align}
where
\begin{align}
R_1(y;M)= &\kappa_0 y^{\theta_\varrho} \sum_{\mu_n \leq M}\frac{b(n)}{\mu_n^{r+\varrho-\theta_\varrho}}\cos(h(y\mu_n)^{\frac{1}{2\a}}+c_0\pi),
          \label{R1-expression} \\[1ex]
R_2(y;M)= &y^{\theta_\varrho+\frac{1}{2\a}}\Re \left\{c_{00}I\left(\lambda_0, M,N,y\right)\right\}, \label{R2-expression} \\[1ex]
R_3(y;M)= &\Dsum_{l+m \neq 0}\Re \left\{c_{lm} I\left(\lambda_0+\frac{l-m}{2\a},M,N,y\right)\right\}
            x^{-l}y^{-l+\theta_{\varrho}+\frac{1}{2\a}+\frac{l-m}{2\a}}, \label{R3-expression} \\[1ex]
R_4(y;M)= &\sum_{j=0}^k \sum_{m=0}^k \Re\left\{c_{jm}' I\left(\lambda_0-\frac{k+m}{2\a},N,\infty,y+\frac{j}{x}\right)\right\}  \label{R4-expression} \\
          & \hspace{5cm} \times x^k\left(y+\frac{j}{x}\right)^{k+\theta_\varrho+\frac{1}{2\a}-\frac{k+m}{2\a}}, \nonumber \\[1ex]
R_5(y;M) \ll& x^{\theta_{\varrho}-\frac{1}{2\a}}M^{\max(\sigma_b^{\ast}-\frac{r}{2}-\frac{3}{4\a}-\frac{\varrho}{2\a}+\frac{\mu'-\mu}{2\a},0)}\log^AM
              \label{R5-expression}\\
          & + x^{-2+\theta_\varrho+\frac{1}{2\a}}M^{\max(\sigma_b^{\ast}-\frac{r}{2}+\frac{1}{4\a}-\frac{\varrho}{2\a}+\frac{\mu'-\mu}{2\a},0)}\log^A M \notag \\
          & +x^{\theta_{\varrho}-\frac{r}{2}}M^{\sigma_b^{\ast}+\omega_1-r-\frac{1}{4\a}-\frac{\varrho}{2\a}+\frac{\mu'-\mu}{2\a}-1} \notag \\
          & + x^{(4\a-1)(\sigma_b^{\ast}+\omega_1)-2k+r-\varrho+2(\mu'-\mu)+\frac{2k}{\a}-2\a r -4\a},  \notag \\[1ex]
R_6(y;M) \ll& \begin{cases}
0, & \mbox{if $b(n)\ge 0\ (n\ge 1)$,}\\
x^{\theta_\varrho} M^{\sigma_b^{\ast}-1+\omega_0-\frac{r}{2}-\frac{1}{4\a}-\frac{\varrho}{2\a}+\frac{\mu'-\mu}{2\a}}, & \mbox{if the condition (B1) holds.}
\end{cases} \label{R6-expression}
\end{align}
Here $\kappa_0, c_0$ and $c_{00}, c_{lm}, c_{jm}'$ are certain real and complex constants, respectively.
\end{thm}


\subsection{ Lemmas for generalized Bessel functions}


Before proving Theorem 7 we shall prepare some preliminary assertions for
the so-called generalized Bessel functions.

By induction one easily verifies ($r-s \notin -\mathbb{N}$)
\begin{align}\label{yrs-integral}
\int_{\mathbf{E}_k}\tilde{y}^{r-s}d\mathrm{Y}_k=x^k\sum_{l=0}^{k}
(-1)^{k-l}\binom kl \frac{(y+ lx^{-1})^{k+r-s}}{(1+r-s)(2+r-s)\cdots(k+r-s)}
\end{align}

For $j \geq 0$, define
\begin{align*}
f_{\varrho,j}(x)=\frac {1}{2\pi i}\int_{\mathcal{C}_{a,b}}\frac{\Gamma(r-s)}{(r+\varrho-s)^j\Gamma(r+\varrho+1-s)}
\frac{\Delta_2(s)}{\Delta_1(r-s)}x^{r+\varrho-s}ds
\end{align*}
and
\begin{equation*}
g_{\varrho,j}(u)=\int_{\mathbf{E}_k}f_{\varrho,j}(\tilde{y}u)d\mathrm{Y_k}.
\end{equation*}
By definition $f_{\varrho}(x)=f_{\varrho,0}(x)$.
We note that $f_{\varrho,j}(x)$ converges absolutely if
\begin{equation*}
a<\frac{r}{2}+\frac{\varrho+j-\mu'+\mu}{2\a}.
\end{equation*}

\begin{lem}
Let $x \leq y \leq (1+\delta)x$, $b$  is a sufficiently large positive constant.
Then we have
\begin{align}
f_{\varrho,j}(x) & \ll x^{\theta_{\varrho}-\frac{j}{2\a}}, \label{f-estimate}\\
\intertext{and}
g_{\varrho,j}(u) & \ll x^{2k}(ux)^{\theta_{\varrho}-\frac{k+j}{2\a}}. \label{g-estimate}
\end{align}
\end{lem}

\proof
The first assertion \eqref{f-estimate} is obtained immediately by replacing $\varrho$ by $\varrho+1$ and taking $M=-j$ in \eqref{J-b}.
($M$ is the parameter in the definition of $F(s),$ see  \eqref{def_function_F}.)

To prove \eqref{g-estimate}, substituting the definitions of $f_{\varrho,j}(x)$ into $g_{\varrho,j}(u)$ and using \eqref{yrs-integral},
we get
\begin{align*}
g_{\varrho,j}(u)
&=x^k\sum_{l=0}^k(-1)^{k-l}\binom{k}{l}\left(y+\frac{l}{x}\right)^{k} \\
& \quad \times \frac{1}{2\pi i} \int_{\mathcal{C}_{a,b}}\frac{\Gamma(r-s)}{(r+\varrho-s)^j\Gamma(r+\varrho+k+1-s)}\frac{\Delta_2(s)}{\Delta_1(r-s)}
\left(\left(y+\frac{l}{x}\right)u\right)^{r+\varrho-s}ds.
\end{align*}
Replacing $\varrho$ by $\lambda=\varrho+k+1$ and taking $M=-j$ in \eqref{J-b}, we get
\begin{align*}
g_{\varrho,j}(u) \ll x^{2k}(ux)^{\theta_{\varrho}-\frac{k+j}{2\a}}.
\end{align*}
\qed

We also need the asymptotic behaviour of the function $\Phi_\varrho(y)$ defined by
\begin{equation*} 
\Phi_{\varrho}(y)=\frac{d}{du}\left(\frac{f_{\varrho}(yu)}{u^{r+\varrho}}\right).
\end{equation*}

\begin{lem}
Let $l$ be a fixed positive integer and $x \leq y \leq (1+\delta)x$. Then there are constants $b_{lm}'$ and $b_{jm}''$ such that
for any positive integer $H$, we have
\begin{align}\label{G-bibun}
\frac{d^l}{dy^l}\Phi_{\varrho}(y)&=\sum_{m=0}^{H}\Re\left\{b_{lm}'\frac{(uy)^{\theta_{\varrho}+\frac{1}{2\a}+\frac{l-m}{2\a}}}{y^lu^{r+\varrho+1}}
e^{-ih(uy)^{\frac{1}{2\a}}}+O\left(\frac{(uy)^{\theta_{\varrho}+\frac{l-H}{2\a}}}{y^lu^{r+\varrho+1}}\right)\right\}  \\[1ex]
& \quad +O\left(y^{r+\varrho-b-l}u^{-1-b}\right), \notag
\intertext{and}
\int_{\mathbf{E}_k}\Phi_{\varrho}(\tilde{y})d\mathrm{Y_k}
&=\sum_{j=0}^{k} \left\{\sum_{m=0}^{H} \Re \left(b_{jm}''
\frac{x^k(y+\frac{j}{x})^{k+\theta_{\varrho}+\frac{1}{2\a}-\frac{k+m}{2\a}}}{u^{r+\varrho+1-\frac{1}{2\a}-\theta_{\varrho}+\frac{k+m}{2\a}}}
e^{-ih((y+\frac{j}{x})u)^{\frac{1}{2\a}}} \right) \right. \label{G-heikin} \\
& \qquad \left. +O\left(\frac{x^{2k+\theta_{\varrho}+\frac{1}{2\a}-\frac{k+H+1}{2\a}}}
{u^{r+\varrho+1-\frac{1}{2\a}-\theta_{\varrho}+\frac{k+H+1}{2\a}}}\right)\right\} +O\left(x^{2k+r+\varrho-b}u^{-1-b}\right). \notag
\end{align}
\end{lem}

\proof
From the definition of $f_{\varrho}(yu)$, we have
\begin{align}
\Phi_{\varrho}(y)=\frac{1}{2\pi i}\int_{\mathcal{C}_{a,b}}\frac{\Gamma(r-s)\Delta_2(s)(-s)}{\Gamma(r+\varrho+1-s)\Delta_1(r-s)}
y^{r+\varrho-s}u^{-s-1}ds, \label{G-shiki}
\end{align}
hence
\begin{align*}
\frac{d^l}{dy^l}\Phi_{\varrho}(y)&=\frac{1}{y^lu^{r+\varrho+1}}\frac{(-1)^{l-1}}{2\pi i}
\int_{\mathcal{C}_{a,b}}\frac{\Gamma(r-s)\Delta_2(s)}{\Gamma(r+\varrho+1-s)\Delta_2(r-s)} \\
& \quad \times s(s-r-\varrho)(s-r-\varrho+1)\cdots(s-r-\varrho+l-1)(yu)^{r+\varrho-s}ds.
\end{align*}
Replacing $\varrho$ by $\varrho+1$ and taking $M=l+1$ in the formula \eqref{J-asym}, we find that
\begin{align*}
\frac{d^l}{dy^l}\Phi_{\varrho}(y)
& \sim \frac{1}{y^lu^{r+\varrho+1}}\Bigl\{\sum_{m=0}^{\infty}b_{lm}(yu)^{\theta_\varrho+\frac{1}{2\a}+\frac{l-m}{2\a}}\cos(h(yu)^{\frac{1}{2\a}}+c_m\pi) \\
& \hspace{3cm}  +O((yu)^{r+\varrho-b})\Bigr\},
\end{align*}
which proves the assertion \eqref{G-bibun}.

Next we prove \eqref{G-heikin}. From \eqref{G-shiki} and \eqref{yrs-integral} we have
\begin{align}\label{G-heikin-2}
\int_{\mathbf{E}_k}\Phi_{\varrho}(\tilde{y})d\mathrm{Y_k}
&=\frac{1}{2\pi i}\int_{\mathcal{C}_{a,b}}\frac{\Gamma(r-s)\Delta_2(s)(-s)}{\Gamma(r+\varrho+1-s)\Delta_1(r-s)}
  u^{-s-1}\int_{\mathbf{E}_k}\tilde{y}^{r+\varrho-s}d\mathrm{Y_k}ds  \\
&=\frac{x^k}{u^{r+\varrho+1}}\sum_{j=0}^{k}(-1)^{k-j}\binom{k}{j}\left(y+\frac{j}{k}\right)^k \notag \\
& \quad \times \frac{1}{2\pi i}\int_{\mathcal{C}_{a,b}}\frac{\Gamma(r-s)\Delta_2(s)(-s)}{\Gamma(r+\varrho+k+1-s)\Delta_1(r-s)}
\left(\left(y+\frac{j}{x}\right)u\right)^{r+\varrho-s}ds. \notag
\end{align}
Replacing $\varrho$ by $\varrho+k+1$ and taking $M=1$ in \eqref{J-asym}, we find that
\begin{align*}
&\frac{1}{2\pi i}\int_{\mathcal{C}_{a,b}}\frac{\Gamma(r-s)\Delta_2(s)(-s)}{\Gamma(r+\varrho+k+1-s)\Delta_1(r-s)}
\left(\left(y+\frac{j}{x}\right)u\right)^{r+\varrho-s}ds \\
& \sim \sum_{m=0}^{\infty} b_{jm}'\left(\left(y+\frac{j}{x}\right)u \right)^{\theta_\varrho +\frac{1}{2\a}-\frac{k+m}{2\a}}
\cos\left(h\left(\left(y+\frac{j}{x}\right)u\right)^{\frac{1}{2\a}}+c_m\pi\right) \\
& \quad  +O\left((yu)^{r+\varrho-b}\right)
\end{align*}
with certain real constants $b_{jm}'$. Substituting this expression to \eqref{G-heikin-2}, we get \eqref{G-heikin}.
\qed


\subsection{ Proof of Theorem 7}


Now we turn to the proof of Theorem 7. Let $M \in [1, \sqrt{N}]$ be a real number as is chosen before, namely.
$M \neq \mu_n$ for $n=1,2,\ldots$ and $\hat{E}(M)= 0$ if $b(n) \geq 0$ for all $n$
or $\hat{E}(M) \ll M^{\sigma_b^{*}-1+\omega_0}$ otherwise. We always suppose that there are
infinitely many $n$ for which $b(n)\gg 1.$

Recall the Tong-type identity:
\begin{align*}
I:= \int_{\mathbf{E}_k}E_{\varrho}\left(y+\frac{\sum_{j=1}^k y_j}{x}\right)d\mathrm{Y}_k
=\sum_{n=1}^{\infty}\frac{b(n)}{\mu_n^{r+\varrho}}\int_{\mathbf{E}_k}f_{\varrho}\left(\left(y+\frac{\sum_{j=1}^k y_j}{x}\right)\mu_n\right)d\mathrm{Y}_k.
\end{align*}
In the right hand side we divide the sum over $n$ into two parts, namely
$$
I=\sum_{\mu_n \leq M} +\sum_{\mu_n>M}=:S_1(M)+S_2(M).
$$

We first treat $S_1(M)$. Let $F(\tilde{y})=f_{\varrho}(\tilde{y}u)$. By Taylor's formula  we have
$$
F(\tilde{y})=F(y)+F'(y_1)(\tilde{y}-y)
$$
with some $y \leq y_1 \leq \tilde{y}$. Then from the definition $f_{\varrho}(x)$ and \eqref{f-asymptotic} of Lemma 4 we have an asymptotic expansion
\begin{align}
F(y)& \sim \sum_{l=0}^{\infty} \kappa_l(yu)^{\theta_\varrho-\frac{l}{2\a}}\cos(h(yu)^{\frac{1}{2\a}}+c_l\pi)+O((yu)^{r+\varrho-b}) \label{F-asymp}
\intertext{and the upper bound}
    & \ll (yu)^{\theta_\varrho}  \notag
\end{align}
if we take $b$ large. Furthermore since $F'(y_1) = u f_{\varrho-1}(y_1u)$, we have
\begin{align*}
|F'(y_1)(\tilde{y}-y)| & \ll \frac{u}{x}(xu)^{\theta_{\varrho-1}}\\
                       & = x^{-2} (xu)^{\theta_{\varrho}+\frac{1}{2\a}}.
\end{align*}
Take out the term corresponding to $l=0$ in \eqref{F-asymp},   we get
\begin{align*}
S_1(M)&=\sum_{\mu_n \leq M}\frac{b(n)}{\mu_n^{r+\varrho}}\left\{\kappa_0(y\mu_n)^{\theta_{\varrho}}\cos(h(y\mu_n)^{\frac{1}{2\a}}+c_0\pi)
       +O((y\mu_n)^{\theta_{\varrho}-\frac{1}{2\a}}) \right. \\
& \hspace{2cm}  \left.  +O((y\mu_n)^{r+\varrho-b})+O\left(x^{-2} (x\mu_n)^{\theta_{\varrho}+\frac{1}{2\a}}\right)\right\}\\[1ex]
&=:R_1(y;M)+V_1+V_2+V_3.
\end{align*}
The error term $V_2$ can be negligible by taking $b$ large. On the hand, by (B2), we have
\begin{align}\label{V1-estimate}
V_1 &\ll x^{\theta_\varrho-\frac{1}{2\a}}\sum_{\mu_n \leq M}|b(n)|\mu_n^{-r-\varrho+\theta_\varrho-\frac{1}{2\a}}  \\
&\ll x^{\theta_\varrho-\frac{1}{2\a}}M^{\max(\sigma_b^{\ast}-\frac{r}{2}-\frac{3}{4\a}-\frac{\varrho}{2\a}+\frac{\mu'-\mu}{2\a},0)}\log^A M. \notag
\end{align}
Similarly we have
\begin{equation}
V_3 \ll x^{-2+\theta_\varrho+\frac{1}{2\a}}M^{\max(\sigma_b^{\ast}-\frac{r}{2}+\frac{1}{4\a}-\frac{\varrho}{2\a}+\frac{\mu'-\mu}{2\a},0)}\log^A M. \label{V3-estimate}
\end{equation}

Since $M \neq \mu_n$ for all $n$, it is easy to see that
if $K(u)$ is  any continuously differentiable function such that
$$K(u)Q^{*}(u)\rightarrow 0,\int_u^\infty Q^{*}(t)K^{\prime}(t)dt\rightarrow 0\ (u\rightarrow \infty), $$ then
 $$\int_M^\infty K(u)d\sum_{\mu_n\le u}b(n)=\int_M^\infty K(u)dA^{*}(u).$$
So for $S_2(M)$ we can  write
\begin{align}  
S_2(M)&=\int_M^\infty \left(\frac{1}{u^{r+\varrho}}\int_{E_k}f_{\varrho}(\tilde{y}u)d\mathrm{Y_k}\right)d\sum_{\mu_n\le u}b(n)\notag\\
&=\int_M^\infty \left(\frac{1}{u^{r+\varrho}}\int_{E_k}f_{\varrho}(\tilde{y}u)d\mathrm{Y_k}\right)dA^{*}(u) \notag \\
&=\int_M^\infty \left(\frac{1}{u^{r+\varrho}}\int_{E_k}f_{\varrho}(\tilde{y}u)d\mathrm{Y_k}\right)d\hat{Q}(u)
 +\int_M^\infty \left(\frac{1}{u^{r+\varrho}}\int_{E_k}f_{\varrho}(\tilde{y}u)d\mathrm{Y_k}\right)d\hat{E}(u) \notag \\
&=\int_M^{\infty}\frac{\hat{Q}'(u)}{u^{r+\varrho}}\int_{E_k}f_{\varrho}(\tilde{y}u)d\mathrm{Y_k}du
   +\left. \frac{\hat{E}(u)}{u^{r+\varrho}}\int_{E_k}f_{\varrho}(\tilde{y}u)d\mathrm{Y_k} \right|_{u=M}^{\infty} \notag \\
& \quad -\int_M^\infty \hat{E}(u)du \int_{\mathbf{E}_k} \frac{d}{du}\left(\frac {f_{\varrho}(\tilde{y}u)}{u^{r+\varrho}}\right)d\mathrm{Y}_k \nonumber\\
&=:S_{21}(M)+S_{22}(M)-S_{23}(M). \nonumber
\end{align}

We shall consider $S_{21}(M)$ first. From the definition of $f_{\varrho,j}(x)$, one easily checks that
\begin{align*}
\frac{d}{du} f_{\varrho,j+1}(yu)=\frac 1u f_{\varrho,j}(yu) \quad (1\le j\le k+1).
\end{align*}
Let $w_j(u)$ be the functions defined by
$$
w_1(u)=\hat{Q}'(u)u^{1-r-\varrho}
$$
and
$$
w_{j+1}(u)=w_j'(u)u\ (j\geq 1).
$$
Then by repeated integration by parts and the definition of $g_{\varrho,j}$ we have
\begin{align}\label{S1}
\int_M^{M_1}\frac{\hat{Q}'(u)}{u^{r+\varrho}}\int_{E_k}f_{\varrho}(\tilde{y}u)d\mathrm{Y_k}du
&=\sum_{j=1}^{k+1}(-1)^{j+1}\,  w_j(u)g_{\varrho,j}(u) \Big|_M^{M_1}   \\
& \quad + (-1)^{k+1}\int_M^{M_1}w_{k+1}'(u)g_{\varrho,k+1}(u)du. \notag
\end{align}
Recall that $\hat{Q}(u)$ is the sum of residues of of all poles of $\psi(s)x^ss^{-1}$ except $0$. Let $\gamma$ be
a pole with order $\nu$ and $\gamma_0$ the pole with the maximal real part.
Then the residue at $\gamma$ has the form $\hat{Q}_\gamma(u)=u^{\gamma}P(\log u)$, where $P(u)$ is
a polynomial of degree $\nu-1$. Let $w_{j,\gamma}(u)$ be functions defined by the same way for
$\hat{Q}_\gamma(u)$. Then we see easily
\begin{align*}
w_{1, \gamma}(u)&=u^{1-r-\varrho}(\gamma u^{\gamma-1}P(\log u)+u^{\gamma}P'(\log u)/u) \\
  & \ll u^{\Re\gamma-r-\varrho}\log^{\nu-1} u.
\end{align*}
Similarly we have
$$
w_{j,\gamma}(u) \ll u^{\Re\gamma-r-\varrho}\log^{\nu-1}u
$$
for $j \geq 2$.
Hence by \eqref{g-estimate}, we have
\begin{align*}
w_{j,\gamma}(M_1)g_{\varrho, j}(M_1) \ll M_1^{\Re\gamma-r-\varrho}(\log M_1)^{\nu-1}\, x^{2k}(xM_1)^{\theta_\varrho-\frac{k+j}{2\a}}.
\end{align*}
Since $k$ is chosen as \eqref{k-condition}, the exponent of $M_1$ is negative, $w_{j,\gamma}(M_1)g_{\varrho, j}(M_1)$,
as well as $w_{j}(M_1)g_{\varrho, j}(M_1)$ tend to $0$ when $M_1 \to \infty$.

Next we consider the term with $u=M$. Applying \eqref{f-estimate} in this case, we get
\begin{align}\label{S21-1}
|w_{j,\gamma}(M)g_{\varrho,j}(M)|& \leq |w_{j,\gamma}(M)|\int_{E_k}|f_{\varrho,j}(\tilde{y}u)|d\mathrm{Y_k}  \\
& \ll M^{\Re\gamma-r-\varrho}(\log M)^{\nu-1}\, (xM)^{\theta_\varrho-\frac{1}{2\a}}\notag \\
&= x^{\theta_\varrho-\frac{j}{2\a}}M^{\Re\gamma-\frac{r}{2}-\frac{1+2j}{4\a}-\frac{\varrho}{2\a}+\frac{\mu'-\mu}{2\a}}(\log M)^{\nu-1}. \notag
\end{align}

For the last integral of \eqref{S1} which corresponds to $\gamma$, \eqref{f-estimate} implies that
\begin{align}\label{S21-2}
\int_M^{M_1}w_{k+1,\gamma}'(u)g_{\varrho,k+1}(u)du & \ll \int_M^{M_1}u^{\Re\gamma-r-\varrho}(\log u)^{\nu-1}\,(xu)^{\theta_\varrho-\frac{k+1}{2\a}}du  \\
& = x^{\theta_\varrho-\frac{k+1}{2\a}}M^{\Re\gamma-\frac{r}{2}-\frac{1}{4\a}-\frac{\varrho}{2}+\frac{\mu'-\mu}{2\a}-\frac{k+1}{2\a}}(\log M)^{\nu-1}. \notag
\end{align}
From \eqref{S21-1} and \eqref{S21-2} we get the estimate from the pole  $\gamma$. It is easily seen that the term $j=1$ coming from $\gamma_0$
gives the maximal estimate, namely we have
\begin{equation*}
S_{21}(M) \ll x^{\theta_\varrho-\frac{1}{2\a}}M^{\Re\gamma_0-\frac{r}{2}-\frac{3}{4\a}-\frac{\varrho}{2\a}+\frac{\mu'-\mu}{2\a}}(\log M)^{A'} 
\end{equation*}
with some integer $A'$. Now we note that $\Re \gamma_0 \leq \sigma_b^{\ast}$, hence $S_{21}(M)$ is absorbed into the right hand side
of \eqref{V1-estimate}.

Next we consider $S_{22}(M)$.
By the trivial bound $\hat{E}(x) \ll x^{\sigma_b^{\ast}+\varepsilon}$ and \eqref{g-estimate}, we see that
\begin{align*}
\frac{\hat{E}(M_1)}{M_1^{r+\varrho}}\int_{\mathbf{E}_k}f_\varrho(\tilde{y}M_1)d\mathrm{Y_k}
& \ll M_1^{\sigma_b^{\ast}-r-\varrho+\varepsilon}x^{2k}(yM_1)^{\theta_\varrho-\frac{k}{2\a}}
\end{align*}
The exponent of $M_1$ is negative by \eqref{k-condition}, hence this term vanishes when $M_1 \to \infty$.

Consider the term $\frac{\hat{E}(M)}{M^{r+\varrho}}\int_{E_k}f_{\varrho}(\tilde{y}M)dY$.
If we assume that $b(n)\ge 0(n\ge 1)$, then from Theorem 6 (ii) we can find $M$ such that $\hat{E}(M)=0$, hence this term becomes zero.
If we assume the condition (B1), there exists $M$ such that $\hat{E}(M) \ll M^{\sigma_b^{\ast}-1+\omega_0}$.
Hence by \eqref{f-estimate}
\begin{align}
S_{22}(M) &\ll x^{\theta_\varrho}M^{\sigma_b^{\ast}-1+\omega_0-(r+\varrho)+\theta_\varrho} \notag \\
          &=   x^{\theta_\varrho}M^{\sigma_b^{\ast}-1+\omega_0-\frac{r}{2}-\frac{1}{4\a}-\frac{\varrho}{2\a}+\frac{\mu'-\mu}{2\a}}.
            \nonumber 
\end{align}
We call the above estimate of $S_{22}(M)$ as $R_6(M)$.

Finally we consider $S_{23}(M)$. Divide the integral as
\begin{align*}
S_{23}(M)=\left(\int_M^N+\int_N^\infty\right)\hat{E}(u)\int_{\mathbf{E}_k} \frac{d}{du}\left(\frac{f_{\varrho}(\tilde{y}u)}{u^{r+\varrho}}\right)d\mathrm{Y_k} du
=:S_{231}+S_{232}.
\end{align*}

To treat $S_{231}$, we substitute Taylor's formula for $\Phi_{\varrho}(\tilde{y})=\frac{d}{du}\left(\frac{f_{\varrho}(\tilde{y}u)}{u^{r+\varrho}}\right)$:
\begin{align}
\Phi_{\varrho}(\tilde{y})=\sum_{l=0}^{J}\left(\frac {d^l}{dy^l}\Phi_{\varrho}(y)\right)\frac {(\tilde{y}-y)^l}{l!}+
\left.\frac {d^{J+1}}{dv^{J+1}}\Phi_{\varrho}(v)\right|_{v=v_0}\frac {(\tilde{y}-y)^{J+1}}{(J+1)!}, \label{G-Taylor}
\end{align}
where $v_0=y+\theta(\tilde{y}-y)$ with some $0 \leq \theta \leq 1$.
Let $V_4$ be the contribution in $S_{231}$ from the last error term of \eqref{G-Taylor}. From \eqref{G-heikin} and $v_0 \sim y \sim x$, we have
\begin{align}
V_4&=\int_M^N \hat{E}(u)du \int_{\mathbf{E}_k} \left.\frac {d^{J+1}}{dv^{J+1}}\Phi_{\varrho}(v)\right|_{v=v_0}
\frac {(\tilde{y}-y)^{J+1}}{(J+1)!}d\mathrm{Y}_k \notag  \\
& \ll \int_M^N |\hat{E}(u)|\frac{(xu)^{\theta_\varrho+\frac{1}{2\a}+\frac {J+1}{2\a}}}{x^{J+1}u^{r+\varrho+1}}\frac{1}{x^{J+1}}du \notag \\
& \ll \int_M^N |\hat{E}(u)|\frac{(xu)^{\theta_\varrho-\frac{r}{2}}}{u^{r+\varrho+1}}(xu)^{\frac{r}{2}+\frac{1}{2\a}}
       \left(\frac{(xu)^{\frac{1}{2\a}}}{x^2}\right)^{J+1}du. \nonumber 
\end{align}
Since $u \leq N=[x^{4\a-1-\varepsilon}]$ and $J=[(4\a^2 r+4\a)\varepsilon^{-1}]$, it follows that
\begin{align*}
(ux)^{\frac r2 +\frac {1}{2\alpha}}\left(\frac {(ux)^{\frac {1}{2\alpha}}}{x^2}\right)^{J+1}
& \leq (x^{4\alpha-1-\varepsilon}x)^{\frac r2 +\frac {1}{2\alpha}}
\left(\frac {(x^{4\alpha-1-\varepsilon}x)^{\frac {1}{2\alpha}}}{x^2}\right)^{J+1} \\
&\leq x^{-\varepsilon(\frac r2+\frac{1}{2\a})} \leq 1.
\end{align*}
Therefore
\begin{align*}
V_4 &\ll x^{\theta_\varrho-\frac{r}{2}}\int_M^N|\hat{E}(u)|u^{-r-\frac{1}{4\a}-\frac{\varrho}{2\a}+\frac{\mu'-\mu}{2\a}-1}du.
\end{align*}
Now by the assumption (B3) we get the upper bound of $V_4$:
\begin{equation} \label{s231-3}
V_4 \ll x^{\theta_\varrho-\frac{r}{2}}M^{\sigma_b^{\ast}+\omega_1-r-\frac{1}{4a}-\frac{\varrho}{2\a}+\frac{\mu'-\mu}{2\a}-1}.
\end{equation}

Next we treat the first $J$-terms of the Taylor expansion of $\Phi_{\varrho}(\tilde{y})$. Their contribution to $S_{231}$ is
\begin{align}
&\sum_{l=0}^{J}\frac {x^{-l}}{l!}\int_M^N \hat{E}(u) \frac {d^l}{dy^l}\Phi_{\varrho}(y) du
 \int_{\mathbf{E}_k}(y_1+\cdots+y_k)^ld\mathrm{Y}_k\nonumber\\
&=\sum_{l=0}^{J}c_{l}^*x^{-l}\int_M^N\hat{E}(u) \frac {d^l}{dy^l}\Phi_{\varrho}(y)du \nonumber\\
&=\sum_{l=0}^{J}\sum_{m=0}^{J}\Re\left\{c_l^* b_{lm}'\int_M^N \hat{E}(u) \frac{(uy)^{\theta_{\varrho}+\frac{1}{2\a}+\frac{l-m}{2\a}}}{(xy)^lu^{r+\varrho+1}}
e^{-ih(uy)^{\frac{1}{2\a}}}du \right\}\nonumber\\
&\quad +O\left(\sum_{l=0}^{J} \int_M^N |\hat{E}(u)| \frac {(ux)^{\theta_\varrho+\frac{l-J}{2\a}}}{x^{2l} u^{r+\varrho+1}}du\right)
       +O\left(\sum_{l=0}^{J}x^{-l+r+\varrho-b-1}\int_M^N |\hat{E}(u)|u^{-1-b}du\right). \nonumber
\end{align}
The last $O$-term can be removed by taking $b$ large. In the first $O$-term, the integrand can be transformed into
\begin{align*}
|\hat{E}(u)|\frac{(ux)^{\theta_{\varrho}-\frac{J}{2\a}}}{u^{r+\varrho+1}}\left(\frac{(ux)^{\frac{1}{2\a}}}{x^2}\right)^l.
\end{align*}
Since $\frac{(ux)^{\frac{1}{2\a}}}{x^2} \leq 1$, this $O$-term is bounded by
$$
\ll \sum_{l=0}^{J}\int_M^N |\hat{E}(u)|\frac{(ux)^{\theta_{\varrho}-\frac{J}{2\a}}}{u^{r+\varrho+1}}du
\ll x^{\theta_{\varrho}-\frac{J}{2\a}}M^{\sigma_b^{\ast}+\omega_1-\frac{r}{2}-\frac{1}{4\a}-\frac{\varrho}{2\a}+\frac{\mu'-\mu}{2\a}-\frac{J}{2\a}-1},
$$
which is contained in the right hand side of \eqref{s231-3}.
Combining these formulas we obtain that
\begin{align}
S_{231}& = \sum_{l=0}^{J}\sum_{m=0}^{J}\Re\left\{c_l^* b_{lm}'\int_M^N \hat{E}(u) \frac{(uy)^{\theta_{\varrho}+\frac{1}{2\a}+\frac{l-m}{2\a}}}{(xy)^lu^{r+\varrho+1}}
e^{-ih(uy)^{\frac{1}{2\a}}}du \right\}\nonumber\\
& \quad +O\left(x^{\theta_\varrho-\frac{r}{2}}M^{\sigma_b^{\ast}+\omega_1-r-\frac{1}{4a}-\frac{\varrho}{2\a}+\frac{\mu'-\mu}{2\a}-1}\right). \nonumber
\end{align}
As for the double sum over $l$ and $m$, let $R_2(y;M)$ denote the  term corresponding to $l=m=0$
and $R_3(y;M)$ the term corresponding to $l+m >0$.
With the notation of the function $I(\lambda, M,N,y)$, we can express these terms as
\begin{align*}
R_2(y;M)&=y^{\theta_\varrho+\frac{1}{2\a}}\Re \left\{c_{00}I\left(\theta_\varrho+\frac{1}{2\a}-r-\varrho-1, M,N,y\right)\right\}, \\
R_3(y;M)&=\Dsum_{l+m \neq 0}\Re \left\{c_{lm} I\left(\theta_\varrho+\frac{1}{2\a}+\frac{l-m}{2\a}-r-\varrho-1,M,N,y\right)\right\}
x^{-l}y^{\theta_\varrho+\frac{1}{2\a}-l+\frac{l-m}{2\a}}.
\end{align*}

\bigskip

Finally we consider $S_{232}$. From \eqref{G-heikin} with $H=k$, we have
\begin{align*}
\int_{\mathbf{E}_k}\Phi_{\varrho}(\tilde{y})d\mathrm{Y_k}
&=\sum_{j=0}^{k} \left\{\sum_{m=0}^{k} \Re \left(b_{jm}''
\frac{x^k(y+\frac{j}{x})^{k+\theta_{\varrho}+\frac{1}{2\a}-\frac{k+m}{2\a}}}{u^{r+\varrho+1-\frac{1}{2\a}-\theta_{\varrho}+\frac{k+m}{2\a}}}
e^{-ih((y+\frac{j}{x})u)^{\frac{1}{2\a}}} \right) \right. \notag \\
& \qquad \left. +O\left(\frac{x^{2k+\theta_{\varrho}+\frac{1}{2\a}-\frac{2k+1}{2\a}}}
{u^{r+\varrho+1-\frac{1}{2\a}-\theta_{\varrho}+\frac{2k+1}{2\a}}}\right)\right\} +O\left(x^{2k+r+\varrho-b}u^{-1-b}\right).
\end{align*}
Substituting this expression to the definition of $S_{232}$ and noting \eqref{xi-1},
the contribution from the second $O$-term becomes
\begin{align*}
&\ll \int_N^{\infty}|\hat{E}(u)|\frac{x^{2k+\theta_{\varrho}+\frac{k}{\a}}}{u^{\frac{r}{2}+\frac{1}{4\a}+1+\frac{\varrho}{2\a}-\frac{\mu'-\mu}{2\a}+\frac{k}{\a}}}du\\
&\ll x^{2k+\theta_{\varrho}+\frac{k}{\a}}N^{\sigma_b^{\ast}+\omega_1-(\frac{r}{2}+\frac{1}{4\a}+1+\frac{\varrho}{2\a}-\frac{\mu'-\mu}{2\a}+\frac{k}{\a})}, \\
&\ll x^{(4\a-1)(\sigma_b^{\ast}+\omega_1)-2k+r-\varrho+2(\mu'-\mu)+\frac{2k}{\a}-2\a r -4\a}.
\end{align*}
The contribution from the last $O$-term is also negligible by taking $b$ large.
Combining these, we deduce that
\begin{align}\label{s232-errorterm}
S_{232}&=\sum_{j=0}^{k}\sum_{m=0}^{k}\Re \left\{c_{jm}'\int_N^{\infty}\hat{E}(u)
\frac{x^k(y+\frac{j}{x})^{k+\theta_{\varrho}+\frac{1}{2\a}-\frac{k+m}{2\a}}}{u^{r+\varrho+1-\frac{1}{2\a}-\theta_{\varrho}+\frac{k+m}{2\a}}}
e^{-ih((y+\frac{j}{x})u)^{\frac{1}{2\a}}}du\right\} \\
& \quad +O\left(x^{(4\a-1)(\sigma_b^{\ast}+\omega_1)-2k+r-\varrho+2(\mu'-\mu)+\frac{2k}{\a}-2\a r -4\a}\right). \notag
\end{align}
The term of sum over $j$ and $m$, which we denote by $R_4(y;M)$, can be expressed as
\begin{align*}
R_4(y;M)=&x^k \sum_{j=0}^k \sum_{m=0}^k \Re\left\{c_{lm}' \left(y+\frac{j}{x}\right)^{k+\theta_\varrho+\frac{1}{2\a}-\frac{k+m}{2\a}} \right. \\
& \times \left. I\left(-r-\varrho-1+\frac{1}{2\a}+\theta_\varrho-\frac{k+m}{2\a},N,\infty,y+\frac{j}{x}\right)\right\}
\end{align*}

We call the sum of remaining terms \eqref{V1-estimate}, \eqref{V3-estimate}, \eqref{s231-3} and \eqref{s232-errorterm} as $R_5(y;M)$.
This completes the proof of Theorem 7.   \qed


\section{Some estimates for the weighted integral of the error term}


In this section we shall give several estimates involving the exponential integral
$I(\lambda, M,N,y) $ defined by \eqref{I-definition}, which is closely related to the mean square bound $\eqref{psi-condition}$ and important in this paper.
For convenience we  use the notation $\int_{(\sigma)}f(s)ds=\int_{\sigma-i\infty}^{\sigma+i\infty}f(s)ds$.

We first give the following preliminary lemma.


\begin{lem} Let  $x\ge 1, \frac 12\le P\le P_1<P_2\le (1+\delta)P$.
 We have
\begin{align}\label{CDT}
\int_{P_1}^{P_2}u^{\lambda} e^{\pm i(t\log u-ux)}du
\ll x^{-\frac 12}P^{\lambda+\frac 12}(1+|t|)^{-\varepsilon}(Px)^{\varepsilon}.
\end{align}
\end{lem}

\proof Let $f(u)=t \log u -ux$. Clearly
$$
f'(u)=\frac tu -x,\quad f''(u)=-\frac {t}{u^2}.
$$

First we suppose $(1-\delta)Px\le t \le (1+2\delta)Px$, we apply the second mean value theorem and the second derivative test
to get
\begin{align}\label{SDT}
&\int_{P_1}^{P_2}u^{\lambda}e^{\pm i(t\log u-ux)}du
\ll P^{\lambda}\max_{P_1\le P_3\le P_4\le P_2}\left|\int_{P_3}^{P_4}e^{\pm i(t\log u-ux)}du\right|  \\
&\ll P^{\lambda}\max_{P_1\le u\le P_2}\sqrt{\frac {u^2}{t}}
\ll x^{-\frac 12}P^{\lambda+\frac 12}, \nonumber
\end{align}
from which the assertion of the lemma follows in this case.

If $(1+2\delta)Px< t < \infty$ or $-\infty<t<(1-\delta)Px$, by the first derivative test we get
\begin{align}\label{FDT}
&\int_{P_1}^{P_2}u^{\lambda}e^{\pm i(t\log u-ux)}du \ll  P^{\lambda}\max_{P_1\le P_3<P_4\le P_2}
\left\{\left|\frac {t}{P_3}-x\right|^{-1},\left|\frac {t}{P_4}-x\right|^{-1} \right\}   \\
&\ll P^{\lambda}x^{-1}(1+|t|)^{-\varepsilon}(Px)^{\varepsilon}
\ (\mbox{or}\ \ \ll P^{\lambda}x^{-1}).\nonumber
\end{align}
The assertion of the lemma can be checked easily by the above two cases.\qed

\begin{lem}
Let $M<N \le x^{A}$ ($A$ a fixed positive number), $w$ be a real number and $0<\mu_n\le \frac M2$. Then we have
\begin{align}\label{proposition-1}
\int_{x}^{(1+\delta)x}I(\lambda,M,N,y)y^w\cos\left(h(\mu_ny)^{\frac {1}{2\alpha}}+c\right)dy
\ll x^{w+1-\frac {3}{4\alpha}+\varepsilon '}\max_{M\le P\le N}P^{\lambda +\sigma^*+1-\frac {3}{4\alpha}}.
\end{align}
\end{lem}

\proof From \eqref{psi-condition} we easily get (by using Cauchy's inequality )
\begin{align}
\int_{-T}^{T}\left|\psi(\sigma^*+it)\right|^2dt &\ll T^{1+\varepsilon}, \label{mean-square} \\
\int_{-T}^{T}\left|\psi(\sigma^*+it)\right|dt &\ll \sqrt{2T\int_{-T}^{T}\left|\psi(\sigma^*+it)\right|^2dt}
\ll T^{1+\frac {\varepsilon}{2}}, \nonumber
\end{align}
and (by using integration by parts)
\begin{align} \label{convergence}
\int_{0}^{\infty}\frac{\left|\psi(\sigma^*\pm it)\right|}
{(1+|t|)^{1+\varepsilon}}dt=\left. \frac {\int_{0}^{t}\left|\psi(\sigma^*\pm iu)\right|du}
{(1+|t|)^{1+\varepsilon}}\right|_{t=0}^{\infty}+(1+\varepsilon)
\int_{0}^{\infty}\frac {\int_{0}^{t}\left|\psi(\sigma^*\pm iu)\right|du}
{(1+|t|)^{2+\varepsilon}}dt\ll 1.
\end{align}

By Perron's formula  and the residue theorem we have
\begin{align}
\hat{E}(u)&=\lim_{T\rightarrow \infty}\frac {1}{2\pi i}\int_{\sigma_b^*+\varepsilon-iT}^{\sigma_b^*+\varepsilon+iT}
\psi(s)\frac {u^s}{s}ds-\hat{Q}(u) \nonumber \\
&=\lim_{T\rightarrow \infty}\frac {1}{2\pi i}\int_{\sigma^*-iT}^{\sigma^*+iT}
\psi(s)\frac {u^s}{s}ds,  \nonumber 
\end{align}
since $\psi (s)s^{-1}\rightarrow 0$ uniformly in the strip $\sigma^*\le \Re s\le \sigma_b^* $ as $t\rightarrow\pm \infty$
(also see page 357 in Ivi\'c \cite{I1}).
Let $\frac 12\le P<Q\le (1+\delta)P$, then from \eqref{I-definition}
\begin{align}\label{I-expression}
I(\lambda,P,Q,y)&=\int_P^Qu^\lambda e^{-ih(uy)^{\frac {1}{2\alpha}}}du \lim_{T\rightarrow \infty}\int_{\sigma^*-iT}^{\sigma^*+iT}
\psi(s)\frac {u^s}{s}ds  \\
&=\lim_{T\rightarrow \infty}\int_{\sigma^*-iT}^{\sigma^*+iT}
\frac {\psi(s)}{s}ds\int_P^Qu^{\lambda+s}e^{-ih(uy)^{\frac {1}{2\alpha}}}du \nonumber \\
&=\int_{\sigma^*-i\infty}^{\sigma^*+i\infty}
\frac {\psi(s)}{s}ds\int_P^Qu^{\lambda+\sigma^*+it}e^{-ih(uy)^{\frac {1}{2\alpha}}}du, \quad (s=\sigma^*+it), \nonumber
\end{align}
where the inversion of order of two integrations is justified
by \eqref{CDT} and \eqref{convergence}.

Assume $0<\mu_n\le \frac P2$, we get
\begin{align}\label{I-integral}
 &\int_{x}^{(1+\delta)x}I(\lambda,P,Q,y)y^w e^{\pm h(\mu_ny)^{\frac {1}{2\alpha}}}dy  \\
 &=\int_{x}^{(1+\delta)x}y^w e^{\pm ih(\mu_ny)^{\frac {1}{2\alpha}}}dy\int_{\sigma^*-i\infty}^{\sigma^*+i\infty}
 \frac {\psi(s)}{s}ds\int_P^Qu^{\lambda+\sigma^*+it}e^{-ih(uy)^{\frac {1}{2\alpha}}}du\nonumber\\
 &=\int_{\sigma^*-i\infty}^{\sigma^*+i\infty}\frac {\psi(s)}{s}ds\int_P^Q\int_{x}^{(1+\delta)x}
   y^w u^{\lambda+\sigma^*+it} e^{i(\pm h(\mu_ny)^{\frac {1}{2\alpha}}-h(uy)^{\frac {1}{2\alpha}})}dydu \nonumber\\
 &=:\int_{\sigma^*-i\infty}^{\sigma^*+i\infty} \frac {\psi(s)}{s}I_1(s)ds, \nonumber
 \end{align}
 where the inversion of order of the integrations is justified by \eqref{CDT} and \eqref{convergence}.
Set $u^{\frac {1}{2\alpha}}=U, U_0=P^{\frac {1}{2\alpha}},
U_1=Q^{\frac {1}{2\alpha}}$, $hy^{\frac {1}{2\alpha}}=Y, X_0=hx^{\frac {1}{2\alpha}}$ and $X_1=h\left((1+\delta)x\right)^{\frac {1}{2\alpha}}$, hence
\begin{align}\label{I1-definition}
I_1(s)=\frac {(2\alpha)^2}{h^{2\alpha(w+1)}}\int_{U_0}^{U_1}
U^{2\alpha (\lambda +1+\sigma^*+it)-1}dU \int_{X_0}^{X_1}
Y^{2\alpha(w+1) -1}e^{i (\pm \mu_n^{\frac {1}{2\alpha}}-U )Y }dY.
\end{align}

We write $\eta=2\alpha(w+1) -1$. Applying integration by parts we get
\begin{align*}
&\int_{X_0}^{X_1}
Y^{\eta}e^{i(\pm \mu_n^{\frac {1}{2\alpha}}-U)Y}dY\\
&=\frac {1}{(\pm\mu_n^{\frac {1}{2\alpha}}-U)i}\left(\left.Y^{\eta }e^{i (\pm \mu_n^{\frac {1}{2\alpha}}-U)Y} \right|_{X_0}^{X_1}
-\eta\int_{X_0}^{X_1}Y^{\eta-1}e^{i(\pm \mu_n^{\frac {1}{2\alpha}}-U)Y}dY \right).
\end{align*}
Substituting the above expression into \eqref{I1-definition} we obtain
\begin{align*}
I_1(s)\ll X_0^{\eta}\max_{X_0\le Y\le X_1}\left|
\int_{U_0}^{U_1}
\frac {U^{2\alpha (\lambda +1+\sigma^*)-1}}{\pm \mu_n^{\frac {1}{2\alpha}}-U} e^{i(2\alpha t\log U-UY)}dU
\right|.
\end{align*}
Now we apply the second mean value theorem and \eqref{CDT} to get
\begin{align*}
I_1(s)&\ll \frac {X_0^{\eta}}{U_0\mp\mu_n^{\frac {1}{2\alpha}}}\max_{X_0\le Y\le X_1, U_0\le U^*\le U_1}\left|
\int_{U_0}^{U^*}
U^{2\alpha (\lambda+1+\sigma^*)-1} e^{i(2\alpha t \log U-UY)}dU
\right|\\
&\ll \frac {X_0^{\eta}}{U_0}X_0^{-\frac 12}
U_0^{2\alpha(\lambda +1+\sigma^*)-\frac 12}(X_0U_0)^{\varepsilon}(1+|t|)^{-\varepsilon}\\
&\ll x^{w+1-\frac {3}{4\alpha}}P^{\lambda +1+\sigma^*-\frac {3}{4\alpha}}(Px)^{\varepsilon}(1+|t|)^{-\varepsilon}.
\end{align*}
We substitute the above estimate into \eqref{I-integral} and then use \eqref{convergence} to get
\begin{align}
 &\int_{x}^{(1+\delta)x}I(\lambda,P,Q,y)y^{w} e^{\pm ih(\mu_ny)^{\frac {1}{2\alpha}}}dy=\int_{\sigma^*-i\infty}^{\sigma^*+i\infty}
\frac {\psi(s)}{s}I_1(s)ds  \nonumber \\
&\qquad \ll x^{w+1-\frac {3}{4\alpha}}P^{\lambda +1+\sigma^*-\frac {3}{4\alpha}}(Px)^{\varepsilon}\int_{-\infty}^{\infty}
\frac {|\psi(\sigma^*+it)|}{(1+|t|)^{1+\varepsilon}}dt \nonumber \\
&\qquad \ll x^{w+1-\frac {3}{4\alpha}}P^{\lambda +1+\sigma^*-\frac {3}{4\alpha}}(Px)^{\varepsilon}. \nonumber 
\end{align}
Now applying a simple splitting argument to the interval $[M,N]$ would suffice to complete the proof of the lemma.\qed

\medskip

\begin{lem}
Let $2(\lambda+\sigma^*)\ne -1$, $M<N\le x^{A} (A \mbox{ a fixed positive number})$ and  $\delta>0$ with
 $(1+\delta)^{\frac {1}{\alpha}}-1<\frac 14$. Then we have
\begin{align}\label{proposition-2}
\int_{x}^{(1+\delta)x}\left|I(\lambda,M,N,y)\right|^2dy\ll x^{1-\frac {1}{\alpha}+\varepsilon '}\max_{M\le P\le N}P^{2(\lambda +\sigma^*+1)-\frac {1}{\alpha}}.
\end{align}
\end{lem}

\proof From \eqref{I-expression} we obtain
\begin{align}\label{IPQ-expression}
&\int_{x}^{(1+\delta)x}\left|I(\lambda,P,Q,y)\right|^2dy  \\
&=\int_{x}^{(1+\delta)x}dy\left(\int_{(\sigma^*)}\frac {\psi(s)}{s}ds
\int_P^Q u^{\lambda+\sigma^*} e^{i\bigl(t\log u- h(uy)^{\frac {1}{2\alpha}}\bigr)}du\right) \nonumber \\
& \quad\ \times \left(\int_{(\sigma^*)}\frac {\psi(\bar{s}_1)}{\bar{s}_1}d\bar{s}_1
\int_P^Q v^{\lambda+\sigma^*} e^{i\bigl(-t_1\log v+ h(vy)^{\frac {1}{2\alpha}}\bigr)}dv\right)\nonumber\\
&=\int_{(\sigma^*)}\int_{(\sigma^*)}\frac {\psi(s)\psi (\bar {s}_1)}{s\bar{s}_1}dsd\bar{s}_1\times I(s,s_1), \nonumber
\end{align}
where the inversion of order of the integrations is justified by \eqref{CDT} and \eqref{convergence},
and we use notation $s=\sigma^* +it$, $s_1=\sigma^* +it_1$ and
\begin{align*}
I(s,s_1):=\int_{x}^{(1+\delta)x}dy
\int_P^Q u^{\lambda+\sigma^*} e^{i\bigl(t\log u- h(uy)^{\frac {1}{2\alpha}}\bigr)}du
\int_P^Q v^{\lambda+\sigma^*} e^{i\bigl(-t_1\log v+ h(vy)^{\frac {1}{2\alpha}}\bigr)}dv.
\end{align*}

Set $h(uy)^{\frac {1}{2\alpha}}=U, h(vy)^{\frac {1}{2\alpha}}=V,
P^{\frac {1}{2\alpha}}=U_0, Q^{\frac {1}{2\alpha}}=U_1,
hy^{\frac {1}{2\alpha}}=Y,hx^{\frac {1}{2\alpha}}=X_0,
h((1+\delta)x)^{\frac {1}{2\alpha}}=X_1
$, $t_2=2\alpha t, t_3=2\alpha t_1$,  $\mu=2\alpha(\lambda +\sigma^*+1)-1$ and $\beta=-4\alpha(\lambda+\sigma^*) -2\alpha -1$. We have
\begin{align}\label{Iss1-D}
I(s,s_1)=&\frac {(2\alpha)^3}{h^{2\alpha}}\int_{X_0}^{X_1}
Y^{\beta}e^{i(t_3-t_2)\log Y}G(Y;t_2,t_3)dY,
\end{align}
where
\begin{align*}
G(Y;t_2,t_3):=
&\int_{U_0Y}^{U_1Y}U^\mu e^{i(t_2\log U-U)}dU \int_{U_0Y}^{U_1Y}V^\mu e^{i(- t_3\log V+V)}dV.
\end{align*}

Since $P<Q\le (1+\delta)P$, so $U_0X_0<U_1X_1\le (1+\delta)
^{\frac {1}{\alpha}}U_0X_0:=(1+\delta_1)U_0X_0$. Here note that $0<\delta_1=(1+\delta)
^{\frac {1}{\alpha}}-1<\frac 14.$

To estimate $I(s,s_1)$ we need to consider the following four cases.

\noindent{\bf (Case 1).}
 $t_2, t_3\notin \left((1-2\delta_1)U_0X_0,(1+4\delta_1)U_0X_0\right)$. In this case
applying \eqref{FDT} to estimate the two integrals in $G(Y;t_2,t_3)$ we obtain
\begin{align*}
G(Y;t_2,t_3)&\ll \max_{X_0\le Y\le X_1}\left|\int_{U_0Y}^{U_1Y}U^\mu e^{i( t_2\log U-U)}dU\right|\times
\max_{X_0\le Y\le X_1}\left|\int_{U_0Y}^{U_1Y}V^\mu e^{i(- t_3\log V+V)}dV \right| \\
&\ll (U_0X_0)^{2\mu}(U_0X_0)^{2\varepsilon}(1+|t_2|)^{-\varepsilon} (1+|t_3|)^{-\varepsilon}.
\end{align*}
We substitute the above estimate into \eqref{Iss1-D} to get
\begin{align}\label{case1-bound}
I(s,s_1)\ll U_0^{2\mu}X_0^{2\mu+\beta+1}(U_0X_0)^{2\varepsilon}(1+|t_2|)
^{-\varepsilon}(1+|t_3|)^{-\varepsilon}.
\end{align}

\noindent{\bf (Case 2).}
$t_2\in \left((1-2\delta_1)U_0X_0,(1+4\delta_1)U_0X_0\right)$, $t_3\notin \left((1-2\delta_1)U_0X_0,(1+4\delta_1)U_0X_0\right)$.
If $|t_2-t_3|\le \delta_1U_0X_0$, then $t_2<(1-\delta_1)U_0X_0$ or $t_2>(1+3\delta_1)U_0X_0$. Clearly the estimate \eqref{case1-bound} also holds.

 Now we assume $|t_2-t_3|>\delta_1U_0X_0$.
 Applying integration by parts to the variable $Y$ in \eqref{Iss1-D} we get
\begin{align}\label{Iss1-expression-1}
&\frac {h^{2\alpha}}{(2\alpha)^3}I(s,s_1)=\int_{X_0}^{X_1}Y^{\beta}e^{i(t_3-t_2)\log Y}G(Y;t_2,t_3)dY  \\
&=\frac {1}{\beta+1+(t_3-t_2)i}\int_{X_0}^{X_1}G(Y;t_2,t_3)d Y^{\beta+1+ i(t_3-t_2)}\nonumber\\
&=\frac{1}{\beta+1+(t_3-t_2)i}\left(\sum_{j=0}^{1}(-1)^{j+1}X_j^{\beta+1} e^{i(t_3-t_2)\log X_j} G(X_j;t_2,t_3)\right.\nonumber\\
&\ \left. +\sum_{j=0}^{1}(-1)^{j}U_j^{\mu+1+it_2}I_{1j} +\sum_{j=0}^{1}(-1)^{j}U_j^{\mu+1-it_3}I_{2j}\right),\nonumber
\end{align}
where
\begin{align*}
I_{1j}:=\int_{X_0}^{X_1}Y^{\beta+\mu+1} e^{i( t_3\log Y-U_jY)}dY
\int_{U_0Y}^{U_1Y}V^\mu e^{i(- t_3\log V+V)}dV\\
I_{2j}:=\int_{X_0}^{X_1}Y^{\beta+\mu+1} e^{i(- t_2\log Y+U_jY)}dY
\int_{U_0Y}^{U_1Y}U^\mu e^{i( t_2\log U-U)}dU.
\end{align*}

Applying \eqref{SDT} to the first integral over $U$ in $G(X_j;t_2,t_3)$, and \eqref{FDT}
to the second integral over $V$ in $G(X_j;t_2,t_3)$ respectively, we obtain
\begin{align}\label{GXj-bound}
G(X_j;t_2,t_3)\ll (U_0X_0)^{\mu+\frac 12}(U_0X_0)^{\mu}
\ll (U_0X_0)^{2\mu+\frac 12}, \quad (j=0,1).
\end{align}
Applying \eqref{FDT} to $I_{1j}$ one has
\begin{align}\label{I1j-bound}
I_{1j}&\ll X_0^{\beta+\mu+2}\max_{X_0\le Y\le X_1}\left|\int_{U_0Y}^{U_1Y}V^\mu e^{i(- t_3\log V+V)}dV\right|  \\
&\ll X_0^{\beta+\mu+2}(U_0X_0)^\mu\ll U_0^\mu X_0^{\beta+2\mu+2}.\notag
\end{align}

To treat $I_{2j}$ we define for $j=0,1$
$$X_{2j}=\frac {t_2}{U_j}-\sqrt{\frac {X_0}{U_j}}\quad X_{3j}=\frac {t_2}{U_j}+\sqrt{\frac {X_0}{U_j}}.$$

First we consider the case when $X_0\le X_{2j}<X_{3j}\le X_1$, we divide the interval $[X_0,X_1]$ into three subintervals and write
\begin{align}\label{I2j-decom}
I_{2j}=\int_{X_0}^{X_{2j}}+\int_{X_{2j}}^{X_{3j}}+\int_{X_{3j}}^{X_1}
:=I_{3j}+I_{4j}+I_{5j}.
\end{align}

We let $\kappa =\beta +\mu+1$ and  $H(Y,t_2):=\int_{U_0Y}^{U_1Y}U^\mu e^{i( t_2\log U-U)}dU$. It follows that from integration by parts
\begin{align}\label{I3j-decom}
I_{3j}&=\int_{X_0}^{X_{2j}}Y^{\kappa}H(Y,t_2) e^{i(- t_2\log Y+U_jY)}dY \\
&=\int_{X_0}^{X_{2j}}\frac {Y^{\kappa+1}}{i(-t_2+U_jY)}H(Y,t_2) de^{i(- t_2\log Y+U_jY)}  \nonumber\\[1ex]
&=I_{6j}+I_{7j}-I_{8j}, \nonumber
\end{align}
where\footnote{In Tong's paper\cite{To2}, page 521,  he missed the term $I_{8j}$.}
\begin{align*}
I_{6j}:&=\left. \frac {Y^{\kappa+1}}{i(-t_2+U_jY)}H(Y,t_2) e^{i(- t_2\log Y+U_jY)}\right|_{X_0}^{X_{2j}},\\
I_{7j}:&=\sum_{l=0}^{1}(-1)^{l}U_l^{\mu+1+it_2}\int_{X_0}^{X_{2j}}
\frac {Y^{\kappa+\mu+1}}{i(t_2-U_jY)}e^{i(U_j-U_l)Y}dY,\\
I_{8j}:&=\int_{X_0}^{X_{2j}}\left(\frac {(\kappa+1)Y^{\kappa}}{i(-t_2+U_jY)}-\frac
{Y^{\kappa+1}U_j}{i(-t_2+U_jY)^2}\right)H(Y,t_2) e^{i(- t_2\log Y+U_jY)}dY.
\end{align*}
By \eqref{SDT} we get
\begin{align*}
I_{6j}&\ll X_0^{\kappa+1}\max_{X_0\le Y\le X_{2j}}\frac {1}{|t_2-U_jY|}\times (U_0X_0)^{\mu+\frac 12}  \\
&\ll X_0^{\kappa+1}\frac {1}{\sqrt{X_0U_0}} (U_0X_0)^{\mu+\frac 12}\ll U_0^{\mu}X_0^{\kappa+\mu+1}.
\end{align*}
We also have
\begin{align*}
I_{7j}\ll U_0^{\mu+1}\int_{X_0}^{X_{2j}}\frac {X_0^{\kappa+\mu+1}}{t_2-U_jY}dY\ll U_0^{\mu}X_0^{\kappa+\mu+1}\log (U_0X_0).
\end{align*}
Applying \eqref{SDT} to the integral $H(Y,t_2)$ in $I_{8j}$ we obtain
\begin{align*}
I_{8j}&\ll \int_{X_0}^{X_{2j}}\left(\frac {X_0^{\kappa}}{t_2-U_jY}+\frac
{X_0^{\kappa+1}U_0}{(t_2-U_jY)^2}\right)(U_0X_0)^{\mu+\frac 12}dY  \\
&\ll U_0^{\mu+\frac 12}X_0^{\kappa+\mu+\frac 12}\int_{X_0}^{X_{2j}}\frac {1}{t_2-U_jY}dY+
U_0^{\mu+\frac 32}X_0^{\kappa+\mu+\frac 32}\int_{X_0}^{X_{2j}}\frac {1}{(t_2-U_jY)^2}dY \\
&\ll U_0^{\mu-\frac 12}X_0^{\kappa+\mu+\frac 12}\log (U_0X_0)+
 U_0^{\mu}X_0^{\kappa+\mu+1}\ll  U_0^{\mu}X_0^{\kappa+\mu+1}.
\end{align*}
Combining the above three estimates and \eqref{I3j-decom} we have
 \begin{align}\label{I3j-bound}
 I_{3j}\ll U_0^{\mu}X_0^{\kappa+\mu+1}\log (U_0X_0).
 \end{align}
Similarly to the estimate for $I_{8j}$ one has
 \begin{align}\label{I4j-bound}
 I_{4j}&\ll \int_{X_{2j}}^{X_{3j}}X_0^{\kappa}(U_0X_0)^{\mu+\frac 12}dY\ll (X_{3j}-X_{2j})U_0^{\mu+\frac 12}X_0^{\kappa+\mu+\frac 12}  \\
 &\ll \sqrt {\frac {X_0}{U_0}}U_0^{\mu+\frac 12}X_0^{\kappa+\mu+\frac 12}\ll
 U_0^{\mu}X_0^{\kappa+\mu+1}. \nonumber
\end{align}
 Similarly to the estimate for $I_{3j}$ we can prove
  \begin{align}\label{I5j-bound}
 I_{5j}\ll U_0^{\mu}X_0^{\kappa+\mu+1}\log (U_0X_0).
 \end{align}
 Substituting \eqref{I3j-bound}, \eqref{I4j-bound} and \eqref{I5j-bound} into \eqref{I2j-decom} we obtain when $X_0\le X_{2j}<X_{3j}\le X_1$
\begin{align*} 
 I_{2j}\ll U_0^{\mu}X_0^{\kappa+\mu+1}\log (U_0X_0).
 \end{align*}
It remains to consider the other four cases: $ X_{2j}<X_{3j}\le X_0< X_1$, $ X_{2j}\le X_0<X_{3j}\le X_1$,  $X_0\le X_{2j}\le X_1<X_{3j}$
and $X_0\le X_1\le X_{2j}<X_{3j}$.  By a very similar approach to the case $X_0\le X_{2j}<X_{3j}\le X_1$,
one may get the same estimate of $I_{2j}$ for the above four cases. For example, when $\le X_{2j}\le X_0<X_{3j}\le X_1$ we use decomposition
$$
I_{2j}=\int_{X_0}^{X_1}=\int_{X_0}^{X_{3j}}+\int_{X_{3j}}^{X_1}.
$$
Hence we always have
\begin{align}\label{I2j-bound-1}
 I_{2j}\ll U_0^{\mu}X_0^{\beta +2\mu+2}\log (U_0X_0).
\end{align}
Since $|t_2-t_3|>\delta_1U_0X_0$ and $t_2\in \left((1-2\delta_1)U_0X_0,(1+4\delta_1)U_0X_0\right)$, it is easy to check
\begin{equation*}
|t_2-t_3|>(U_0X_0)^{1-2\varepsilon}(1+|t_2|)^{-\varepsilon}
(1+|t_3|)^{-\varepsilon}.
\end{equation*}
 Substituting \eqref{GXj-bound}, \eqref{I1j-bound} and \eqref{I2j-bound-1} into \eqref{Iss1-expression-1},
and noting the above estimate, in (Case 2) we finally obtain
\begin{align} \label{case2-bound}
 I(s,s_1)\ll  U_0^{2\mu}X_0^{2\mu+\beta+1}(U_0X_0)^{3\varepsilon}(1+|t_2|)^{-\varepsilon}(1+|t_3|)^{-\varepsilon}.
\end{align}

\noindent{\bf (Case 3).}
$t_2\notin \left((1-2\delta_1)U_0X_0,(1+4\delta_1)U_0X_0\right)$, $t_3\in \left((1-2\delta_1)U_0X_0,(1+4\delta_1)U_0X_0\right)$.

We can use a similar approach as (Case 2) to estimate $I(s,s_1)$
and get the same estimate as \eqref{case2-bound}.

\noindent{\bf (Case 4).}
$t_2, t_3\in \left((1-2\delta_1)U_0X_0,(1+4\delta_1)U_0X_0\right)$.

Note that $\beta+1+ i(t_3-t_2)\neq 0$. We first transform $I(s,s_1)$ into the expression \eqref{Iss1-expression-1}.
Then we use the same approach as $I_{2j}$ in (Case 2) to treat $I_{2j}$ and estimate $I_{1j}$ respectively.
 Lastly by the second derivative test \eqref{SDT} to each integral in $G(X_j;t_2,t_3)$, we get
\begin{align*}
G(X_j;t_2,t_3)\ll (U_0X_0)^{\mu+\frac 12}(U_0X_0)^{\mu+\frac 12}
\ll (U_0X_0)^{2\mu+1}, \quad (j=0,1).
\end{align*}
Combining these estimates and  \eqref{Iss1-expression-1} we obtain
\begin{align*} 
I(s,s_1)\ll \frac {U_0^{2\mu+1}X_0^{\beta+2\mu+2}}{1+|t_2-t_3|}(U_0X_0)^{\varepsilon}
\end{align*}

Now let $h(Px)^{\frac {1}{2\alpha}}=T$. Note that $2\mu +\beta +1=2\alpha -2$. Combining four estimates
in the above cases and \eqref{IPQ-expression} we get
\begin{align}\label{IPQ-bound-1}
&\int_{x}^{(1+\delta)x}\left|I(\lambda,P,Q,y)\right|^2dy
=\int_{(\sigma^*)}\int_{(\sigma^*)}
\frac {\psi(s)\psi (\bar {s}_1)}{s\bar{s}_1}I(s,s_1) dsd\bar{s}_1  \\
&\ll x^{1-\frac {1}{\alpha}}P^{\frac {\mu}{\alpha}}T^{3\varepsilon}\int_{-\infty}^{\infty}
\int_{-\infty}^{\infty}\frac {|\psi (\sigma^*+it)\psi (\sigma^*-it_1)|}{(1+|t|)(1+|t_1|)}
\frac {dtdt_1}{(1+|t|)^{\varepsilon}(1+|t_1|)^{\varepsilon}}\nonumber\\
&\quad \quad+ x^{1-\frac {1}{2\alpha}}P^{\frac {2\mu+1}{2\alpha}}T^{\varepsilon}\int_{(1-2\delta_1)T}^{(1+4\delta_1)T}
\int_{(1-2\delta_1)T}^{(1+4\delta_1)T}
\frac {|\psi (\sigma^*+it) \psi (\sigma^*-it_1)|}{tt_1}\frac {dtdt_1}{1+|t-t_1|}\nonumber\\
&=:x^{1-\frac {1}{\alpha}}P^{\frac {\mu}{\alpha}}T^{3\varepsilon}W_1+x^{1-\frac {1}{2\alpha}}P^{\frac {2\mu+1}{2\alpha}}T^{\varepsilon}W_2.\nonumber
\end{align}
It follows from \eqref{convergence} that
\begin{align}\label{W01}
W_1=\int_{-\infty}^{\infty}\frac {|\psi (\sigma^*+it)|}{(1+|t|)^{1+\varepsilon}}dt
\times\int_{-\infty}^{\infty}
\frac {|\psi (\sigma^*-it_1)|}{(1+|t_1|)^{1+\varepsilon}}dt_1\ll 1.
\end{align}
Making the change of variables $t=u+v, t_1=u-v$,  using Cauchy's inequality and noting \eqref{mean-square}, it is easy to check
\begin{align}\label{W02}
W_2&\ll T^{-2}\int_{0}^{2T} \int_{0}^{2T} \frac {|\psi (\sigma^*+it)\psi (\sigma^*-it_1)|}{1+|t-t_1|}dtdt_1  \\
&\ll T^{-2}\int_{-T}^{T}\frac {1}{1+|v|}
\left(\int_{0}^{2T}|\psi(\sigma^*+i(u+v))
\psi(\sigma^*-i(u-v))|du\right)dv \nonumber\\
&\ll T^{-2}\int_{-T}^{T}\frac {1}{1+|v|}
\left(\int_{0}^{2T}|\psi(\sigma^*+i(u+v))|^2du \int_{0}^{2T}|
\psi(\sigma^*-i(u-v))|^2du\right)^{\frac 12} dv \nonumber \\
&\ll T^{-1+\varepsilon}\int_{-T}^{T}\frac {1}{1+|v|}dv\ll T^{-1+2\varepsilon}.\nonumber
\end{align}
Now substituting \eqref{W01} and \eqref{W02}
into \eqref{IPQ-bound-1} we get with $\varepsilon_1=\frac {3\varepsilon}{2\alpha}$
\begin{align}\label{IPQ-bound-2}
\int_{x}^{(1+\delta)x}\left|I(\lambda,P,Q,y)\right|^2dy
\ll x^{1-\frac {1}{\alpha}+\varepsilon_1}P^{\frac {\mu}{\alpha}+\varepsilon_1}.
\end{align}

Applying a simple splitting argument to the interval $[M,N]$ with
$M_0=M, M_r=N, M_{j+1}\le (1+\delta )M_{j} \ (j=0,\ldots, r-1; r\ll \log x)$ and Cauchy's inequality, it follows that from \eqref{IPQ-bound-2}
\begin{align}\label{IMN-bound}
&\int_{x}^{(1+\delta)x}\left|I(\lambda,M,N,y)\right|^2dy\ll
\int_{x}^{(1+\delta)x}\left(\sum_{j=0}^{r-1}
\left|I(\lambda,M_j,M_{j+1},y)\right|\right)^2dy  \\
&=\sum_{j=0}^{r-1}\sum_{l=0}^{r-1}\int_{x}^{(1+\delta)x}
\left|I(\lambda,M_j,M_{j+1},y)\right|
\left|I(\lambda,M_l,M_{l+1},y)\right|dy\nonumber\\
&\ll  \sum_{j=0}^{r-1}\sum_{l=0}^{r-1}
\sqrt{\int_{x}^{(1+\delta)x}\left|I(\lambda,M_j,M_{j+1},y)\right|^2dy
\times \int_{x}^{(1+\delta)x}\left|I(\lambda,M_l,M_{l+1},y)\right|^2dy}
\nonumber\\
&=\left(\sum_{j=0}^{r-1}\sqrt{\int_{x}^{(1+\delta)x}
\left|I(\lambda,M_j,M_{j+1},y)\right|^2dy}\right)^2
\ll x^{1-\frac {1}{\alpha}+\varepsilon_1}\left(\sum_{j=0}^{r-1}M_j^{\frac {\mu}{2\alpha}+\frac {\varepsilon_1}{2}}\right)^2.\nonumber
\end{align}
Recall $\mu=2\alpha(\lambda +\sigma^*+1)-1$. Now Lemma 9 follows from \eqref{IMN-bound} at once.\qed

\begin{lem}
Let $2(\lambda+\sigma^*)\ne -1$, $2(\lambda +\sigma^*+1)<\frac {1}{\alpha}$,  $M\ge 1$ ,  and  $\delta>0$ with
 $(1+\delta)^{\frac {1}{\alpha}}-1<\frac 14$. Then we have
\begin{align}\label{proposition-3}
\int_{x}^{(1+\delta)x}\left|I(\lambda,M,\infty,y)\right|^2dy\ll x^{1-\frac {1}{\alpha}+\varepsilon '}M^{2(\lambda +\sigma^*+1)-\frac {1}{\alpha}}.
\end{align}
\end{lem}

\proof
We take $M_j=M(1+\delta_1)^j,\  (j=0,1,\cdots,r-1)$. The lemma is immediately obtained from \eqref{IMN-bound} and let $r\rightarrow \infty$.
\qed

\begin{rem}
In Tong's original proof (similar to Lemma 9) he used the second mean value theorem for complex integral, but we don't.
Moreover, in both Lemma 9 and Lemma 10 the condition $2(\lambda+\sigma^*)\ne -1$ is not essential.
In fact, we can use the second mean value theorem to remove this restriction.
\end{rem}


\section{\bf The proof of Theorem 1 }


In this section we shall prove Theorem 1.

Suppose that ${\cal L}(s)\in {\cal S}_{real}^{\theta}$ is a function of degree
$d\ge 2$ such that $0\le \theta\le 1/2-1/2d.$ So  ${\cal L}(s)$ satisfies the functional equation \eqref{functional-equation} in Section 1.
We assume $\sum_{j=1}^{L} \b_j$ is real where $\b_j$ are the constants in \eqref{gamma-factors}.
From \eqref{def-degree} and \eqref{def-mu}--\eqref{def-thetarho}, we have obviously
$d=2\a, \mu=\mu', \nu=\nu', \lambda=\lambda', \theta_{\varrho}=1/2-1/2d+\varrho(1-1/d)$.
We also have $\varphi(s)={\cal L}(s)$ and $\psi(s)=\omega Q^{2s-1}{\cal L}(s)$, hence
$b(n)=\omega Q^{-1}a(n)$, $\lambda_n=n$ and $\mu_n=Q^{-2}n$.
Suppose that $1/2\le \sigma^{*}<1$ satisfies \eqref{meanL} and \eqref{sigmastarcondition}.

Suppose $T$ is a large parameter and $\delta>0$ is a small  positive constant. We shall evaluate the integral $\int_T^{(1+\delta)T}E^2(y)dy.$


\subsection{Tong's formula of $E(y)$}\


We first show that
\begin{equation}
\int_1^T|E(y)|dy\ll T^{1+\frac{d^2-1}{2d^2}+\varepsilon}. \label{Eq8-1}
\end{equation}

 Since ${\cal L}(s)$ is a function in $  {\cal S}_{real}^{\theta}$ of degree $d$, we have trivially that
${\cal L}(it)\ll (1+|t|)^{d/2+\varepsilon/2}$ and thus
\begin{equation}
\int_0^T|{\cal L}(it)|^2dt\ll T^{d+1+\varepsilon}.  \label{Eq8-2}
\end{equation}

From \eqref{meanL}, \eqref{Eq8-2} and Lemma 8.3 of Ivi\'c \cite{I1} we get for $0\le\sigma\le\sigma^{*}$
\begin{equation} \label{Eq8-3}
\int_0^T |{\cal L}(\sigma+it)|^2dt\ll T^{\frac{(d+1)(\sigma^{*}-\sigma)}{\sigma^{*}}+\frac{\sigma}{\sigma^{*}}+\varepsilon}.
\end{equation}

Define $\hat{E}(y)=A(y)-\Res_{s=1}{\cal L}(s)y^s s^{-1}.$ Then
$$E(y)= \hat{E}(y)-{\cal L}(0).$$
By Perron's formula we have for some $\sigma^{*}<\sigma<1$ that
$$ \hat{E}(y)=\frac{1}{2\pi i}\int_{\sigma-i\infty}^{\sigma+i\infty}{\cal L}(s)y^s s^{-1}ds.$$

From \eqref{Eq8-3} and integration by parts we see that
$$\int_{-\infty}^{\infty}|{\cal L}(\sigma_0+it)|^2|\sigma_0+it|^{-2}dt\ll 1$$
for $\sigma_0=\frac{(d-1)\sigma^{*}+2\sigma^{*}\varepsilon}{d}$.  So we have
$$
\hat{E}(y)=\frac{1}{2\pi i}\int_{\sigma_0-i\infty}^{\sigma_0+i\infty}{\cal L}(s)y^s s^{-1}ds.
$$

Replacing $y$ in the above formula by $1/y$, and then using Parseval's identity (see the formula (A.5) of Ivi\'c \cite{I1})
we have
$$
\frac{1}{2\pi}\int_{-\infty}^{\infty}|{\cal L}(\sigma_0+it)|^2|\sigma_0+it|^{-2}dt
=\int_0^\infty | \hat{E}(1/y)|^2 y^{ 2\sigma_0-1}dy=\int_0^\infty | \hat{E}(y)|^2y^{-2\sigma_0-1}dy.
$$
Thus we get for any $D>1$ that
$$\int_D^{2D}| \hat{E}(y)|^2dy\ll D^{2\sigma_0+1},$$
which with Cauchy's inequality gives
$$\int_D^{2D}| \hat{E}(y)| dy\ll D^{ \sigma_0+1}.$$
By the splitting argument we get
$$\int_0^{T}| \hat{E}(y)| dy\ll T^{ \sigma_0+1}.$$
Recalling   $E(y)= \hat{E}(y)-{\cal L}(0)$ and (1.5) we get
$$\int_0^{T}|E(y)| dy\ll \int_0^{T}|\hat{E}(y) | dy+T\ll T^{ \sigma_0+1}\ll T^{1+\frac{d^2-1}{2d^2}+\varepsilon},$$
namely \eqref{Eq8-1} holds.

Now we use Theorem 7. Take $\varrho=0.$ From the first estimate of \eqref{meanC} we take $\omega_0=\theta+\varepsilon$.
From \eqref{Eq8-1} we  take $\omega_1=(d^2-1)/2d^2<1/2$.  Suppose $k>1$ is a large but fixed integer which satisfies \eqref{k-condition}.
Take $x=T, N=[T^{2d-1-\varepsilon}]$.  $M$ is a real number such that $T^\varepsilon\ll M\ll \sqrt N$
and $\hat{E}(M)=0$ if $b(n)\ge 0 \ \  (n\ge 1)$ and $\hat{E}(M)\ll M^{\theta+\varepsilon}$ otherwise.

The truncated Tong-type formula of $E(y)$ was given in Theorem 7. Since we are not interested in the exact value of constants
like $\kappa_0$, we may assume $Q=1$ for simplicity of notation.  Hence we may take $\mu_n=n$ in the rest of this section.
Then from \eqref{E-I-expression}--\eqref{R6-expression} of Theorem 7 we have
\begin{equation}
E(y)=\sum_{j=1}^7R_j(y), \label{Eq8-4}
\end{equation}
where
\begin{align*}
R_1(y)&=\kappa_0y^{\frac{d-1}{2d}}\sum_{n\le M}\frac{a(n)}{n^{\frac{d+1}{2d}}}\cos\left(h(yn)^{1/d}+c_0\pi\right),\\
R_2(y)&=y^{\frac{d+1}{2d}}\Re \left(c_{00} I\left(-\frac{3}{2}+\frac{1}{2d},M,N,y\right)\right),\nonumber\\
R_3(y)&=\sum_{\substack{0\le l,m\le J \\l+m>0}} \Re \left(c_{lm} I\left(-\frac{3}{2}+\frac{1}{2d}+\frac{l-m}{d},M,N,y\right)\right)
        T^{-l}y^{-l+\frac{1}{2}+\frac{1}{2d}+\frac{l-m}{d}},\nonumber\\
R_4(y)&=\sum_{0\le l,m\le k}\Re \left(c_{lm}^{\prime}I\left(-\frac{3}{2}+\frac{1}{2d}-\frac{k+m}{d},N,\infty,y+\frac{l}{T}\right)\right)
        T^{k}\left(y+\frac lT\right)^{k+\frac{1}{2}+\frac{1}{2d}-\frac{k+m}{d}},\\[1ex]
R_5(y)& \ll T^{\frac 12-\frac{3}{2d}}M^{\max(\frac 12-\frac{3}{2d}, 0)}\log^A T+T^{-\frac 32+\frac{1}{2d}}M^{ \frac 12+\frac{1}{2d} }\\[1ex]
      & \quad + T^{- \frac{1}{2d}}M^{\omega_1 -1-\frac{1}{2d} }+T^{(2d-1)(1+\omega_1)-2k+1+4k/d-3d}\\[1ex]
      & \ll T^{\frac 12-\frac{3}{2d}}M^{\max(\frac 12-\frac{3}{2d}, 0)}\log^A T,\\[1ex]
R_6(y)& \begin{cases} = 0 & \mbox{if $b(n) \geq 0$} \\
                       \ll T^{ \frac 12-\frac{1}{2d}}M^{\theta- \frac 12-\frac{1}{2d}+\varepsilon} & \mbox{otherwise},
        \end{cases} \\
R_7(y)&=E(y)-\int_{\mathbf{E}_k}E(\tilde{y})d\mathrm{Y}_k.
\end{align*}

For simplicity, we write
\begin{equation*}
E(y)=K_1+K_2,
\end{equation*}
where
\begin{eqnarray*}
&&K_1=R_1(y)+R_2(y),\\
&&K_2=R_3(y)+R_4(y)+R_5(y)+R_6(y) +R_7(y).
\end{eqnarray*}


\subsection{Evaluation of $\int_T^{(1+\delta)T}K_1^2dy$}


In this subsection we shall evaluate  $\int_T^{(1+\delta)T}K_1^2dy$.

We write
\begin{align}\label{decomp}
R_1^2(y)&=\kappa_0^2y^{\frac{d-1}{d}}\sum_{n,m\le M}\frac{a(n)a(m)}{(nm)^{\frac{d+1}{2d}}}\cos\left(h(yn)^{1/d}+c_0\pi\right)
          \cos\left(h(ym)^{1/d}+c_0\pi\right)   \\
        &=\frac{\kappa_0^2}{2}y^{\frac{d-1}{d}}\sum_{n,m\le M}\frac{a(n)a(m)}{(nm)^{\frac{d+1}{2d}}}
           \left(\cos\left(h y^{1/d}(n^{1/d}-m^{1/d})\right)+\cos\left(h y^{1/d}(n^{1/d}+m^{1/d})+2c_0\pi\right)\right)\nonumber\\
        &=W_1(y)+W_2(y)+W_3(y),  \nonumber
\end{align}
where
\begin{align*}
W_1(y)&=\frac{\kappa_0^2}{2}y^{\frac{d-1}{d}}\sum_{n\le M}\frac{a^2(n)}{n^{\frac{d+1}{d}}},\\
W_2(y)&=\frac{\kappa_0^2}{2}y^{\frac{d-1}{d}}\sum_{\substack{n,m\le M \\ n \not= m}}
          \frac{a(n)a(m)}{(nm)^{\frac{d+1}{2d}}}\cos\left(h y^{1/d}(n^{1/d}-m^{1/d})\right),\\
W_3(y)&=\frac{\kappa_0^2}{2}y^{\frac{d-1}{d}}\sum_{n,m\le M}\frac{a(n)a(m)}{(nm)^{\frac{d+1}{2d}}}\cos\left(h y^{1/d}(n^{1/d}+m^{1/d})+2c_0\pi\right).
\end{align*}

For $W_1(y),$ we have
\begin{align}\label{W1}
\int_T^{(1+\delta)T}W_1(y)dy
&=\frac{\kappa_0^2}{2}\sum_{n\le M}\frac{a^2(n)}{n^{\frac{d+1}{d}}}\int_T^{(1+\delta)T}y^{\frac{d-1}{d}}dy\\
&=\frac{\kappa_0^2}{2}\sum_{n=1}^\infty\frac{a^2(n)}{n^{\frac{d+1}{d}}}\int_T^{(1+\delta)T}y^{\frac{d-1}{d}}dy+O(T^{2-1/d+\varepsilon}M^{-1/d}), \notag
\end{align}
where in the last step we used the second estimate of \eqref{meanC}.

By the first derivative test we have
\begin{align*}
\int_T^{(1+\delta)T}W_2(y)dy
&\ll T^{2-2/d}\sum_{\substack{n,m\le M \\ n \not= m}}\frac{|a(n)a(m)|}{(nm)^{\frac{d+1}{2d}}|n^{1/d}-m^{1/d}|}\\
&= T^{2-2/d}(\Sigma_1+\Sigma_2),
\end{align*}
say, where
\begin{align*}
\Sigma_1&=\sum_{\substack{n,m\le M \\ |n^{1/d}-m^{1/d}|\ge (mn)^{1/2d}/10d}} \frac{|a(n)a(m)|}{(nm)^{\frac{d+1}{2d}}|n^{1/d}-m^{1/d}|},\\
\Sigma_2&=\sum_{\substack{n,m\le M \\ 0<|n^{1/d}-m^{1/d}|< (mn)^{1/2d}/10d}} \frac{|a(n)a(m)|}{(nm)^{\frac{d+1}{2d}}|n^{1/d}-m^{1/d}|}.
\end{align*}

For $\Sigma_1$ we have
\begin{align*}
\Sigma_1& \ll \sum_{\substack{n,m\le M \\ |n^{1/d}-m^{1/d}|\ge (mn)^{1/2d}/10d}}\frac{|a(n)a(m)|}{(nm)^{\frac{d+2}{2d}}}\\
        & \ll \left(\sum_{n\leq M}\frac{|a(n)|}{n^{\frac{d+2}{2d}}}\right)^2\ll M^{1-2/d+\varepsilon} ,
\end{align*}
where we used the estimate $\sum_{n\le y}|a(n)|\ll y^{1+\varepsilon/2},$
which follows from the second estimate of \eqref{meanC} and Cauchy's inequality.

Next we consider the sum $\Sigma_2$. By Lagrange's mean value theorem we have $|n^{1/d}-m^{1/d}|\gg n_0^{1/d-1}|n-m| $ for some $n_0$ between $n$ and $m$.
If $m, n$ are contained in the sum $\Sigma_2,$ then $m\asymp n$ and hence
$$
|n^{1/d}-m^{1/d}| \gg  (mn)^{1/2d-1/2}|m-n|.
$$
Thus 
\begin{align*}
\Sigma_2 &\ll \sum_{n \not= m}\frac{|a(n)a(m)|}{(nm)^{1/d}|n-m|} \\
& \ll \sum_{n \neq m} \frac{1}{|n-m|} \left\{\left(\frac{a(n)}{n^{1/d}}\right)^2+ \left(\frac{a(m)}{m^{1/d}}\right)^2\right\}\\
& \ll \sum_{n \neq m} \frac{1}{|n-m|} \frac{a^2(n)}{n^{2/d}}\\
& \ll \sum_{n \leq M} \frac{a^2(n)}{n^{2/d}} \sum_{m \neq n} \frac{1}{|n-m|} \\
& \ll M^{1-2/d+\varepsilon}.
\end{align*}

From the above two estimates we get
\begin{eqnarray}
\int_T^{(1+\delta)T}W_2(y)dy\ll T^{2-2/d+\varepsilon}M^{1-2/d}. \label{W2}
\end{eqnarray}

We have similarly that
\begin{eqnarray}
\int_T^{(1+\delta)T}W_3(y)dy\ll T^{2-2/d+\varepsilon}M^{1-2/d}. \label{W3}
\end{eqnarray}

From \eqref{decomp}, \eqref{W1}, \eqref{W2} and \eqref{W3} we have
\begin{align}\label{R1}
\int_T^{(1+\delta)T}R_1^2(y)dy
=& \frac{\kappa_0^2}{2}\sum_{n=1}^\infty\frac{a^2(n)}{n^{\frac{d+1}{d}}}\int_T^{(1+\delta)T}y^{\frac{d-1}{d}}dy  \\
 &  +O(T^{2-1/d+\varepsilon}M^{-1/d} +T^{2-2/d+\varepsilon}M^{1-2/d} ). \nonumber
\end{align}

From \eqref{sigmastarcondition} we have $2\sigma^{*}-1-1/d<0,$ so by Lemma 9 we have
\begin{align}\label{R2}
\int_T^{(1+\delta)T}R_2^2(y)dy
&\ll T^{\frac{d+1}{d}}\int_T^{(1+\delta)T}\left|I\left(-\frac 32+\frac{1}{2d}, M, N, y\right)\right|^2dy \\
&\ll T^{2-1/d+\varepsilon}\max_{M\leq P\le N}P^{2\sigma^{*}-1-1/d}\nonumber\\
&\ll T^{2-1/d+\varepsilon} M^{2\sigma^{*}-1-1/d}. \nonumber
\end{align}

By the definitions of $R_1(y)$ and $R_2(y)$ we have
\begin{align*}
& \int_T^{(1+\delta)T}R_1(y)R_2(y)dy\\
&= \Re \kappa_0 c_{00} \int_T^{(1+\delta)T}yI\left(-\frac 32+\frac{1}{2d}, M, N, y\right)\sum_{n\leq M}
   \frac{a(n)}{n^{\frac{d+1}{2d}}}\cos(h(yn)^{1/d}+c_0\pi)dy\\
&= \Re \kappa_0 c_{00} (I_1+I_2),
\end{align*}
where
\begin{align*}
I_1 &=\int_T^{(1+\delta)T}yI\left(-\frac 32+\frac{1}{2d}, M, N, y\right)\sum_{n\le M/2}
\frac{a(n)}{n^{\frac{d+1}{2d}}}\cos(h(yn)^{1/d}+c_0\pi)dy,\\
I_2 &=\int_T^{(1+\delta)T}yI\left(-\frac 32+\frac{1}{2d}, M, N, y\right)\sum_{M/2<n\le M}
\frac{a(n)}{n^{\frac{d+1}{2d}}}\cos(h(yn)^{1/d}+c_0\pi)dy.
\end{align*}

By Lemma 8 we get
\begin{align*}
I_1&=\sum_{n\le M/2}
\frac{a(n)}{n^{\frac{d+1}{2d}}}\int_T^{(1+\delta)T}yI\left(-\frac 32+\frac{1}{2d}, M, N, y\right)\cos(h(yn)^{1/d}+c_0\pi)dy\\
&\ll \sum_{n\le M/2} \frac{|a(n)|}{n^{\frac{d+1}{2d}}}\left|\int_T^{(1+\delta)T}yI\left(-\frac 32+\frac{1}{2d}, M, N, y\right)
      \cos(h(yn)^{1/d}+c_0\pi)dy\right|\\
&\ll \sum_{n\le M/2}  \frac{|a(n)|}{n^{\frac{d+1}{2d}}}T^{2-3/2d+\varepsilon}\max_{M\leq P\leq N}P^{\sigma^{*}-1/2-1/d}\\
&\ll \sum_{n\le M/2} \frac{|a(n)|}{n^{\frac{d+1}{2d}}}T^{2-3/2d+\varepsilon}M^{\sigma^{*}-1/2-1/d}\\
&\ll T^{2-3/2d+\varepsilon}M^{\sigma^{*}-3/2d}.
\end{align*}

By Cauchy's inequality
$$I_2\ll T(V_1V_2)^{1/2},$$
where
\begin{eqnarray*}
&&V_1=\int_T^{(1+\delta)T}\left|I\left(-\frac 32+\frac{1}{2d}, M, N, y\right)\right|^2dy,\\
&&V_2=\int_T^{(1+\delta)T} \left(\sum_{M/2<n\le M}
\frac{a(n)}{n^{\frac{d+1}{2d}}}\cos(h(yn)^{1/d}+c_0\pi)\right)^2dy.
\end{eqnarray*}

By Lemma 9 we get
$$V_1\ll T^{1-2/d+\varepsilon}M^{2\sigma^{*}-1-1/d}.$$

By using the same approach of $\int_T^{(1+\delta)T}R_1^2(y)dy,$ we get
$$V_2\ll T^{1+\varepsilon}M^{-1/d} +T^{1-1/d+\varepsilon}M^{1-2/d}.$$

Thus
$$I_2\ll T^{2-1/d+\varepsilon}M^{\sigma^{*}-1/2-1/d}+T^{2-3/2d+\varepsilon}M^{\sigma^{*}-3/2d}.$$

From the estimates of $I_1$ and $I_2$
we get
\begin{equation} \label{R1R2}
\int_T^{(1+\delta)T}R_1(y)R_2(y)dy\ll T^{2-1/d+\varepsilon}M^{\sigma^{*}-1/2-1/d}+T^{2-3/2d+\varepsilon}M^{\sigma^{*}-3/2d}.
\end{equation}

Combining \eqref{R1}, \eqref{R2} and \eqref{R1R2} we get
\begin{align}\label{K1}
\int_T^{(1+\delta)T}K_1^2dy& =\frac{\kappa_0^2}{2}\sum_{n=1}^\infty\frac{a^2(n)}{n^{\frac{d+1}{d}}}\int_T^{(1+\delta)T}y^{\frac{d-1}{d}}dy\\
& \quad   +O(T^{2-1/d+\varepsilon}M^{2\sigma^{*}-1-1/d}+T^{2-2/d+\varepsilon}M^{1-2/d} ), \nonumber
\end{align}
where we used the following estimates (since $\sigma^{*}\geq 1/2$)
\begin{eqnarray*}
&&T^{2-1/d+\varepsilon}M^{-1/d} \ll T^{2-1/d+\varepsilon}M^{2\sigma^{*}-1-1/d},\\[1ex]
 &&T^{2-1/d+\varepsilon}M^{\sigma^{*}-1/2-1/d}\ll   T^{2-1/d+\varepsilon}M^{2\sigma^{*}-1-1/d},\\
 &&T^{2-3/2d+\varepsilon}M^{\sigma^{*}-3/2d}\ll
 \left(T^{2-1/d+\varepsilon}M^{2\sigma^{*}-1-1/d}\right)^{1/2}\left(T^{2-2/d+\varepsilon}M^{1-2/d} \right)^{1/2}.
 \end{eqnarray*}


\subsection{Upper bound of $\int_T^{(1+\delta)T}K_2^2dy$ }


By Cauchy's inequality and Lemma 9, we have
\begin{align*}
& \int_T^{(1+\delta)T}R_3^2(y)dy\\
&\ll \sum_{\substack{0\le l,m\le J \\[2pt] l+m>0}} T^{-2l}\int_T^{(1+\delta)T}
     y^{-2l+1+\frac 1d+\frac{2l-2m}{d}}\left|I\left(-\frac 32+\frac{1}{2d}+\frac{l-m}{d}, M, N, y\right)\right|^2dy\\
&\ll \sum_{\substack{0\le l,m\le J \\[2pt] l+m>0}} T^{-4l+1+\frac 1d+\frac{2l-2m}{d}}\int_T^{(1+\delta)T}
     \left|I\left(-\frac 32+\frac{1}{2d}+\frac{l-m}{d}, M, N, y\right)\right|^2dy\\
&\ll \sum_{\substack{0\le l,m\le J \\[2pt] l+m>0}} T^{-4l+1+\frac 1d+\frac{2l-2m}{d}+1-\frac 2d+\varepsilon}
     \max_{M\leq P\leq N}P^{2\sigma^{*}-1-\frac 1d+\frac{2l-2m}{d}}\\
&=\Sigma_3+\Sigma_4,
\end{align*}
say, where
\begin{align*}
\Sigma_3 &=\sum_{\substack{0\le l\le m\le J \\[2pt] l+m>0}} T^{-4l+1+\frac 1d+\frac{2l-2m}{d}+1-\frac 2d+\varepsilon}
\max_{M\leq P\leq N}P^{2\sigma^{*}-1-\frac 1d+\frac{2l-2m}{d}},\\
\Sigma_4 &=\sum_{0\le m<l\le J}T^{-4l+1+\frac 1d+\frac{2l-2m}{d}+1-\frac 2d+\varepsilon}\max_{M\leq P\leq N}P^{2\sigma^{*}-1-\frac 1d+\frac{2l-2m}{d}}.
\end{align*}
For $\Sigma_3$ we have by \eqref{sigmastarcondition} that
\begin{align*}
\Sigma_3 &\ll \sum_{\substack{0\le l\le m\le J \\[2pt] l+m>0}} T^{-4l+1+\frac 1d+\frac{2l-2m}{d}+1-\frac 2d+\varepsilon} M^{2\sigma^{*}-1-\frac1d+\frac{2l-2m}{d}}\\
&\ll T^{2-1/d+\varepsilon}M^{2\sigma^{*}-1-1/d}\sum_{1\le m\le J}T^{-2m/d}M^{-2m/d}\\
&\quad +T^{2-1/d+\varepsilon}M^{2\sigma^{*}-1-1/d}\sum_{1\le l\le J}T^{-4l+2l/d}M^{2l/d}\sum_{l\le m\le J}T^{-2m/d}M^{-2m/d}\\
&\ll T^{2-3/d+\varepsilon}M^{2\sigma^{*}-1-3/d}.
\end{align*}
For $\Sigma_4$ we have (since $\sigma^{*}\ge 1/2$ and recall $N=[T^{2d-1-\varepsilon}]$)
\begin{align*}
\Sigma_4 &\ll \sum_{0\le m<l\le J}T^{-4l+1+\frac 1d+\frac{2l-2m}{d}+1-\frac 2d+\varepsilon} N^{2\sigma^{*}-1-\frac 1d+\frac{2l-2m}{d}}\\
&\ll T^{2-1/d+\varepsilon}N^{2\sigma^{*}-1-1/d}\sum_{0\le m\le J}T^{-2m/d}N^{-2m/d}\sum_{m+1\le l\le J}
     T^{-4l+\frac{2l}{d}} N^{\frac{2l}{d}}\\
&\ll T^{2-1/d+\varepsilon}N^{2\sigma^{*}-1-1/d}.
\end{align*}
From the above estimates we get
\begin{eqnarray} \label{R3}
\int_T^{(1+\delta)T}R_3^2(y)dy\ll T^{2-3/d+\varepsilon}M^{2\sigma^{*}-1-3/d}+T^{2-1/d+\varepsilon}N^{2\sigma^{*}-1-1/d}.
\end{eqnarray}

By Cauchy's inequality and Lemma 10, we get (recall that $N=[T^{2d-1-\varepsilon}]$)
\begin{align}\label{R4}
&  \int_T^{(1+\delta)T}R_4^2(y)dy \\
& \ll \sum_{0\le l,m\le k}T^{2k}\int_T^{(1+\delta)T}(y+\frac lT)^{2k+1+\frac 1d-\frac{2k+2m}{d}}
      \left|I\left(-\frac{3}{2}+\frac{1}{2d}-\frac{k+m}{d},N,\infty,y+\frac{l}{T}\right)\right|^2dy\nonumber\\
&\ll \sum_{0\le l,m\le k}T^{4k+1+\frac 1d-\frac{2k+2m}{d}+1-\frac 2d+\varepsilon}N^{2\sigma^*-1-\frac 1d-\frac{2k+2m}{d}}\nonumber\\
&\ll  T^{4k+2-\frac 1d-\frac{2k}{d}+\varepsilon}N^{2\sigma^*-1-\frac 1d-\frac{2k}{d}}\ll
T^{2-\frac 1d+\varepsilon}N^{2\sigma^*-1-\frac 1d}. \nonumber
\end{align}

For $R_5(y)$ we have trivially that
\begin{align}
\int_T^{(1+\delta)T}R_5^2(y)dy
\ll \left\{\begin{array}{ll} T^{1/2},&\mbox{if $d=2,$}\\[1ex]
                             T^{2-3/d+\varepsilon}M^{1-3/d},& \mbox{if $d\ge 3.$}
           \end{array}\right. \label{R5}
\end{align}

For $R_6(y)$ we have trivially that
\begin{align}
\int_T^{(1+\delta)T}R_6^2(y)dy \ll
T^{2-1/d+\varepsilon}M^{2\theta-1-1/d+\varepsilon}. \label{R6}
\end{align}

Next we estimate $\int_T^{(1+\delta)T}R_7^2(y)dy.$ Trivially we have
\begin{align*}
\int_T^{(1+\delta)T}R_7^2(y)dy
&\ll \sum_{T\leq n\leq (1+\delta)T+1}\int_{n-1}^{n}R_7^2(y)dy\\
&\ll \sum_{T\leq n\leq (1+\delta)T+1}\int_{n-1}^{n-k/T}R_7^2(y)dy+\sum_{T\leq n\leq (1+\delta)T+1}\int_{n-k/T}^{n}R_7^2(y)dy .
\end{align*}
It is easy to see that
\begin{eqnarray*}
R_7(y)=E(y)-\int_{\mathbf{E}_k}E(\tilde{y})d\mathrm{Y}_k=\int_{\mathbf{E}_k}(E(y)-E(\tilde{y}))d\mathrm{Y}_k.
\end{eqnarray*}
Suppose $n-1<y\leq n.$ If $n-1<y\le \tilde{y}<n,$ then by the definition of $E(u)$ we get
$$E(\tilde{y})-E(y)=Q(y)-Q(\tilde{y})\ll T^{-1}\max_{y\ll u\ll \tilde{y}}Q^{\prime}(u)\ll T^{-1}\log^{m_{\cal L}} T;$$
otherwise if $\tilde{y}\geq n,$ then
$$E(\tilde{y})-E(y)=Q(y)-Q(\tilde{y})+a(n)/2\ll T^{-1}\log^{m_{\cal L} } T+|a(n)|.$$

Thus we have
\begin{eqnarray*}
&&\sum_{T\leq n\leq (1+\delta)T+1}\int_{n-1}^{n-k/T}R_7^2(y)dy\ll  T^{-1}\log^{2m_{\cal L} } T,\\
&&\sum_{T\leq n\leq (1+\delta)T+1}\int_{n-k/T}^{n}R_7^2(y)dy\ll T^{-1}\log^{2 m_{\cal L}} T+T^{-1}\sum_{n\ll T}a^2(n)\ll   T^\varepsilon.
\end{eqnarray*}

Hence we get
\begin{equation} \label{R7}
\int_T^{(1+\delta)T}R_7^2(y)dy\ll   T^\varepsilon.
\end{equation}

The first term in the right hand side of \eqref{R3} is absorbed in the right hand side of \eqref{R5}. Furthermore
from the assumption $M \ll N^{1/2}$ and $1/2 \leq \sigma^{*}<1/2+1/2d$, we have
$$
N^{2\sigma^{*}-1-1/d} \ll M^{2(2\sigma^*-1- 1/d)}
$$
and also
$$
M^{2\theta-1-1/d} \ll M^{2(2\sigma^*-1- 1/d)}
$$
by $\theta < 1/2-1/2d$. Hence from \eqref{R3}--\eqref{R7} we have
\begin{align}
\int_T^{(1+\delta)T}K_2^2dy
&\ll T^{2-1/d+\varepsilon}M^{2(2\sigma^{\ast}-1-1/d)}
         +  \left\{\begin{array}{ll} T^{1/2},&\mbox{if $d=2,$}\\[1ex]
                                   T^{2-3/d+\varepsilon}M^{1-3/d},& \mbox{if $d\ge 3.$}
                 \end{array}\right.   \label{K2}
\end{align}


\subsection{Proof of Theorem 1}


Now we take $M$   such that
\begin{align}
& M \asymp T^{\frac{1}{2d(1-\sigma^{*})-1}}, \label{eq8-21}\\
\intertext{and}
&\hat{E}(M)\left\{\begin{array}{ll} =0,&\mbox{if $b(n)\ge 0\ (n\ge 1)$,}\\[1ex]
                                    \ll T^{\theta+\varepsilon},& \mbox{otherwise.}
\end{array}\right.\nonumber
\end{align}

By noting that $\sigma^{*}=1/2$ for $d=2$, we can see easily that $M \ll N^{1/2}$ for any $d \geq 2$.
Then the formula \eqref{K1} becomes
\begin{equation} \label{last-K1}
\int_T^{(1+\delta)T}K_1^2dy =\frac{\kappa_0^2}{2}\sum_{n=1}^\infty\frac{a^2(n)}{n^{\frac{d+1}{d}}}\int_T^{(1+\delta)T}y^{\frac{d-1}{d}}dy
+O(T^{2-\frac{3-4\sigma^{*}}{2d(1-\sigma^{*})-1}+\varepsilon}).
\end{equation}

From \eqref{K2}, \eqref{last-K1} and Cauchy's inequality we get
\begin{align}\label{last-K1K2}
\int_T^{(1+\delta)T}K_1K_2dy
&\ll T^{2-1/d+\varepsilon}M^{(2\sigma^*-1- 1/d)}  
     +  \left\{\begin{array}{ll} T^{\varepsilon},  &\mbox{if $d=2$}\\[1ex]
                                   T^{2-2/d+\varepsilon} M^{(1-3/d)/2},& \mbox{if $d\ge 3$}
                 \end{array}\right.  \\[1ex]
& \ll T^{2-\frac{3-4\sigma^{*}}{2d(1-\sigma^{*})-1}+\varepsilon}, \nonumber
\end{align}
where the last inequality follows from the choice of $M$ \eqref{eq8-21}.
It is also seen easily that
\begin{align}
\int_T^{(1+\delta)T}K_2^2dy \ll T^{2-\frac{3-4\sigma^{*}}{2d(1-\sigma^{*})-1}+\varepsilon} \label{last-K2}
\end{align}
by (\ref{K2}) and the choice of $M$.
Hence from \eqref{last-K1}, (\ref{last-K1K2}) and (\ref{last-K2}) we get
\begin{equation}
  \int_T^{(1+\delta)T}E^2(y)dy  =\frac{\kappa_0^2}{2}\sum_{n=1}^\infty\frac{a^2(n)}{n^{\frac{d+1}{d}}}\int_T^{(1+\delta)T}y^{\frac{d-1}{d}}dy
 +O(T^{2-\frac{3-4\sigma^{*}}{2d(1-\sigma^{*})-1}+\varepsilon} ). \label{Esquare}
\end{equation}
Finally from \eqref{Esquare} we get
\begin{align*}
\int_0^TE^2(y)dy
&=\sum_{0\le j\le \frac{\log T}{\delta}}\int_{\frac{T}{(1+\delta)^{j+1}}}^{\frac{T}{(1+\delta)^j}}E^2(y)dy+O(T^{2\delta})\\
&=\sum_{0\le j\le \frac{\log T}{\delta}}  \frac{\kappa_0^2}{2}\sum_{n=1}^\infty\frac{a^2(n)}{n^{\frac{d+1}{d}}}
    \int_{\frac{T}{(1+\delta)^{j+1}}}^{\frac{T}{(1+\delta)^j}}y^{\frac{d-1}{d}}dy \\
&\quad +\sum_{0\le j\le \frac{\log T}{\delta}}O\left(\left(\frac{T}{(1+\delta)^j}\right)^{2-\frac{3-4\sigma^{*}}{2d(1-\sigma^{*})-1}
       +\varepsilon}\right)\nonumber +O(T^{2\delta})  \\
&=\frac{\kappa_0^2}{2}\sum_{n=1}^\infty\frac{a^2(n)}{n^{\frac{d+1}{d}}}  \int_{0}^{T}y^{\frac{d-1}{d}}dy
     +O\left(T^{2-\frac{3-4\sigma^{*}}{2d(1-\sigma^{*})-1}+\varepsilon}\right).
\end{align*}


\subsection{Proofs of Corollaries \ref{cor-1}     and \ref{cor-2}}


We first prove Corollary \ref{cor-1}. Suppose that $0\le \theta\le 1/4$ is a
real number and ${\cal L}(s)\in {\cal S}_{real}^{\theta}$ is a function of degree $2$.
In this case it is well-known that
\begin{equation} \label{meansquare_degree2}
\int_0^T|{\cal L}(1/2+it)|^2dt\ll T^{1+\varepsilon},
\end{equation}
namely, we have $\sigma^{*}=1/2<3/4$. So Corollary \ref{cor-1} follows from Theorem 1.

\medskip

Next we prove Corollary \ref{cor-2}. Suppose that $0\le \theta\le 1/3$ is a real number and
${\cal L}(s)\in {\cal S}_{real}^{\theta}$ is a function of degree $3$ which can be written as
${\cal L}(s)={\cal L}_1(s){\cal L}_2(s)$,  where $   {\cal L}_1(s)\in  {\cal S}_{real}^{\theta}$
is a function of degree $1$ and  ${\cal L}_2(s)\in {\cal S}_{real}^{\theta}$ is a function of degree $2.$
We can show that
\begin{equation} \label{MSQ-3-version1}
\int_0^T |{\cal L}(5/8+it)|^{2}dt\ll T^{1+\varepsilon}.
\end{equation}
Hence we can take $\sigma^{*}=5/8$ in \eqref{meanL} and the assertion \eqref{unconditional} follows from Theorem 1.
Furthermore if we assume \eqref{assump-d2},  we can show that
\begin{equation} \label{MSQ-3-version2}
\int_0^T |{\cal L}(7/12+it)|^{2}dt\ll T^{1+\varepsilon}.
\end{equation}
Hence we can take $\sigma^{*}=7/12$ in this case and we get a better assertion \eqref{conditional}.


\medskip

In order to prove \eqref{MSQ-3-version1} or \eqref{MSQ-3-version2}, we follow a method of Ivi\'{c} \cite{I8}.
In \cite{I8}, Ivi\'c treated the automorphic $L$-function attached to a cusp form over $SL_2({\Bbb Z})$.
Our proof for the general case is due to Ivi\'c with some modifications. So here we give a detailed  proof.

\begin{lem}
Let $t_1<\cdots <t_R$ be real numbers such that $T\le t_r\le 2T$ for $r=1,2,\ldots,R$ and $|t_r-t_s|\ge \log^4 T$
for $1\le s\not= r\le R.$ Suppose $V$ and $1\ll M\ll T^C \ (C>0)$ are large parameters such that
$$
T^\varepsilon <V\le \left|\sum_{M<n\le 2M} a(n)n^{-\sigma-it_r}\right|
$$
with
\begin{equation}
\sum_{M\le n\le 2M}|a(n)|^2\ll M^{1+\varepsilon},  \label{keisu-2jyouwa} 
\end{equation}
then we have
$$
R\ll(M^{2-2\sigma}V^{-2}+TV^{-\frac{2}{3-4\sigma}})T^\varepsilon
$$
for $1/2< \sigma\le 2/3$.
\end{lem}

\medskip

\noindent {\bf Remark to Lemma 11.}
This is the case $1/2<\sigma\le 2/3$ of Lemma 8.2 of Ivi\'c \cite{I1}, where Ivi\'c supposed the condition
$a(n)\ll M^\varepsilon$ for $M<n\le 2M.$ However we see easily from Ivi\'c's argument that this condition can be relaxed to
the condition \eqref{keisu-2jyouwa}. 

\medskip


\begin{lem}  \label{lemma_degree2}
Let $0\le \theta\le 1/3$ be a fixed real number and let $\mathfrak{L}(s)$ be an element of ${\cal S}_{real}^{\theta}$ with degree 2. Then we have
\begin{equation}
\int_0^T|\mathfrak{L}(\sigma+it)|^{1/(1-\sigma)}dt \ll T^{1+\varepsilon} \label{MSQ2-version1}
\end{equation}
for $1/2<\sigma<2/3$. Especially we have
\begin{equation}
\int_0^T|\mathfrak{L}(5/8+it)|^{8/3}dt \ll T^{1+\varepsilon}. \label{MSQ2-version1-rei}
\end{equation}
Furthermore if we assume
\begin{equation}
\int_0^T|\mathfrak{L}(1/2+it)|^{6}dt \ll T^{2+\varepsilon}, \label{six-power}
\end{equation}
we have
\begin{equation}
\int_0^T|\mathfrak{L}(\sigma+it)|^{2/(3-4\sigma)}dt \ll T^{1+\varepsilon} \label{MSQ2-version2}
\end{equation}
for $1/2<\sigma<5/8$. Especially we have
\begin{equation}
\int_0^T|\mathfrak{L}(7/12+it)|^{3}dt \ll T^{1+\varepsilon}. \label{MSQ2-version2-rei}
\end{equation}
\end{lem}

\proof
We follow the approach of Ivi\'c  \cite{I8}.  Suppose $T<t_1<t_2<\cdots<t_R<2T$ such that
\begin{equation} \label{LargeValue}
|\mathfrak{L}( \sigma+it_r)|\geq V\gg  T^\varepsilon
\end{equation}
and
\begin{equation*}
|t_r-t_s|\gg \log^4 T \ (1\le r\not= s\le R).
\end{equation*}

Suppose that
$$
\mathfrak{L}(s)=\sum_{n\ge 1}a(n)n^{-s}, \ \ \Re s>1.
$$
From the formula (see (A7) of Ivi\'c \cite{I1})
$$
e^{-v}=(2\pi i)^{-1}\int_{b-i\infty}^{b+i\infty}\Gamma(w)v^{-w}dw\  (b,v>0),
$$
we have
\begin{equation} \label{L-rep}
\sum_{n\ge 1}a(n)e^{-n/Y}n^{-s}=\frac{1}{2\pi i}\int_{2-i\infty}^{2+i\infty}\Gamma(w)Y^w\mathfrak{L}(s+w)dw,
\end{equation}
where $1\ll Y\ll T^A$ is a parameter to be chosen later and $A>0$ is an arbitrary constant.
Let $s$ be any one of $\sigma+it_r$ with  $\sigma>1/2$ in \eqref{L-rep}.
Moving the line of integration  to $1/2-\sigma$ in \eqref{L-rep} we get by the residue theorem that
\begin{eqnarray*}
\mathfrak{L}(s)=\sum_{n\ge 1}a(n)e^{-n/Y}n^{-s}
-\frac{1}{2\pi i}\int_{1/2-\sigma-i\infty}^{1/2-\sigma+i\infty}\Gamma(w)Y^w\mathfrak{L}(s+w)dw+O(T^{-B}),
\end{eqnarray*}
where $B>1$ is a large fixed constant. So we get
\begin{eqnarray*}
\mathfrak{L}(\sigma+it_r) \ll \left|\sum_{n\le Y\log^2 T}a(n)e^{-n/Y}n^{-\sigma-it_r}\right|
+ Y^{1/2-\sigma}\int_{-\log^2 T}^{\log^2 T} | \mathfrak{L}(1/2+it_r+iv)|dv+1,
 \end{eqnarray*}
which implies that if \eqref{LargeValue} holds, we must have
\begin{align}\label{Case1}
V&\ll \left|\sum_{n\le Y\log^2 T}a(n)e^{-n/Y}n^{- \sigma-it_r}\right|  \\
 &\ll \log T \max_{N_1\le \frac{Y\log^2 T}{2}}\left|\sum_{N_1<n\le 2N_1}a(n)e^{-n/Y}n^{- \sigma-it_r}\right| \nonumber
\intertext{or}
V & \ll  Y^{1/2-\sigma}\int_{-\log^2 T}^{\log^2 T} | \mathfrak{L}(1/2+it_r+iv)|dv. \label{Case2}
\end{align}
Let $\mathcal{A}_1$ be the set of points $t_r$ which satisfy \eqref{Case1} and $\mathcal{A}_2$ the set of points $t_r$
which satisfy \eqref{Case2}. We also let $R_1=|\mathcal{A}_1|$ and $R_2=|\mathcal{A}_2|$, whence $R=R_1+R_2$.

For $R_1$ we have
\begin{equation}  \label{Ivic-8-2}
R_1  \ll T^{\varepsilon} (Y^{2-2\sigma}V^{-2}+TV^{-2/(3-4\sigma)})
\end{equation}
by Lemma 11.

Now we are going to estimate $R_2$. Since $t_r \in \mathcal{A}_2$ satisfies \eqref{Case2} we have by Cauchy's inequality
\begin{align*}
V &\ll Y^{1/2-\sigma} \log T \left(\int_{-\log^2 T}^{\log^2 T} |\mathfrak{L}(1/2+it_r+iv)|^2dv\right)^{1/2}.
\end{align*}
Squaring both sides and summing up over $t_r \in \mathcal{A}_2$, we get
\begin{align}\label{R2-firstestimate}
R_2 & \ll Y^{1-2\sigma} V^{-2}(\log T)^2 \sum_{t_r \in \mathcal{A}_2} \int_{-\log^2 T}^{\log^2 T} |\mathfrak{L}(1/2+it_r+iv)|^2dv \\
& \ll Y^{1-2\sigma} V^{-2}(\log T)^2 \int_{T-\log^2 T}^{2T+\log^2 T} |\mathfrak{L}(1/2+iv)|^2dv \nonumber \\
& \ll Y^{1-2\sigma} V^{-2} T^{1+\varepsilon}, \nonumber
\end{align}
where in the last inequality we used the mean square estimate \eqref{meansquare_degree2} for $\mathfrak{L}(1/2+it)$ of degree 2.
Hence as a first estimate we get, from \eqref{Ivic-8-2} and \eqref{R2-firstestimate},
\begin{equation} \label{from_meansquare}
R \ll T^{\varepsilon}\left(Y^{2-2\sigma}V^{-2}+TV^{-2/(3-4\sigma)}+Y^{1-2\sigma}V^{-2}T\right).
\end{equation}

We shall derive another estimate for $R_2$. It is convenient to consider the conditional case first, so we assume \eqref{six-power}.
This time we use H\"older's inequality in \eqref{Case2}. Thus we have
\begin{align*}
V \ll Y^{1/2-\sigma} (\log T)^{5/3} \left(\int_{-\log^2 T}^{\log^2 T} |\mathfrak{L}(1/2+it_r+iv)|^6dv\right)^{1/6},
\end{align*}
and
\begin{align}\label{R2-secondesitmate}
R_2 & \ll Y^{3-6\sigma} V^{-6}(\log T)^{10} \sum_{t_r \in \mathcal{A}_2} \int_{-\log^2 T}^{\log^2 T} |\mathfrak{L}(1/2+it_r+iv)|^6dv  \\
& \ll Y^{3-6\sigma} V^{-6}(\log T)^{10} \int_{T-\log^2 T}^{2T+\log^2 T} |\mathfrak{L}(1/2+iv)|^6dv \nonumber \\
& \ll Y^{3-6\sigma} V^{-6} T^{2+\varepsilon}. \nonumber
\end{align}
Hence \eqref{Ivic-8-2} and \eqref{R2-secondesitmate} give the second estimate
\begin{equation} \label{from_sixpower}
R \ll T^{\varepsilon}\left(Y^{2-2\sigma}V^{-2}+TV^{-2/(3-4\sigma)}+Y^{3-6\sigma}V^{-6}T^2\right).
\end{equation}
If $V \leq T^{(3-4\sigma)/4}$, we apply \eqref{from_meansquare} with $Y=T$. Then
\begin{align*}
R &\ll (Y^{2-2\sigma}V^{-2}+TV^{-\frac{2}{3-4\sigma}}+TY^{1 -2\sigma}V^{-2})T^\varepsilon  \\
&\ll (T^{2-2\sigma}V^{-2}+TV^{-\frac{2}{3-4\sigma}} )T^\varepsilon \\
&\ll T^{1+\varepsilon}V^{-\frac{2}{3-4\sigma}}  
\end{align*}
by recalling the condition $V\le T^{(3-4\sigma)/4}.$
Next if $V \geq T^{(3-4\sigma)/4}$, we apply \eqref{from_sixpower} with $Y=T^{2/(4\sigma-1)}V^{-4/(4\sigma-1)}+1$ and get
\begin{align*}
R &\ll (Y^{2-2\sigma}V^{-2}+T V^{-\frac{2}{3-4\sigma}}+ T^2 Y^{3-6\sigma}V^{-6} )T^\varepsilon  \\
&\ll ( TV^{-\frac{2}{3-4\sigma}}+T^{\frac{4-4\sigma}{4\sigma-1}}V^{-\frac{6}{4\sigma-1}} )T^\varepsilon \\
&\ll T^{1+\varepsilon}V^{-\frac{2}{3-4\sigma}} 
\end{align*}
if $T^{\frac{5-8\sigma}{4\sigma-1}}\le V^{ \frac{4(5-8\sigma)}{(4\sigma-1)(3-4\sigma)}},$
which is true when $V \geq  T^{(3-4\sigma)/4} $ and $1/2<\sigma\le 5/8.$
Summing up, if $1/2<\sigma \leq 5/8$,
$$
R \ll TV^{-2/(3-4\sigma)}
$$
holds true for all $V$, completing the proof of \eqref{MSQ2-version2}.

Next we shall prove \eqref{MSQ2-version1}. We use the fourth-power mean value theorem:
\begin{equation} \label{four-power}
\int_0^T|\mathfrak{L}(1/2+it)|^4dt \ll T^{2+\varepsilon},
\end{equation}
which is known to be true unconditionally (see e.g. Kanemitsu, Sankaranarayanan and Tanigawa \cite{KST}).
By the similar way we have
\begin{equation} \label{from_fourpower}
R \ll T^{\varepsilon}\left(Y^{2-2\sigma}V^{-2}+TV^{-2/(3-4\sigma)}+Y^{2-4\sigma}V^{-4}T^2\right)
\end{equation}
from \eqref{four-power}. If $V \leq T^{1-\sigma}$, we take $Y=T$ in \eqref{from_meansquare} and get
\begin{equation*}
R \ll T^{\varepsilon}\left(T^{2-2\sigma}V^{-2}+TV^{-2/(3-4\sigma)}\right) \ll T^{1+\varepsilon} V^{-1/(1-\sigma)}.
\end{equation*}
On the other hand if $V \geq T^{1-\sigma}$,  we take $V=(TV^{-1})^{1/\sigma}$ in \eqref{from_fourpower} and get
\begin{equation*}
R \ll T^{2-2\sigma)/\sigma}V^{-2/\sigma}+TV^{-2/(3-4\sigma)} \ll TV^{-1/(1-\sigma)}.
\end{equation*}
Therefore the estimate
$$
R \ll T^{1+\varepsilon} V^{-1/(1-\sigma)}
$$
holds true for all $V$ for $1/2<\sigma<2/3$. This completes the proof of \eqref{MSQ2-version1}.

\bigskip

Now we prove \eqref{MSQ-3-version1} and \eqref{MSQ-3-version2}.
Recall that   ${\cal L}(s)={\cal L}_1(s){\cal L}_2(s)$,  where ${\cal L}_1(s)\in  {\cal S}_{real}^{\theta}$
is a   function of degree $1$ and  $   {\cal L}_2(s)\in  {\cal S}_{real}^{\theta}$ is a function of degree $2.$
 By the work of Kaczorowski and Perelli \cite{KP1}, it is known that the functions
of degree one in the Selberg class are the Riemann zeta-function $\zeta(s)$ and the shifts
$L(s+i\vartheta, \chi)$ of Dirichlet $L$-functions attached to primitive characters $\chi$ with
$\vartheta\in {\Bbb R}.$
So we have
\begin{equation*}
\int_0^T |{\cal L}_1(5/8+it)|^8dt\ll T \log^{4} T
\end{equation*}
and 
\begin{equation*}
\int_0^T |{\cal L}_2(5/8+it)|^{8/3}dt\ll T^{1+\varepsilon}.
\end{equation*}
Hence by H\"older's inequality and \eqref{MSQ2-version1-rei}
\begin{align*}
\int_0^T |{\cal L}(5/8+it)|^2dt &=\int_0^T |{\cal L}_1(5/8+it){\cal L}_2(5/8+it)|^2dt\\
&\ll \left( \int_0^T |{\cal L}_1(5/8+it)|^{2 \cdot 4} dt\right)^{1/4}  \left( \int_0^T |{\cal L}_2(5/8+it)|^{2\cdot 4/3} dt\right)^{3/4}\\[1ex]
&\ll T^{1+\varepsilon}.
\end{align*}

Assume that \eqref{assump-d2} holds. Then by
\begin{equation*}
\int_0^T |{\cal L}_1(7/12+it)|^{6}dt\ll  T^{1+\varepsilon}
\end{equation*}
and
\begin{equation*}
\int_0^T |{\cal L}_2(7/12+it)|^{3}dt\ll  T^{1+\varepsilon}
\end{equation*}
(this comes from \eqref{MSQ2-version2-rei}), we have
\begin{equation*}
\int_0^T |{\cal L}(7/12+it)|^2dt \ll T^{1+\varepsilon}
\end{equation*}
similarly. This completes the proof of Corollary 2.

\section{Proofs of Theorem 2 and Theorem 3}


We shall follow Heath-Brown's argument \cite{He1} to prove Theorem 2 and Theorem 3.
We first quote some results  from \cite{He1}. The following Hypothesis (H), Lemma 13 and Lemma 14 are  Hypothesis (H),
Theorem 5 and Theorem 6 of \cite{He1}, respectively.

\medskip

\noindent {\bf Hypothesis (H)}: Let ${\cal M}(t)$ be a real valued function, $a_1(t),
a_2(t),\cdots,$ be continuous real valued functions of period 1,
and suppose there are non-zero constants $\gamma_1,\gamma_2,\cdots$ such that
$$
\lim_{N\rightarrow \infty}\limsup_{T\rightarrow\infty}
\frac{1}{T}\int_0^T\min\left(1,\left|{\cal M}(t)-\sum_{n\leq N}a_n(\gamma_nt)\right|\right)dt=0.
$$

\medskip

\begin{lem} \label{lem9-1}
Suppose ${\cal M}(t)$ satisfies {\rm (H)} and suppose that the constants $\gamma_i$ are linearly independent over ${\Bbb Q}$.
Suppose further
\begin{eqnarray*}
&&\int_0^1a_n(t)dt=0\ \ \ (n\in {\Bbb N}),\\
&&\sum_{n=1}^\infty\int_0^1a_n^2(t)dt<\infty,
\end{eqnarray*}
and there is a constant $\mu>1$ such that
\begin{eqnarray*}
&&\max_{t\in [0,1]}|a_n(t)|\ll n^{1-\mu},\\
&&\lim_{n\rightarrow\infty}n^\mu\int_0^1a_n^2(t)dt=\infty.
\end{eqnarray*}
Then ${\cal M}(t)$ has a distribution function $f(\alpha)$ with the properties described  in Theorem 2.
\end{lem}

\medskip

\begin{lem} \label{lem9-2}
Suppose ${\cal M}(t)$ satisfies {\rm (H)} and that
\begin{eqnarray*}
\int_0^T|{\cal M}(t)|^Kdt\ll T
\end{eqnarray*}
holds for some positive number $K$. Then for any real number $k\in [0,K),$ the limit
\begin{eqnarray*}
\lim_{T\rightarrow\infty}\frac 1T\int_0^T|{\cal M}(t)|^Kdt
\end{eqnarray*}
exists.
\end{lem}

\bigskip

Suppose $T\leq y\leq (1+\delta)T,$ define
\begin{align*}
{\cal M}(y)&=\kappa_0^{-1}(2\pi/h)^{-\frac{d-1}{2}}y^{-\frac{d-1}{2}}E((2\pi)^dh^{-d}y^d), \\
a_n(y)&= \frac{q_d(n)}{n^{\frac{d+1}{2d}}}\sum_{m=1}^\infty\frac{a(nm^d)}{m^{\frac{d+1}{2}}}\cos(2\pi ym+c_0\pi), \\
\intertext{and}
\gamma_n&=n^{1/d},
\end{align*}
where $q_d(n)$ is $1$ if $n$ is $d$-free and $0$ otherwise.

It is easy to see that $a_n(y)$ satisfies all conditions of Lemma \ref{lem9-1} for any fixed constant
$1+1/d<\mu<3/2+1/2d.$

Now we suppose $M_0\asymp T^{\frac{d}{2d(1-\sigma^{*})-1}}$ such that $E(M_0)\ll T^{d\theta+\varepsilon}$ and
$N_0=T^{d(2d-1)-\varepsilon}$. By \eqref{Eq8-4} we have
\begin{eqnarray*}
&&{\cal M}(y)={\cal M}_1(y)+{\cal M}_2(y),\\
&&{\cal M}_1(y)= \sum_{1\leq n\leq M_0}\frac{a(n)}{n^{\frac{d+1}{2d}}}\cos\left(2\pi yn^{1/d}
+c_0\pi\right),\\
&&{\cal M}_2(y)=\kappa_0^{-1}(2\pi/h)^{-\frac{d-1}{2}}y^{-\frac{d-1}{2}}\sum_{j=2}^7 R_j((2\pi)^dh^{-d}y^d),
\end{eqnarray*}
where $R_j(y)$ is defined in Subsection 8.1.

It is easy to see that for any integer $N_1\le M_0$ we have
\begin{align*}
&|{\cal M}(y)-\sum_{n\leq N_1}a_n(\gamma_n y)|\\
& \ll \left| \sum_{n\leq M_0}{\!\!}^{^ *}\, \frac{a(n)}{n^{\frac{d+1}{2d}}}\cos(2\pi yn^{1/d}+c_0\pi)\right|
      +\sum_{n\leq N_1}\frac{1}{n^{\frac{d+1}{2d}}}\sum_{m>(M_0n^{-1})^{1/d}}\frac{|a(nm^d)|}{m^{(d+1)/2}} +|{\cal M}_2(y)|\\
& \ll \left| \sum_{n\leq M_0}{\!\!}^{^\ast}\, \frac{a(n)}{n^{\frac{d+1}{2d}}}\cos(2\pi yn^{1/d}+c_0\pi)\right|
      +N_1^{1-1/d}M_0^{\theta+1/2d-1/2+\varepsilon}+|{\cal M}_2(y)|,
\end{align*}
where $  \sum^{*}$ means that $n$ has a $d$-free kernel great that $N_1,$ and for the sum involving $a(nm^d)$ in the above formula
 we used the first estimate of \eqref{meanC} and the condition $\theta<1/2-1/2d.$

Similarly to the approach of \eqref{R1} we have
\begin{align*}
& \int_T^{(1+\delta)T}\left| \sum_{n\leq M_0} {\!\!}^{^ *}\,\frac{a(n)}{n^{\frac{d+1}{2d}}}\cos(2\pi yn^{1/d}+c_0\pi)\right|^2dy\\
& \ll T \sum_{n\leq M_0} {\!\!}^{^ *}\,\frac{a^2(n)}{n^{\frac{d+1}{d}}}
      +\sum_{\stackrel{n\not= m}{n,m\le M_0}}\frac{|a(n)a(m)|}{(nm)^{\frac{d+1}{2d}}|n^{1/d}-m^{1/d}|}\\
& \ll T \sum_{n> N_1}\frac{a^2(n)}{n^{\frac{d+1}{d}}}+M_0^{1-2/d}\log^\eta M_0\\
& \ll TN_1^{-1/d+\varepsilon} +M_0^{1-2/d+\varepsilon}.
\end{align*}

Let ${\cal M}_2^{*}(y)=\sum_{j=2}^7 R_j(y), $ then ${\cal M}_2(y)={\cal M}_2^{*}((2\pi)^dh^{-d}y^d).$
From \eqref{R2}, \eqref{R3}--\eqref{R7} we get
$$
\int_T^{(1+\delta)T}{\cal M}_2^{*2}(y)dy\ll \sum_{j=2}^7\int_T^{(1+\delta)T}R_j^2(y)dy
\ll T^{1-\frac{d+1-2d\sigma^{*}}{2d(1-\sigma^{*})-1}+\varepsilon},
$$
which combining a change of variable implies that
\begin{eqnarray*}
 \int_T^{(1+\delta)T}{\cal M}_2^2(y)dy\ll T^{1/d-\frac{d+1-2d\sigma^{*}}{2d(1-\sigma^{*})-1}+\varepsilon}\ll T^{-c}
\end{eqnarray*}
for some positive constant $c>0.$

From the above estimates and Cauchy's inequality we get
\begin{eqnarray*}
\limsup_{T\rightarrow \infty}\frac 1T\int_T^{2T} \Bigl|{\cal M}(y)-\sum_{n\leq
N_1}a_n(\gamma_n y)\Bigr|^2dy\ll N_1^{-1/d+\varepsilon}
\end{eqnarray*}
and whence Hypothesis (H) follows. From Lemma 13 with $\mu=5/4+3/4d$ we
get Theorem 2.

From Theorem 1 and integration by parts we get easily that
$$\int_0^T\left(y^{-\frac{d-1}{2d}}E(y)\right)^2dy\ll T,$$
hence Theorem 3 follows from Lemma \ref{lem9-2} by taking ${\cal M}(y)=y^{-\frac{d-1}{2d}}E(y).$


\section{\bf Some applications of our main results}


In this section we consider some applications of our results.

\subsection{Some functions of degree $2$}

{\bf  Example 1:} The Dirichlet divisor problem.
The first example is $\zeta^2(s)\in {\cal S}^{0}_{real}.$
Since
$$
\zeta^2(s)=\sum_{n=1}^{\infty} d(n)n^{-s} \qquad  \Re s >1,
$$
where $d(n)$ is the divisor function, this is the Dirichlet divisor problem. This problem is fundamental for
all theory and hence has a long history. Let us review some of  them.
 Dirichlet first proved that the error term
$$
\Delta(x): =\sideset{}{^{\prime}}\sum_{n\leq x}d(n)-x\log x-(2\gamma-1)x,\hspace{3mm}x\geq 2
$$
satisfies
$\Delta(x)=O(x^{1/2}).$
The latest result is due to Huxley \cite{Hux2}, who showed that
$$
\Delta(x)\ll x^{131/416}(\log x)^{26947/8320}.
$$
Cram\'er \cite{Cr} proved the classical mean-square result
$$
\int_1^T\Delta^2(x)d x \;=\;\frac{(\zeta(3/2))^4}{6\pi^2 \zeta(3)}T^{3/2}+F(T)
$$
with $F(T)=O(T^{5/4+\varepsilon}).$
The estimate $F(T)=O(T^{5/4+\varepsilon})$ was improved to $O(T\log^5 T)$ by Tong \cite{To2},
from which the present work originates. Tong's result was improved to $O(T\log^4 T)$ by
Preissmann \cite{P}, and recently to $O(T\log^3 T \log\log T)$ by
Lau and Tsang \cite{LT2}. For the asymptotic formulae of higher power moments of $\Delta(x)$
see, for example, the papers of Ivi\'c and Sargos \cite{IS},
Tsang \cite{Ts1} and Zhai  \cite{Z}.
For a survey of  the Dirichlet divisor problem, see for example, Kr\"atzel \cite{Kra}, Ivi\'c \cite{I1} or Tsang \cite{TS}.

Corollary 1 gives $F(T)=O(T^{1+\varepsilon})$. Since our setting of the problem is of general nature, it is weaker than the
estimate of the form $F(T)\ll T\log^{c} T,$ but is much stronger than Cram\'er's result $F(T)=O(T^{5/4+\varepsilon}).$

\medskip

{\bf Example 2:} The automorphic $L$-function attached to holomorphic cusp forms over $SL_2({\Bbb Z})$.
Let $f$ be a primitive holomorphic cusp form  of weight $k \ge 1$ for the full modular group
$SL_2({\Bbb Z}).$  (For the sake of simplicity we only consider the form for $SL_2(\mathbb{Z})$.)
Let
\begin{equation}
f(z) =\sum_{n=1}^{\infty} \lambda_f(n)n^{(k-1)/2}e(nz), \label{eq10-1}
\end{equation}
where $e(z)=e^{2\pi i z}$, be its normalized Fourier expansion at the cusp $\infty$. Then the automorphic $L$-function
$$
{\cal L}(f,s)=\sum_{n=1}^{\infty} \lambda_f(n)n^{-s}
$$
is an $L$-function of degree $2$ satisfying the functional equation
$$
(2\pi)^{-s}\Delta(s){\cal L}(f,s)=(-1)^{k/2}(2\pi)^{-(1-s)}\Delta(1-s){\cal L}(f,1-s)
$$
with the gamma factor $\Delta(s)=\Gamma\left(s+\frac{k-1}{2}\right)$.
Deligne \cite{D} proved that $|\lambda_f(n)|\le d(n)$, so we have ${\cal L}(f,s)\in {\cal S}_{real}^0$.
It is well-known that $$A(y)=\sum_{n\le y}\lambda_f(n)\ll y^{1/3} $$ and
(see Rankin \cite{Ran} or Selberg \cite{Se})
$$\sum_{n\le y}|\lambda_f(n)|^2=c_f y+O(y^{3/5}) $$ for some positive constant $c_f.$
Walfisz \cite{Wa} proved that
\begin{equation*}
\int_1^TA^2(y)dy=C_AT^{3/2}+O(T\log^2 T)
\end{equation*}
holds for some positive constant $C_A>0$.
From Corollary 1  we have a weaker result
\begin{equation*}
\int_1^TA^2(y)dy=C_AT^{3/2}+O(T^{1+\varepsilon}).
\end{equation*}

\medskip

{\bf  Example 3:} The automorphic $L$-function attached to Maass forms for $SL_2({\Bbb Z})$. \\
Let $\varphi$ be a primitive Maass form for $SL_2({\Bbb Z})$, which
is an eigenfunction of the Laplace operator with an eigenvalue
$\lambda= 1/4 +r^2,$ where $r \in \mathbb{R}.$  Write its Fourier expansion at infinity in the form
$$
\varphi(z)=\sqrt{v}\sum_{n\in {\Bbb Z}\setminus \{0\}}\rho(n)K_{ir}(2\pi |n|v)e(nu) \quad (\ z=u+iv, u\in {\Bbb R},v>0),
$$
where $K_{ir}$ is the modified Bessel function of the third kind. The corresponding automorphic $L$-function is
$$
{\cal L}(\varphi,s)=\sum_{n=1}^{\infty} \rho(n)n^{-s},
$$
which is a function of degree $2$, entire on ${\Bbb C},$ with the functional equation
$$
\pi^{-s}\Delta(s){\cal L}(\varphi,s)=(-1)^{\delta}\pi^{-(1-s)}\Delta(1-s){\cal L}(\varphi,1-s),
$$
where $\Delta(s)=\Gamma\left((s+ \delta+ir)/2\right)\Gamma\left((s+ \delta-ir)/2\right),$
and $\delta$ is the parity of $\varphi$ defined by $\delta=0$ if $ \varphi$ is even and $\delta=1$ if $ \varphi$ is odd.
Kim and Sarnak \cite{KS} proved that
\begin{equation}
|\rho(n)|\le d(n)n^{7/64}. \label{eq10-2}
\end{equation}
 The Rankin-Selberg method implies that
\begin{equation}
\sum_{n\le y}|\rho(n)|^2\ll y .  \label{eq10-3}
\end{equation}
So we have ${\cal L}(\varphi,s)\in {\cal S}_{real}^{ \theta} $ with $\theta=7/64.$

Hafner and Ivi\'c \cite{HI} proved that
$$
A_r(y)=\sum_{n\le y}\rho(n)\ll y^{2/5}.
$$
Meurman \cite{Me2} proved that
$$
 A_r(y)\ll
y^{\frac{1+\mu}{3}+\varepsilon},
$$
where $\mu\ge 0$ is any real number such that $\rho(n)\ll n^\mu.$

From Corollary 1  we get
\begin{equation*}
\int_1^TA_r^2(y)dy=C_{A_r}T^{3/2}+O(T^{1+\varepsilon}),
\end{equation*}
where $C_{A_r}>0$ is a positive constant.


\subsection{Some functions of degree $3$}


{\bf Example 1:} Let ${\cal D}_{r}$ denote the set of all Dirichlet $L$-functions corresponding to real primitive Dirichlet
characters,   ${\cal D}_{c}$ denote the set of all Dirichlet $L$-functions corresponding to complex primitive Dirichlet
characters,   $f$ denote a primitive holomorphic cusp form and $\varphi$ denote a primitive Maass form for $SL_2({\Bbb Z})$
for some fixed  eigenvalue $\lambda= 1/4 +r^2 (r \in \mathbb{R}),$ respectively. Define
\begin{eqnarray*}
&&{\cal D}_{1f}=\{{\cal L}(f,s)L(s,\chi):  \ L(s,\chi)\in {\cal D}_r\},\\
&&{\cal D}_{1\varphi}=\{{\cal L}(\varphi,s)L(s,\chi):\ L(s,\chi)\in {\cal D}_r\},\\
&&{\cal D}_r^3=\{L(s,\chi_1)L(s,\chi_2)L(s,\chi_3):\ L(s,\chi_j) \in {\cal D}_r,j=1,2,3  \},\\
&&{\cal D}_r{\cal D}_c\overline{{\cal D}_c}=\{L(s,\chi_1)L(s,\chi_2)L(s, \overline{\chi_2}):\ L(s,\chi_1) \in {\cal D}_r, L(s,\chi_2) \in {\cal D}_c  \}.
\end{eqnarray*}

Let ${\cal L}(s)$ denote any function in the set ${\cal D}_{1f}\cup {\cal D}_{1\varphi} \cup {\cal D}^3\cup {\cal D}_r{\cal D}_c \overline{{\cal D}_c}.$
It is easy to see that ${\cal L}(s)$ is a function of degree $3$, which satisfies the conditions of Corollary 2.
 We consider only the case ${\cal L}(s)=L(s,\chi){\cal L}(\varphi,s), $ where
 $\chi$ is a real primitive Dirichlet character modulo some positive integer.  Write
$$
{\cal L}(s)=\sum_{n\ge 1}a(n)n^{-s},\ \Re s>1.
$$
Then $a(n)=\sum_{d|n}\rho(d)\chi(n/d).$ From the estimate \eqref{eq10-2} we have
$a(n)\ll d^2(n)n^{7/64}.$ By Cauchy's inequality we have
$$
a^2(n)\le d(n)\sum_{d|n}\rho^2(d),
$$ which combining (\ref{eq10-3}) gives
\begin{align*}
\sum_{n\leq y}a^2(n)
&\ll y^{\varepsilon/2}\sum_{dn\le y}\rho^2(d) \ll y^{\varepsilon/2}\sum_{n\le y}\sum_{d\le y/n}\rho^2(d)
\ll y^{1+\varepsilon},
\end{align*}
where we used the bound $d(n)\ll n^{\varepsilon/2}.$ Whence \eqref{meanC} holds.

As an analogue of the sixth moment of the automorphic $L$-function attached to a cusp form \cite{J2}, we have
also $$\int_1^T |{\cal L}(\varphi,1/2+it)|^6dt\ll T^{2+\varepsilon}.$$
This estimate can be  proved by Jutila's method analogously. So from Lemma \ref{lemma_degree2} we have
$$\int_1^T |{\cal L}(7/12+it)|^2dt\ll T^{1+\varepsilon},$$
namely, \eqref{meanL} holds for $\sigma^{*}=7/12,$ which obviously satisfies \eqref{sigmastarcondition}.   So from
 Corollary 2 we have
 $$\int_1^TE^2(y)dy=C_ET^{5/3}+O(T^{14/9+\varepsilon}),$$
where $C_E$ is a positive constant.

\medskip

{\bf Example 2:} Let $K/{\Bbb Q}$ be a number field of degree $\kappa$ with discriminant $d_K$. The Dedekind  zeta-function defined by
$$
\zeta_K(s)=\sum_{\mathfrak{a}\in \mathfrak{L}\setminus\{0\}}\frac{1}{(N\mathfrak{a})^s}
$$
is an $L$-function of degree $\kappa$ with the functional equation
$$
(\pi^{-\kappa}|d_K|)^{s/2}\Delta(s)\zeta_K(s)=(\pi^{-\kappa}|d_K|)^{(1-s)/2}\Delta(1-s)\zeta_K(1-s)
$$
with the gamma factor $\Delta(s)=\Gamma^{r_1+r_2}\left(\frac{s}{2}\right)\Gamma^{r_2}\left(\frac{s+1}{2}\right)$,
where $r_1$ is the number of real embeddings of $K$ and $r_2$ the number of pairs of
complex embeddings so that $ \kappa = r_1 + 2r_2.$ The Dedekind  zeta-function $\zeta_K(s)$ has a simple at $s=1.$

From now we suppose that  $\kappa=3.$ Let
$E_K(y): =\sum_{n\le y}a(n)-\lambda y,$
where $a(n)$ denotes the number of ideals in $\mathfrak{L}$ of norm $n$ and $\lambda $ is the residue
of $\zeta_K(s)$ at $s=1.$
 Weber and Landau showed that
$$
E_K(y)=O\left(y^{1/2}\right).$$

In 1964, Chandrasekharan and Narasimhan \cite{CN1} studied the mean square of $E_K(y)$
 and proved
$$\int_1^T E_K^2(y)dy=
 O\left(T^{5/3}\log^{3}T\right). $$

In 1999, Y.-K. Lau \cite{La} improved the above result   to
 $$\int_1^T E_K^2(y)dy= O\left(T^{5/3}\log^{2}T\right).$$

Fomenko \cite{Fo1} followed Tong's approach and proved that the asymptotic formula
 \begin{align*}
\int_1^TE_{K}^2(y)dy= C_{K}T^{5/3}+O\left(T^{8/5+\varepsilon}\right)
\end{align*}
holds for some positive constants $C_{K}>0$ with the help of the estimate
\begin{equation}
\int_1^T |\zeta_K(\sigma+it)|^2\ll T^{1+\varepsilon} \label{eq10-5}
\end{equation}
with $\sigma=5/8.$

When $\kappa=3$, we can show  that the Dedekind zeta-function $ \zeta_K(s)$ can be written as
$ \zeta_K(s)=\zeta(s){\cal L}_2(s) $ such that where ${\cal L}_2(s) \in {\cal S}_{real}^{0} $ is a function of degree $2 $
 satisfying the condition \eqref{assump-d2}.
 When $K/{\Bbb Q}$ is a normal extension, from Lemma 1 of  W. M\"uller \cite{Mu}
we know that
$\zeta_K(s)=\zeta(s)L(s,\chi_1)L(s,\overline{\chi_1})$ for some primitive Dirichlet character $\chi_1.$
So  we have $ \zeta_K(s)=\zeta(s){\cal L}_2(s),$ where
 ${\cal L}_2(s)=L(s,\chi_1)L(s,\overline{\chi_1})$ satisfies \eqref{assump-d2} via Meurman  \cite{Me3}.

 When  $K/{\Bbb Q}$ is a not normal extension, from the formula (1) of Fomenko \cite{Fo1} we have
 $\zeta_K(s)=\zeta(s)L(s,F),$ where
$F$ is a  holomorphic primitive cusp form of weight $1$ with respect to  $\Gamma_0(|d_K|) $ and
$L(s,F)$ is its corresponding automorphic $L$-function of degree $2.$ We have the estimate
$$\int_1^T|L(1/2+it, F)|^6dt\ll T^{2+\varepsilon}.$$
The proof is also due to Jutila \cite{J2}, who proved it in the case of a holomorphic cusp form
$F$ of even weight with respect to $SL_2({\Bbb Z}).$ His proof can be applied to $L(s,F),$ too.

From Lemma \ref{lemma_degree2} we  see that  (\ref{eq10-5}) holds for $\sigma=7/12.$ So
from Corollary  2 we get
 \begin{align*}
\int_1^TE_{K}^2(y)dy= C_{K}T^{5/3}+O\left(T^{14/9+\varepsilon}\right),
\end{align*}
which is a slight improvement to Fomenko's result.


\medskip


\subsection{Some examples of degree $4$}


{\bf Example 1.}
Let ${\cal D}_{r}$ denote the set of all Dirichlet $L$-functions corresponding to real primitive Dirichlet
characters,   ${\cal D}_{c}$ denote the set of all Dirichlet $L$-functions corresponding to complex primitive Dirichlet
characters,
  $f$ denote a primitive holomorphic cusp form and $\varphi$ denote a primitive Maass form for $SL_2({\Bbb Z})$
for some fixed  an eigenvalue $\lambda= 1/4 +r^2 (r \in \mathbb{R}),$ respectively. Define
\begin{eqnarray*}
&&{\cal D}_{2f}=\{{\cal L}(f,s)L(s,\chi_1)L(s,\chi_2):  \ L(s,\chi_j) \in {\cal D}_r, \quad j=1,2 \},\\
&&{\cal D}_{2 \varphi}=\{{\cal L}(\varphi,s)L(s,\chi_1)L(s,\chi_2):  \ L(s,\chi_j) \in {\cal D}_r, \quad j=1,2\} ,\\
&&{\cal D}_r^4=\{\prod_{j=1}^4 L(s,\chi_j) :\ L(s,\chi_j) \in {\cal D}_r, 1\le j\le 4\},\\
&&{\cal D}_r^2{\cal D}_c\overline{{\cal D}_c}=\{L(s,\chi_1)L(s,\chi_2)L(s,\chi_3)L(s,\overline{\chi_3}) :
\ L(s,\chi_1), L(s,\chi_2) \in {\cal D}_r, L(s,\chi_3) \in {\cal D}_c \},\\
&&({\cal D}_c\overline{{\cal D}_c})^2=\{\prod_{j=1}^2 L(s,\chi_j)L(s,\overline{\chi_j}) :
\   L(s,\chi_j) \in {\cal D}_c, j=1,2 \}.
\end{eqnarray*}

Suppose ${\cal L}(s)$ is any function in the set
$$ {\cal D}_{2f}\cup{\cal D}_{2 \varphi}\cup {\cal D}^4 \cup {\cal D}_r^2{\cal D}_c\overline{{\cal D}_c}  \cup ({\cal D}_c\overline{{\cal D}_c})^2
 \cup\{  {\cal L}^2(f,s),
{\cal L}^2( \varphi,s), {\cal L}(f,s){\cal L}( \varphi,s) \}.$$
Obviously ${\cal L}(s)$  is a function of  degree $4.$ From Lemma \ref{lemma_degree2} it is easy to see  that
\begin{equation} \label{eq10-7}
\int_1^T|{\cal L}(5/8+it)|^2dt \ll T^{1+\varepsilon}.
\end{equation}

For the mean square of $E(y)$ corresponding to ${\cal L}(s)$ we have
\begin{equation}
T^{7/4}\ll \int_1^T E^2(y)dy\ll T^{7/4+\varepsilon}. \label{eq10-8}
\end{equation}

The left inequality in \eqref{eq10-8} follows from Theorem 5 directly. The proof of   the right inequality in \eqref{eq10-8}
is the same as the case $k=4$ in  Theorem 13.4 of Ivi\'c \cite{I1} with the help of \eqref{eq10-7}.
So we omit the details. We note that the upper bound in \eqref{eq10-8} is already proved in Ivi\'c \cite{I8}
when ${\cal L}(s)={\cal L}^2(f,s).$

 If we could find a number $1/2\le \sigma_0<5/8$ such that
  $$\int_1^T|{\cal L}(\sigma_0+it)|^2dt \ll T^{1+\varepsilon},$$
then from Theorem 1 the asymptotic formula
\begin{equation}\int_1^TE^2(y)dy=C_4T^{7/4}+O(T^{7/4-\delta_{\cal L}})\end{equation}
would hold, where  $C_4>0$ and $ \delta_{\cal L}>0$   are positive constants. The formula (10.7) implies that the estimate
(10.6) is best possible.

\medskip

{\bf Example 2:}   The error term in the Rankin-Selberg problem.

Let $f$ be a primitive holomorphic cusp form  of weight $k$ for  the full modular group
$SL_2({\Bbb Z})$ and its  normalized Fourier expansion at the cusp $\infty$ can be written as \eqref{eq10-1}.
The corresponding Rankin-Selberg zeta function defined by
$$
{\cal L}(f\otimes f, s)=\sum_{n=1}^{\infty} c(n)n^{-s},\ c(n)=\sum_{m^2|n}|\lambda_f(n/m^2)|^2,\ (\Re s>1)
$$
is a function of degree $4$ satisfying the functional equation
$$
(2\pi)^{-2s}\Delta(s){\cal L}(f\otimes f,s)=(2\pi)^{-2(1-s)}\Delta(1-s){\cal L}(f \otimes f,1-s)
$$
with the gamma factor $\Delta(s)=\Gamma(s)\Gamma(s+k-1)$.

Rankin \cite{Ran} and Selberg \cite{Se} proved independently that
\begin{align*}
 \sum_{n\le y}c(n)=C_fy+E(y), \  E(y)=O(y^{3/5}),
\end{align*}
where $C_f>0$ is a positive constant.  From Deligne's estimate $|\lambda_f(n)|\le d(n)$ we have
\begin{equation}
c(n)\le d^2(n)d(1,2;n)\le d^3(n)\ll n^{\varepsilon}, \label{eq10-11}
\end{equation}
where $d(1,2;n)=\sum_{n=m^2d}1\le d(n).$    We write
$c(n)=\sum_{n=m^2d}|\lambda_f(d)|^2.$ By Cauchy's inequality we have
$$
c^2(n)\le d(1,2;n)\sum_{n=m^2d}|\lambda_f(d)|^4\ll n^{\varepsilon}\sum_{n=m^2d}|\lambda_f(d)|^4,
$$
where we used $d(1,2;n) \ll n^\varepsilon.$ So we have
\begin{align}\label{eq10-12}
\sum_{n\le y}c^2(n)
&\ll y^\varepsilon\sum_{m^2d\le y}|\lambda_f(d)|^4  \\
&\ll y^\varepsilon\sum_{m\le y^{1/2}}\sum_{d\le y/m^2}|\lambda_f(d)|^4\nonumber\\
&\ll y^\varepsilon\sum_{m\le y^{1/2}}  y/m^2 \ll y^{1+\varepsilon},   \nonumber
\end{align}
where we used the estimate in L\"u \cite{Lv}
$$
\sum_{n\le y}|\lambda_f(n)|^4\ll y.
$$
The estimates (\ref{eq10-11}) and (\ref{eq10-12}) show  that \eqref{meanC} holds.

Ivi\'c \cite{I3} proved that
\begin{align*}
\int_{1}^{T}E^2(y)dy\ll T^{  1+2\beta+\varepsilon},
\end{align*}
where $\beta=2/(5-2\mu(1/2))$ with $\mu(1/2)$ satisfying $\zeta(1/2+it)\ll (|t|+1)^{\mu(1/2)+\varepsilon}.$ Huxley's best bound\cite{Hux3}
$\mu(1/2)\leq 32/205$ implies $\beta\leq 410/961=0.4266\cdots.$  The Lindel\"of Hypothesis implies $\beta\leq 0.4.$
In \cite{I2}  Ivi\'c even conjectured
\begin{align}
\int_{1}^{T}E^2(y)dy\ll T^{\frac 74+\varepsilon}.\label{eq10-13}
\end{align}
The estimate (\ref{eq10-13}), if it is true, is best possible. Since from Theorem 5 we have
\begin{align*}
\int_1^TE^2(y)dy\gg T^{\frac 74}.
\end{align*}

If we could find a $\sigma_0<\frac 58$ such that
$$
\int_1^T |{\cal L}(f\otimes f, \sigma_0+it)|^2dt\ll T^{1+\varepsilon},
$$
then from  Theorem 1 we would have the asymptotic formula
\begin{align*}
\int_1^TE^2(y)dy= C_fT^{\frac 74}+O\left(T^{\frac 74-\delta_f}\right)
\end{align*}
for some positive constants $C_f>0, \delta_f>0.$

\begin{flushleft}
Xiaodong Cao \\
Department of Mathematics and Physics, \\
Beijing Institute of Petro-Chemical Technology,\\
Beijing, 102617, P. R. China \\
e-mail: caoxiaodong@bipt.edu.cn

\bigskip

Yoshio Tanigawa\\
Graduate School of Mathematics,\\
Nagoya University, \\
Nagoya, 464-8602, Japan\\
e-mail: tanigawa@math.nagoya-u.ac.jp

\bigskip

Wenguang Zhai\\
Department of Mathematics, \\
China University of Mining and Technology, \\
Beijing 100083, P. R. China\\
e-mail: zhaiwg@hotmail.com
\end{flushleft}

\end{document}